\documentclass{amsart}

\DeclareMathOperator{\Lie}{Lie}
\DeclareMathOperator{\val}{val}

\DeclareMathOperator{\dom}{dom}
\DeclareMathOperator{\hgt}{ht}
\DeclareMathOperator{\id}{id}
\DeclareMathOperator{\Ad}{Ad}
\DeclareMathOperator{\defect}{def}
\DeclareMathOperator{\im}{im}
\DeclareMathOperator{\ad}{ad}
\DeclareMathOperator{\Int}{Int}

\newcommand{\Sdom}{\Sigma(\mu)_{M\text{-}\dom}}
\newcommand{\Smax}{\Sigma(\mu)_{M\text{-}\max}}
\newcommand{\Z}{\mathbb Z}
\newcommand\gfrac[2]{\genfrac{}{}{0pt}{}{#1}{#2}}

\usepackage{amscd}
\usepackage{amssymb}
\usepackage[all]{xy}
\usepackage{graphicx}

\numberwithin{equation}{subsection}

\newtheorem{theorem}{Theorem}[subsection]
\newtheorem{corollary}[theorem]{Corollary}
\newtheorem{lemma}[theorem]{Lemma}
\newtheorem{proposition}[theorem]{Proposition}

\theoremstyle{definition}
\newtheorem{defn}[theorem]{Definition}

\newtheorem{conjecture}[theorem]{Conjecture}

\begin{document}
\title{Dimensions of some affine Deligne-Lusztig varieties}

\author[U. G\"{o}rtz]{Ulrich G\"{o}rtz}
\address{Ulrich G\"{o}rtz\\Mathematisches Institut\\Universit\"{a}t Bonn\\Beringstr. 1\\53115
Bonn\\Germany} \email{ugoertz@math.uni-bonn.de}

\author[T. J. Haines]{Thomas J. Haines}
\address{Thomas J. Haines\\Mathematics Department\\ University of Maryland\\ College Park, MD
20742-4015} \email{tjh@math.umd.edu}
\thanks{Haines was partially supported by NSF Grant DMS-0303605 and a Sloan Research
Fellowship}

\author[R. E. Kottwitz]{Robert E. Kottwitz}
\address{Robert E. Kottwitz\\Department of Mathematics\\ University of Chicago\\ 5734 University
Avenue\\ Chicago, Illinois 60637}
\email{kottwitz@math.uchicago.edu}
\thanks{Kottwitz was partially supported by NSF Grant DMS-0245639}

\author[D. C. Reuman]{Daniel C. Reuman}
\address{Daniel C. Reuman\\Laboratory of Populations\\ Rockefeller University\\ 1230
York Ave.\\ New York, NY 10021} \email{reumand@rockefeller.edu}
\thanks{Reuman was supported by NSF grants DEB-9981552 and
DMS-0443803}


\maketitle

\section{Introduction}

 Let $k$ be a finite field with $q$ elements,
and let $\bar k$ be an algebraic
closure of $k$.  We consider the field $L:=\bar k((\epsilon))$
and its subfield
$F:=k((\epsilon))$. We write $\sigma:x\mapsto x^q$ for the Frobenius
automorphism of
$\bar k/k$, and we also regard $\sigma$ as an automorphism of $L/F$ in the
usual way, so that $\sigma(\sum a_n\epsilon^n)=\sum \sigma(a_n)\epsilon^n$.
We write
$\mathfrak o$  for the valuation ring $\bar k[[\epsilon]]$ of~$L$.

Let $G$ be a split connected reductive group over $k$, and let $A$ be a split
maximal torus of~$G$. Put $\mathfrak a:=X_*(A)_ \mathbb R$. Write
$W$ for the Weyl group of
$A$ in
$G$.   Fix a Borel subgroup
$B=AU$ containing $A$ with unipotent radical~$U$.
For $\lambda \in X_*(A)$ we
write
$\epsilon^\lambda$ for the element of $A(F)$ obtained as the image of
$\epsilon
\in \mathbb G_m(F)$ under the homomorphism $\lambda:\mathbb G_m \to A$.

This paper concerns the dimensions of certain affine Deligne-Lusztig
varieties, both in the affine Grassmannian and in the affine flag manifold.
We begin with the affine Grassmannian.

Put
$K:=G(\mathfrak o)$.  We denote by $X$ the affine Grassmannian $X=G(L)/K$ and
by $x_0$ its obvious base-point. The group $G(L)$ acts by left translation on
$X$.
By the Cartan decomposition $G(L)$ is the disjoint union of the subsets
$K\epsilon^\mu K$, with $\mu$ running over the dominant elements in
$X_*(A)$. For
$b \in G(L)$ and a dominant coweight $\mu \in X_*(A)$ the affine
Deligne-Lusztig variety
$X_\mu(b)=X^G_\mu(b)$ is the locally closed subset of
$X$ defined by
\begin{equation}
X_\mu(b):=\{x \in G(L)/K:x^{-1}b\sigma(x) \in K\epsilon^\mu K\}.
\end{equation}
For $g \in G(L)$ it is clear that $x\mapsto gx$ yields an isomorphism
$X_\mu(b)
\to X_\mu(g b\sigma(g)^{-1})$, so the isomorphism class of $X_\mu(b)$ depends
only on the
$\sigma$-conjugacy class of $b$.

Let $\mathbb D$ be the diagonalizable group over $F$ with character group
$\mathbb Q$. Just as in \cite{kottwitz85}, the element $b$ determines a
homomorphism
$\nu_b:\mathbb D \to G$ over $L$, and $b$ is said to be \emph{basic} if
$\nu_b$ factors through the center of $G$. For $g \in G(L)$ we have
\begin{equation}\label{eq.sigconj}
\nu_{gb\sigma(g)^{-1}}=\Int(g) \circ \nu_b,
\end{equation}
where $\Int(g)$ denotes the
inner automorphism $x \mapsto gxg^{-1}$ of $G$ over $L$, and since it is
harmless to $\sigma$-conjugate $b$, we may as well assume that $\nu_b$
factors through $A$, and that the corresponding element $\bar \nu_b \in
X_*(A)_\mathbb Q \subset \mathfrak a$ is dominant.
Following \cite{rapoport-richartz96} we refer to $\bar\nu_b$ as the
\emph{Newton point} of $b$. The centralizer $M_b$ of
$\nu_b$ in $G$ is then a Levi subgroup of $G$ over $F$.

Just as in \cite{kottwitz85}, \cite[1.12]{rapoport-zink96},
\cite{kottwitz97}, there is an inner form $J$ of
$M_b$ whose $R$-valued points (for any $F$-algebra $R$) are given by
\[
J(R)=\{g \in G(R \otimes_F L): g^{-1}b\sigma(g)=b \}.
\]
The group $J(F)$ acts by left multiplication on $X_\mu(b)$. Note that $J(F)
\subset M_b(L)$ because of \eqref{eq.sigconj}.

Let  $\Lambda_G$ denote the quotient of $X_*(A)$ by the coroot lattice for
$G$. We denote by $p_G$ the canonical surjection $X_*(A)
\twoheadrightarrow \Lambda_G$. There is a canonical homomorphism
$\eta_G:G(L) \twoheadrightarrow \Lambda_G$, which is trivial on
$K=G(\mathfrak o)$ and hence induces a surjection, also denoted
$\eta_G$, from $X=G(L)/K$ to $\Lambda_G$. The fibers of $\eta_G:X
\twoheadrightarrow \Lambda_G$ are the connected components of the affine
Grassmannian $X$.
As in \cite[7.7]{kottwitz97}, the restriction of the homomorphism
$\eta_{M_b}:M_b(L)
\twoheadrightarrow \Lambda_{M_b}$ to $J(F)$ is surjective, which implies that
the restriction of $\eta_G:G(L)\twoheadrightarrow \Lambda_G$ to $J(F)$ is
also surjective. Using the action of $J(F)$ on $X_\mu(b)$, we then see that
the intersections $X_\mu(b) \cap \eta_G^{-1}(\lambda)$ ($\lambda \in
\Lambda_G$) of $X_\mu(b)$ with the various connected components of $X$ are
all isomorphic to each other.

As in \cite{kottwitz-rapoport02} (see also \cite{leigh02},
\cite{kottwitz03}), there is a simple criterion for $X_\mu(b)$ to be
non-empty (see Proposition \ref{prop.red.basic}). Here we mention only that
when $b$ is basic, $X_\mu(b)$ is non-empty if and only if
$\eta_G(b)=p_G(\mu)$. When
$X_\mu(b)$ is non-empty, there is a conjectural formula for its dimension,
due to Rapoport \cite{rapoport02}, which we now recall, reformulating
it slightly (see \cite{kottwitz05})  along the way.

This formula involves a non-negative integer $\defect_G(b)$ attached to
$b$.  By
definition
$\defect_G(b)$ is the
$F$-rank of $G$ minus the $F$-rank of $J$.
Clearly $\defect_G(b)$ depends only on the $\sigma$-conjugacy class of
$b$.
(As usual the $F$-rank of $G$ is
the common dimension of all maximal $F$-split tori in $G$.) We write $\rho
\in X^*(A)_{\mathbb Q}$ for the half-sum of the positive roots.

\begin{conjecture}[Rapoport] Assume that $X_\mu(b)$ is non-empty. Then its
dimension is given by
\[
\dim X_\mu(b)=\langle \rho, \mu-\bar \nu_b \rangle-\frac{1}{2}\defect_G(b).
\]
\end{conjecture}
In \cite{Reu04} Reuman proves the conjecture for $G=SL_2, SL_3, Sp_4$ and
$b=1$. In \cite{mier} Mierendorff proves (in the context of $\mathbb Q_p$
rather than $k((\epsilon))$) that Rapoport's conjecture is true for $GL_n$
and minuscule
$\mu$.

In this paper we prove Rapoport's conjecture for all $b \in A(L)$ (see
Theorem \ref{main_result_thm}, noting that $\defect_G(b)=0$ when $b \in
A(L)$), and in fact we show in this case that $X_\mu(b)$ is equidimensional
(see Proposition \ref{equidim.prop}), answering a question of Rapoport.
Moreover, returning to general elements $b
\in G(L)$, we show in Theorem \ref{thm581} that, if $M$ is a Levi subgroup of
$G$ such that
$b \in M(L)$ and $b$ is basic in $M(L)$, and if Rapoport's conjecture holds
for $(M,b)$, then Rapoport's conjecture holds for $(G,b)$.

Consequently, in order to prove Rapoport's conjecture in general, it would
be enough to prove it for \emph{superbasic} elements $b$, by which we mean
those for which no $\sigma$-conjugate is contained in a proper Levi
subgroup of
$G$. As we verify in
\ref{subsec.superbasic}, superbasic elements are very special, and for
simple groups they exist only in type $A_n$. The upshot is that it would be
enough to prove Rapoport's conjecture for basic elements $b \in GL_n(L)$
such that the valuation of $\det(b)$ is relatively prime to $n$.

Now we turn to affine Deligne-Lusztig varieties inside the affine flag
manifold $G(L)/I$, where $I$ is the Iwahori subgroup of $G(L)$ obtained from
an alcove $\mathbf a_1$ in the apartment associated to $A$. Given $b \in
G(L)$ and an element $x \in \tilde W=W\ltimes X_*(A)$, we get the affine
Deligne-Lusztig variety
\[
X_x(b)=\{g \in G(L)/I:g^{-1}b\sigma(g) \in IxI \}.
\]
For the groups $G=SL_2, SL_3, Sp_4$ and
$b=1$ the dimension of $X_x(b)$ was computed by Reuman \cite{Reu04}. For
suitably general $x$ (those in the union of the ``shrunken Weyl chambers'')
Reuman gives a simple formula for $\dim X_x(b)$ and conjectures that it
holds in general (for $b=1$ and suitably general $x$).

In this article we prove a formula (see Theorem \ref{thm.631}) for $\dim
X_x(b)$ when $b \in A(L)$. Unfortunately this formula does not suffice to
establish Reuman's conjecture, since it involves the unknown dimension of the
intersection of $I$- and $U(L)$-orbits in the affine flag manifold. (The
reason we had better luck with the affine Grassmannian is that the
dimensions of intersections of $G(\mathfrak o)$- and $U(L)$-orbits in the
affine Grassmannian are known, thanks to Mirkovi\'c-Vilonen \cite{MV1}.)

However there is an algorithm for computing the dimensions of such
intersections, and in section \ref{sec.computations} we describe the results
of computer calculations made using this algorithm. Reuman's conjecture (see
subsection \ref{sub.comp.2}) turns out to hold in all cases checked by the
computer. In the case of rank 2 groups, the results can be presented in the
form of pictures. Figures 1,2,3 show the dimensions for $b=1$ and $A_2$,
$C_2$, $G_2$ respectively, while figures 11,12 show dimensions for two
elements $b \in A(L)$, $b \ne 1$, one for type $A_2$, and one for type $C_2$.

The results in these last two figures support Conjecture \ref{conj.bne1}, an
extension of Reuman's conjecture to elements $b \in A(L)$. We finish this
introduction by mentioning that there is mounting evidence that an analog of
Reuman's conjecture holds for all $b \in G(L)$ (for suitably general $x \in
\tilde W$).

\section{Affine Deligne-Lusztig varieties inside the affine Grassmannian}
\subsection{Further preliminaries}\label{subsec.notG}
 For  $\nu$ in $X_*(A)$ or
$\mathfrak a$ we write
$\nu_{\dom}$ for the unique dominant element in the $W$-orbit of $\nu$. For
cocharacters
$\mu,\nu \in X_*(A)$ we say that $\nu \le \mu$ if $\mu-\nu$ is a non-negative
integral linear combination of positive coroots.

Any $b \in A(L)$ is $\sigma$-conjugate to an element of the form
$\epsilon^\nu$.    By \cite{kottwitz-rapoport02} the set
$X_\mu(\epsilon^\nu)$ is non-empty if and only if $\nu_{\dom} \le \mu$; we now
assume that this is the case.  We are going to calculate the dimension of
$X_\mu(\epsilon^\nu)$, using the obvious fact that  $X_\mu(\epsilon^\nu)$ is
preserved by the action of
$A(F)$.

\subsection{Topology on $X$} We view $X$ as an ind-scheme in the usual way.
Each $K$-orbit on
$X$ is finite dimensional, and we denote by $Z_n$ the union of all $K$-orbits
having dimension less than or equal to $n$. Each $Z_n$ is a projective variety,
 and the  increasing family $Z_0 \subset Z_1 \subset Z_2 \subset Z_3 \subset
\dots$ exhausts
$X$. We put the direct limit topology on $X$, so that a subset $Y$ of $X$ is
closed (respectively, open) if and only if for all $n$ the intersection
$Y\cap Z_n$ is closed (respectively, open) in the Zariski topology on $Z_n$.
Each $Z_n$ is closed in
$X$. If
$Y$ is locally closed in
$X$, then each intersection $Y \cap Z_n$ is locally closed in $Z_n$.

\subsection{Dimensions of locally closed subsets of $X$}  We write $\mathcal
Z$ for the family of subsets $Z$ of $X$ for which there exists
$n$ such that $Z$ is a closed subset of $Z_n$. Each $Z \in \mathcal Z$ is a
projective variety, and the family $\mathcal Z$ is stable under the action of
$G(L)$. For any locally closed subset $Y$ of $X$ we put
\begin{equation}
\dim(Y):=\sup\{\dim(Y\cap Z): Z \in \mathcal Z\}.
\end{equation} Of course $\dim(Y)$ might be $+\infty$, as happens for example
when $Y=X$. Clearly
$\dim(gY)=\dim(Y)$ for all $g \in G(L)$.

\subsection{$U(L)$-orbits on $X$} For any
$\lambda
\in X_*(A)$ we put
$x_\lambda:=\epsilon^\lambda x_0 \in X$. We write $X_\lambda$ for the
$U(L)$-orbit of $x_\lambda$; by the Iwasawa decomposition we have
\begin{equation} X=\coprod_{\lambda \in X_*(A)} X_\lambda,
\end{equation} set-theoretically. The sets
\begin{equation} X_{\le \lambda}:=\coprod_{\lambda':\lambda' \le \lambda}
X_{\lambda'}
\end{equation} are closed in $X$, and the $U(L)$-orbits $X_\lambda$ are
locally closed in $X$. Taking $\lambda=0$ we get $X_0=U(L)/U(\mathfrak o)$.

\subsection{Dimensions of $A(F)$-stable locally closed subsets of
$X$}\label{subsec.af.stable}
 Let $Y$ be an $A(F)$-stable locally closed subset of
$X$. We claim that
\begin{equation}
\dim(Y)=\dim(Y \cap X_\lambda)
\end{equation} for any $\lambda \in X_*(A)$.  Clearly
$X_\lambda=\epsilon^\lambda X_0$, and because $Y$ is $A(F)$-stable and hence
satisfies $Y=\epsilon^\lambda Y$, we have
$Y \cap X_\lambda=\epsilon^\lambda (Y \cap X_0)$, so that
\begin{equation}\label{eq.dimeq}
\dim(Y \cap X_\lambda)=\dim(Y \cap X_0)
\end{equation} for all $\lambda$. Thus we need only show that
\begin{equation}
\dim(Y)=\dim(Y \cap X_0).
\end{equation}

 The inequality $\dim(Y\cap X_0) \le \dim (Y)$ is clear. For the reverse
inequality we must show that
\begin{equation}
\dim(Y\cap Z) \le \dim(Y\cap X_0)
\end{equation} for any $Z \in \mathcal Z$. It is well-known (and easy) that
there exists a finite subset $S$ of $X_*(A)$ such that $Z$ is contained in
\begin{equation*}
\bigcup_{\lambda \in S} X_\lambda.
\end{equation*} Therefore
\begin{equation}\label{eq.dimest}
\begin{split}
\dim(Y\cap Z) &=\sup\{ \dim(Y \cap Z \cap X_\lambda):\lambda \in S \}\\ &\le
\sup\{ \dim(Y  \cap X_\lambda):\lambda \in S \}.
\end{split}
\end{equation}
 Combining \eqref{eq.dimest} and \eqref{eq.dimeq}, we see that
$\dim(Y\cap Z) \le \dim(Y \cap X_0)$, as desired.

Since $X_\mu(\epsilon^\nu)$ is $A(F)$-stable, the remarks above show that its
dimension is the same as that of
\[ X_\mu(\epsilon^\nu) \cap X_0=\{u \in U(L)/U(\mathfrak o) :
u^{-1}\epsilon^\nu \sigma(u)
\in K\epsilon^\mu K \}.
\]

\subsection{Root subgroups of $U$} For any positive root $\alpha$ we write
$U_\alpha$ for the root subgroup of $U$ corresponding to $\alpha$.  Enumerate
the positive roots (in any order) as
$\alpha_1,\dots,\alpha_r$. Then, as is well-known, the map $(u_1,\dots,u_r)
\mapsto u_1\dots u_r$ is an isomorphism
\[
\prod_{i=1}^rU_{\alpha_i} \to U
\] of algebraic varieties over $k$. We now fix, for each positive root
$\alpha$, an isomorphism $U_\alpha\simeq \mathbb G_a$ over $k$. Thus we may
identify
$U(L)$ with $L^r$ and  $U(\mathfrak o)$ with $\mathfrak o^r$.

\subsection{Subgroups $U_n$ of $U(\mathfrak o)$} For any $n\ge 0$ the ring
homomorphism $\mathfrak o \twoheadrightarrow \mathfrak o/\epsilon^n
\mathfrak o$ induces a surjective group homomorphism
$ U(\mathfrak o) \twoheadrightarrow U(\mathfrak o/\epsilon^n \mathfrak o),
$ whose kernel we denote by $U_n$. Thus we have a descending chain
$U(\mathfrak o)=U_0 \supset U_1 \supset U_2 \supset \dots$ of normal subgroups
of $U(\mathfrak o)$. Under our identification of $U(\mathfrak o)$ with
$\mathfrak o^r$ the subgroup $U_n$ becomes identified with
$(\epsilon^n\mathfrak o)^r$.

\subsection{Subgroups $U(m)$ of $U(L)$}\label{sub.subU(m)} We now fix a
dominant regular coweight $\lambda_0 \in X_*(A)$; thus
$\langle \alpha,\lambda_0 \rangle > 0$ for every positive root $\alpha$ of
$A$. Put
$a := \epsilon ^{\lambda_0}$, and then for $m \in \mathbb Z$ define a subgroup
$U(m)$ of $U(L)$ by
\[ U(m):=a^mU(\mathfrak o)a^{-m}.
\] There is a chain of inclusions
\[
\dots \supset U(-2) \supset U(-1) \supset U(0) \supset U(1) \supset U(2)
\supset
\cdots
\] and moreover
\[ U(L)=\bigcup_{m \in \mathbb Z} U(m).
\] Under our identification of $U(L)$ with $L^r$, the subgroup $U(m)$ becomes
identified with
\[
\prod_{i=1}^r \epsilon^{m\langle \alpha_i,\lambda_0 \rangle } \mathfrak o.
\]

Clearly the filtrations $U_0 \supset U_1\supset U_2 \supset \dots $ and $U(0)
\supset U(1)\supset U(2) \supset \dots $  define the same topology on
$U(\mathfrak o)$.

\subsection{Dimensions of locally closed subsets $Y$ of $X_0$} Recall that
\[X_0=U(L)/U(\mathfrak o).
\]
 Therefore $X_0$ can be written as the increasing
union
\[ X_0=\bigcup _{m\ge 0}U(-m)/U(0)
\] of its closed subspaces $U(-m)/U(0)$.

Now consider  a locally closed subset $Y$ of $X_0$. Then $Y$ is locally closed
in $X$  and by definition its dimension is
$\sup\{\dim(Y \cap Z) : Z \in \mathcal Z\}$. It is clear that for any
$Z \in \mathcal Z$ there exists $m\ge 0$ such that $Z
\cap X_0 \subset U(-m)/U(0)$. It is equally clear that for any $m\ge 0$ there
exists $Z \in \mathcal Z$ such that $U(-m)/U(0) \subset Z\cap X_0$. We
conclude that
\begin{equation}
\dim(Y)=\sup\{\dim Y \cap (U(-m)/U_0) : m \ge 0 \}.
\end{equation}

We of course are particularly interested in the dimension of the locally
closed subset
\begin{equation} Y_{\mu,\nu}:=X_\mu(\epsilon^\nu)\cap X_0=\{u \in
U(L)/U(\mathfrak o):u^{-1}\epsilon^\nu \sigma(u) \in K\epsilon^\mu K\},
\end{equation} defined for any coweights $\mu,\nu$  with  $\mu$  dominant.  In
other words
\begin{equation}\label{eq.def.Ymunu} Y_{\mu,\nu}=f_\nu^{-1}\bigl(K\epsilon^\mu
K\epsilon^{-\nu} \cap U(L)\bigr)/U(\mathfrak o),
\end{equation}
 where $f_\nu$ is the map $U(L) \to U(L)$ defined by $f_\nu(u):=
u^{-1}\epsilon^\nu\sigma(u)\epsilon^{-\nu}$.

\subsection{Dimensions for admissible and  ind-admissible subsets of
$U(L)$}\label{sub.dimaisU} In view of \eqref{eq.def.Ymunu} we see that it
would be useful to introduce a notion of dimension for suitable subsets $V$
of $U(L)$ and then to compute the dimension of
$f_\nu^{-1}V$ in terms of that of  $V$.

For $m,n\ge 0$ the quotient $U(-m)/U(n)$ is the set of $\bar k$-points of an
algebraic variety over $k$. We say that a subset $V$ of $U(-m)$ is
\emph{admissible} if there exists
$n\ge 0$ and a locally closed subset $V'$ of $U(-m)/U(n)$ such that $V$ is the
full inverse image of $V'$ under $U(-m) \twoheadrightarrow U(-m)/U(n)$. We say
that a subset $V$ of $U(L)$ is \emph{admissible} if there exists
$m\ge 0$ such that $V$ is an admissible subset of $U(-m)$. We say that a
subset $V$ of $U(L)$ is
\emph{ind-admissible} if $V \cap U(-m)$ is admissible for all $m \ge 0$.
Obviously admissible subsets are also ind-admissible.

For any admissible subset $V$ of $U(L)$ we choose $n\ge 0$ such that $V$ is
preserved by right multiplication by $U(n)$ and then put
\begin{equation}
\dim V:=\dim (V/U(n))-\dim(U(0)/U(n));
\end{equation} this is clearly independent of the choice of $n$. In this
definition we could equally well have used  the subgroups $U_n$ instead of
$U(n)$.   Clearly $\dim U(\mathfrak o)=0$, and of course
$\dim V$ can be negative.

For any ind-admissible subset $V$ of $U(L)$ we put
\begin{equation}
\dim V:=\sup\{\dim V \cap U(-m): m \ge 0 \}.
\end{equation} Of course $\dim(V)$ can be $+\infty$, as happens in case
$V=U(L)$.

\subsection{Warm-up exercise and key proposition}\label{sec.warmup} The
following familiar lemma is what makes the next proposition work.
\begin{lemma}\label{lem.warm} Let $a,b$ be non-negative integers and consider
the group homomorphism
$f:\mathfrak o
\to
\mathfrak o$ defined by $f(x):=\epsilon^a \sigma(x)-\epsilon^b x$. Then the
image  of $f$ is the subgroup $\epsilon^c\mathfrak o$, where
$c:=\min\{a,b\}$.
\end{lemma}
\begin{proof} Clearly  it is enough to treat the case in which $c=0$. Then at
least one of $a,b$ is zero. If $(a,b) \ne (0,0)$, then $f$ preserves the
filtration
\[\mathfrak o \supset \epsilon \mathfrak o \supset \epsilon^2 \mathfrak
o\supset
\dots\] of
$\mathfrak o$ and induces an isomorphism on the associated graded group.
Therefore
$f$ is bijective in this case.

If $(a,b)=(0,0)$, then $f$ preserves the filtration above and induces on each
successive quotient $\bar k$ the map $x \mapsto \sigma(x)-x=x^q-x$, which is
surjective since $\bar k$ is algebraically closed. Therefore $f$ maps
$\mathfrak o$ onto $\mathfrak o$.
\end{proof}

Now we come to the key proposition. It involves two dominant coweights $\nu$
and
$\nu'$. Define homomorphisms $\phi,\psi:U(L) \to U(L)$ by
$\phi(u):=\epsilon^{\nu'}u\epsilon^{-\nu'}$ and
$\psi(u)=\epsilon^\nu\sigma(u)\epsilon ^{-\nu}$. The dominance of $\nu,\nu'$
implies that both $\phi,\psi$ preserve the normal subgroups $U_n$ of
$U(\mathfrak o)$ and hence induce homomorphisms $\phi_n,\psi_n:U(\mathfrak
o/\epsilon^n
\mathfrak o) \to U(\mathfrak o/\epsilon^n\mathfrak o) $.

We use $\phi,\psi$ to define a right action, denoted $*$, of $U(\mathfrak o)$
on itself: the action of an element $u \in U(\mathfrak o)$ upon an element $u'
\in U(\mathfrak o)$ is given by
$u'*u:=
\phi(u)^{-1}u'\psi(u)$.  Similarly, we use $\phi_n,\psi_n$ to define a right
action, again denoted $*$, of $U(\mathfrak o/\epsilon^n\mathfrak o)$ on
itself: the action of an element $u \in U(\mathfrak o/\epsilon^n\mathfrak o)$
upon an element $u' \in U(\mathfrak o/\epsilon^n\mathfrak o)$ is given by
$u'*u:=
\phi_n(u)^{-1}u'\psi_n(u)$.

We regard $U(\mathfrak o/\epsilon^n \mathfrak o)$ as the set of $\bar
k$-points of an algebraic group over $k$ in the usual way. In particular the
underlying variety is simply an affine space of dimension $n\dim U$. In the
next proposition the dimensions are those of varieties over $\bar k$.  In the
proposition, and throughout the paper, we write sums over the set of positive
roots as $\sum_{\alpha > 0}$.

\begin{proposition}\label{prop.key}  Let $\nu$ and $\nu'$ be dominant
coweights, and   let $n$ be a non-negative integer large enough that \[n\ge
\min\{\langle
\alpha,\nu\rangle, \langle \alpha,\nu' \rangle \} \] for every positive root
$\alpha$. Then we have:
\begin{enumerate}
\item The codimension of the $U(\mathfrak o/\epsilon^n\mathfrak o)$-orbit of
$1 \in U(\mathfrak o/\epsilon^n\mathfrak o)$ is equal to
\begin{equation*}
\sum_{\alpha >0} \min \{\langle \alpha,\nu \rangle,\langle \alpha,\nu'
\rangle \},
\end{equation*} which in turn is equal to the dimension of the stabilizer
$S_n$ of $1 \in U(\mathfrak o/\epsilon^n\mathfrak o)$ in
$U(\mathfrak o/\epsilon^n\mathfrak o)$. Note that $S_n=\{u \in U(\mathfrak
o/\epsilon^n\mathfrak o):\phi_n(u)=\psi_n(u) \}$.
\item The $U(\mathfrak o)$-orbit of $1 \in U(\mathfrak o)$ coincides with the
inverse image under $U(\mathfrak o) \twoheadrightarrow U(\mathfrak
o/\epsilon^n\mathfrak o)$ of the $U(\mathfrak o/\epsilon^n\mathfrak o)$-orbit
of $1 \in U(\mathfrak o/\epsilon^n\mathfrak o)$. In particular the
$U(\mathfrak o)$-orbit of
$1 \in U(\mathfrak o)$ contains $U_n$.
\end{enumerate}
\end{proposition}
\begin{proof} To simplify notation we sometimes write $H$ for the group
$U(\mathfrak o/\epsilon^n \mathfrak o)$. We begin by proving the first
statement.  The codimension of the $H$-orbit of $1 \in H$ is
\[
\dim H -\dim(H/S_n)=\dim S_n.
\] Under our identification $\mathbb G_a^r \simeq U$, the subgroup $S_n$ of
$U(\mathfrak o/\epsilon^n \mathfrak o)$ goes over to
\[
\prod_{i=1}^r \{x \in \mathfrak o/\epsilon^n \mathfrak o : \epsilon ^{\langle
\alpha_i,\nu' \rangle}  x = \epsilon ^{\langle \alpha_i,\nu \rangle}
\sigma(x) \},
\] and the dimension of the $i$-th factor in this product is equal to the
codimension of the image of the homomorphism
$x \mapsto \epsilon ^{\langle \alpha_i,\nu \rangle}
\sigma(x) -  \epsilon ^{\langle
\alpha_i,\nu' \rangle}  x $ from $\mathfrak o/\epsilon^n \mathfrak o$ to
itself,  and by Lemma \ref{lem.warm} (and our hypothesis on $n$) this
codimension is obviously
$ \min \{\langle \alpha_i,\nu \rangle,\langle \alpha_i,\nu' \rangle \}$.

Now we prove the second statement.  For $u \in U(\mathfrak o)$ we write
$\bar u$ for the image of
$u$ under $U(\mathfrak o)
\twoheadrightarrow H$. It is clear that the $*$-actions on $U(\mathfrak o)$ and
$H$ are compatible, in the sense that
\begin{equation}\label{eq.compat}
\overline{u*u'}=\bar u * \bar u',
\end{equation} from which it follows that  $U(\mathfrak o)
\twoheadrightarrow H$ maps the $U(\mathfrak o)$-orbit of $1 \in U(\mathfrak
o)$ into the $H$-orbit of $1 \in H$.  We must prove that if $u \in U(\mathfrak
o)$  has the property that $\bar u$ lies in the $H$-orbit of
$1$, then $u$ lies in the
$U(\mathfrak o)$-orbit of $1$. Replacing $u$ by $u*u'$ for suitable $u'\in
U(\mathfrak o)$, we may assume that
$\bar u=1$.  In other words it remains only to show that any $u \in U_n$ lies
in the $U(\mathfrak o)$-orbit of $1 \in U(\mathfrak o)$. Some care is needed
since the
$*$-action of
$U(\mathfrak o)$ on itself does not preserve the subgroup $U_n$.

Until now we have been working with any ordering $\alpha_1,\dots,\alpha_r$ of
the positive roots. Now let us order them so that $\hgt(\alpha_1) \le
\hgt(\alpha_2)
\le \dots \le \hgt(\alpha_r)$, where $\hgt(\alpha)$ is the number of simple
roots needed in order to write the positive root $\alpha$ as a sum of simple
roots.   We then get a decreasing chain
\[ U=U[1] \supset U[2] \supset U[3] \supset \dots \supset U[r+1]=\{1\}
\] of normal subgroups $U[j]$ of $U$, with $U[j]$ defined as the subgroup
consisting of elements whose projections onto the first $j-1$ root subgroups
$U_{\alpha_1},\dots,U_{\alpha_{j-1}}$ are all equal to $1$.

Now return to our element $u \in U_n$, and let $u_1$ denote the image of $u$ in
$U[1](\mathfrak o)/U[2](\mathfrak o)=U_{\alpha_1}(\mathfrak o)$.    By Lemma
\ref{lem.warm} (and our hypothesis on
$n$) we may choose
$v
\in U_{\alpha_1}(\mathfrak o)$ so that
$u_1*v =1$. It follows that $u*v \in U[2](\mathfrak o)$. We claim that $u*v
\in U_n$. For this we must check that
$\overline{u*v}=1$, which is clear since
\[\overline{u*v}=\bar u *\bar v=1*\bar v=\bar u_1 * \bar
v=\overline{u_1*v}=1.\]  Replacing
$u$ by
$u*v$, we may now assume that
$u
\in U[2](\mathfrak o) \cap U_n$. Using $U_{\alpha_2}$ the same way we just
used $U_{\alpha_1}$,  we may push $u$ down into
$U[3](\mathfrak o) \cap U_n$. Continuing in this way, we eventually end up
with $u
\in U[r+1](\mathfrak o)=\{1\}$, at which point we are done.
\end{proof}

\subsection{Formula for the dimension of $f^{-1}_{\nu,\nu'}V$} Now let
$\nu,\nu' \in X_*(A)$ and define a map $f_{\nu,\nu'}:U(L) \to U(L)$ by
\begin{equation}
f_{\nu,\nu'}(u):=(\epsilon^{\nu'}u\epsilon^{-\nu'})^{-1}\cdot(\epsilon^\nu
\sigma(u)\epsilon^{-\nu})=\epsilon^{\nu'}u^{-1}\epsilon^{\nu-\nu'}
\sigma(u)\epsilon^{-\nu}.
\end{equation} Note that when $\nu'=0$, we get back the map $f_\nu$ considered
earlier. Moreover, when $\nu,\nu'$ are both dominant and $u \in U(\mathfrak
o)$, we have
$f_{\nu,\nu'}(u)=1*u$, with $*$ denoting the action of $U(\mathfrak o)$ on
itself introduced in
\ref{sec.warmup}.

Now let $\lambda \in X_*(A)$. We will also  need  the conjugation map
$c_\lambda:U(L) \to U(L)$ defined by
\begin{equation} c_\lambda(u):=\epsilon^\lambda u \epsilon^{-\lambda}.
\end{equation}

\begin{proposition}\label{prop.main} Let $\lambda,\nu,\nu' \in X_*(A)$. Then
we have:
\begin{enumerate}
\item $f_{\nu,\nu'}\, c_\lambda=f_{\lambda+\nu,\lambda+\nu'}=c_\lambda\,
f_{\nu,\nu'}$.
\item Let $V$ be an admissible {\textup(}respectively,
ind-admissible\textup{)} subset of
$U(L)$. Then
$c^{-1}_\lambda V$ is admissible {\textup(}respectively,
ind-admissible\textup{)} and
\begin{equation*}
\dim c^{-1}_\lambda V=\dim V + \sum_{\alpha >0} \langle \alpha,\lambda
\rangle.
\end{equation*}
\item Let $V$ be an admissible subset of $U(L)$. Then $f^{-1}_{\nu,\nu'}V$ is
ind-admissible and
\begin{equation*}
\dim f^{-1}_{\nu,\nu'} V=\dim V + \sum_{\alpha >0} \min\{\langle \alpha,\nu
\rangle, \langle \alpha,\nu' \rangle\}.
\end{equation*}
\end{enumerate}
\end{proposition}
\begin{proof} The first statement of the proposition is an easy calculation.
Next we prove the second statement. Assume that $V$ is an admissible subset of
$U(L)$, and pick $n$ big enough that $V$ is stable under right multiplication
by $U_n$. Since
$c_\lambda$ induces an isomorphism from the  variety
$c^{-1}_\lambda V/c^{-1}_\lambda U_n$ to the variety $V/U_n$, we see that
\[
\dim c^{-1}_\lambda V-\dim c^{-1}_\lambda U_n=\dim V-\dim U_n,
\] from which it follows that it is enough to prove the second statement of
the proposition for $V=U_n$. Using the root subgroups $U_\alpha$, we are
reduced to the obvious fact that
\[
\dim \epsilon^a\mathfrak o/\epsilon^b\mathfrak o=b-a.
\]

Finally we prove the third statement of the proposition. Using the first two
parts of the proposition, one sees easily that the third part is true for
$(\nu,\nu')$ if and only if it is true for $(\nu+\lambda,\nu'+\lambda)$.
Therefore we may assume without loss of generality that both $\nu,\nu'$ are
dominant. From now on we abbreviate $f_{\nu,\nu'}$ to $f$. Note that the
dominance of $\nu,\nu'$ implies that $f(U_0)\subset U_0$; indeed, for $u \in
U_0$ we have $f(u)=1*u$ in the notation of \ref{sec.warmup}, so that
$f(U_0)$ coincides with the $U_0$-orbit of
$1 \in U_0$, something we have analyzed in Proposition \ref{prop.key}. We write
$f_0$ for the map $f_0:U_0 \to U_0$ obtained by restriction from $f$.

As in Proposition \ref{prop.key} we now choose a non-negative integer $n_0$
large enough that
 \[n_0 \ge \min\{\langle
\alpha,\nu\rangle, \langle \alpha,\nu' \rangle \} \] for every positive root
$\alpha$.  Using the first two parts of Proposition \ref{prop.main} again, we
see that the third part of the proposition is true for $V$ if and only if it
is true for
$\epsilon^\lambda V
\epsilon^{-\lambda}$. Therefore we may assume without loss of generality that
$V
\subset U_{n_0}$.

Under this assumption we are going to prove that for all $m \ge 0$ the
intersection $U(-m) \cap f^{-1}V$ is admissible of dimension
\[
\dim V+\sum_{\alpha >0} \min \{\langle \alpha,\nu \rangle,\langle \alpha,\nu'
\rangle \},
\]  which will be enough to prove the proposition.

Recall that $U(-m)=\epsilon^{-m\lambda_0}U_0\epsilon^{m\lambda_0}$. Once again
using the first two parts of this proposition, we see that
\[ U(-m)\cap f^{-1}V=\epsilon^{-m\lambda_0}\bigl( U_0 \cap
f^{-1}(\epsilon^{m\lambda_0}V\epsilon^{-m\lambda_0})\bigr)\epsilon^{m\lambda_0}
\] and  that
\[
\dim U(-m)\cap f^{-1}V - \dim V=\dim U_0 \cap
f^{-1}(\epsilon^{m\lambda_0}V\epsilon^{-m\lambda_0})-\dim
\epsilon^{m\lambda_0}V\epsilon^{-m\lambda_0}.
\]  Since $\epsilon^{m\lambda_0}V\epsilon^{-m\lambda_0}$ is still contained in
$U_{n_0}$, we may without loss of generality take $m=0$.

Thus we are now reduced to proving that if $V \subset U_{n_0}$, then
$f_0^{-1}V$ is admissible of dimension
\[
\dim V+\sum_{\alpha >0} \min \{\langle \alpha,\nu \rangle,\langle \alpha,\nu'
\rangle \}.
\] Choose $n \ge n_0$ such that $V$ is invariant under right multiplication by
$U_n$, and once again put $H:=U_0/U_n$.  Then $V$ is obtained as the inverse
image under
\[ p:U_0 \twoheadrightarrow H
\] of a locally closed subset $\bar V$ of $H$ contained in $U_{n_0}/U_n
\subset H$.

Writing $\bar f$ for the map $H \to H$ defined by $\bar f(h):=1*h$, we have
$pf_0=\bar fp$, from which it follows that
\[ f_0^{-1}V=p^{-1}\bar f^{-1} \bar V.
\]  In particular $f_0^{-1}V$ is admissible of dimension
\[
\dim U_n+\dim\bar f^{-1}\bar V.
\]

Recall from before that $S_n$ denotes the stabilizer group-scheme of $1 \in H$
for the
$*$-action of $H$ on itself.  Now the morphism $\bar f$ induces an isomorphism
from the fpqc quotient
$S_n\backslash H$ to the $H$-orbit of $1 \in H$ (see Lemma
\ref{orbit_lemma}), and by the second part of Proposition
\ref{prop.key}, applied to $n_0$,  we have
\[
\bar V \subset U_{n_0}/U_n \subset \bar f(H).
\]  Therefore by Lemma \ref{orbit_lemma}
\[
\dim \bar f^{-1}\bar V=\dim \bar V+\dim S_n,
\] from which it follows that
\[
\dim f_0^{-1} V=\dim U_n+\dim \bar V+\dim S_n=\dim V+\dim S_n.
\] Using the formula for $\dim S_n$ given in the first part of Proposition
\ref{prop.key}, we see that
\[
\dim f_0^{-1} V=\dim V+\sum_{\alpha >0} \min \{\langle \alpha,\nu
\rangle,\langle
\alpha,\nu'
\rangle \},
\]  as desired.
\end{proof}

\subsection{Formula for $\langle \rho,\nu-\nu_{\dom} \rangle$} Let $\rho \in
X^*(A)_\mathbb Q$ be the half-sum of the positive roots of $A$. Recall that
$\nu_{\dom}$ denotes the unique dominant element in the $W$-orbit of
$\nu$.
\begin{lemma}\label{lem.rho.dom} For any $\nu \in \mathfrak a$ there is an
equality
\begin{equation}
\langle \rho,\nu-\nu_{\dom} \rangle=\sum_{\alpha >0} \min\{\langle \alpha,\nu
\rangle ,0 \}.
\end{equation}
\end{lemma}
\begin{proof} The lemma follows immediately from the equations
\begin{equation}
\langle \rho,\nu \rangle= \sum_{\alpha > 0} \langle \alpha,\nu \rangle/2
\end{equation} and
\begin{equation}
\langle \rho,\nu_{\dom} \rangle= \sum_{\alpha >0} |\langle \alpha,\nu
\rangle|/2.
\end{equation}
\end{proof}

\subsection{Mirkovi\'c-Vilonen dimension formula}
 Recall that
$X_\nu$ denotes the $U(L)$-orbit of $x_\nu=\epsilon^\nu x_0$ in $X$. For
dominant
$\mu
\in X_*(A)$  we also consider the $K$-orbit $X^\mu:=Kx_\mu$ of $x_\mu$ in
$X$.  We are interested in the intersection $X^\mu\cap X_\nu$. We assume that
$\nu_{\dom} \le \mu$, since otherwise the intersection is empty (see
\cite[4.4.4]{bruhat-tits72}).  As before $\rho$ is the half-sum of the
positive roots. Then we have
\begin{proposition}[Mirkovi\'c-Vilonen]\label{orig.mirk.vil}  For any dominant
$\mu \in X_*(A)$ and any $\nu \in X_*(A)$ such that $\nu_{\dom}
\le \mu$ there is an equality
\begin{equation}\label{eq.mv}
\dim X^\mu \cap X_\nu=\langle \rho,\mu+\nu \rangle.
\end{equation}
\end{proposition} This is proved in \cite{MV1,MV2}.  Another proof in the
present context of finite residue fields is given in section
\ref{reduction.section}, see Remark \ref{alt.mirk.vil}.

We will need the following  reformulation of their result.
\begin{proposition}\label{prop.mirk.vil} For any dominant $\mu \in X_*(A)$ and
any $\nu \in X_*(A)$ such that $\nu_{\dom}
\le \mu$ there is an equality
\begin{equation}
\dim K\epsilon^\mu K \epsilon^{-\nu} \cap U(L) =\langle \rho,\mu-\nu \rangle.
\end{equation}
\end{proposition}
\begin{proof} The variety $X^\mu \cap X_\nu$ is isomorphic to the quotient
\[(K\epsilon^\mu K
\epsilon^{-\nu} \cap U(L))/\epsilon^\nu U_0\epsilon^{-\nu},
\]
 and $K\epsilon^\mu K
\epsilon^{-\nu} \cap U(L)$ is obviously an admissible subset of $U(L)$. Thus,
according to our discussion of  dimension theory for admissible subsets of
$U(L)$, we have
\begin{equation}\label{eq.kmuKnu}
\dim X^\mu \cap X_\nu=\dim K\epsilon^\mu K
\epsilon^{-\nu} \cap U(L) - \dim \epsilon^\nu U_0\epsilon^{-\nu}.
\end{equation} From the second part of Proposition \ref{prop.main} we obtain
\begin{equation}\label{eq.sc.nu}
\dim \epsilon^\nu U_0\epsilon^{-\nu}=-2\langle \rho,\nu \rangle.
\end{equation} Combining \eqref{eq.mv}, \eqref{eq.kmuKnu}, \eqref{eq.sc.nu},
we get the desired result.
\end{proof}

\subsection{Main result} Now we are ready to calculate the dimension of
$X_\mu(\epsilon^\nu)$. Recall that
$\rho$ denotes the half-sum of the positive roots of $A$.

\begin{theorem} \label{main_result_thm} Let $\mu,\nu$ be coweights of $A$.
Assume that $\mu$ is dominant and  that
$\nu_{\dom} \le \mu$. Then
\[
\dim X_\mu(\epsilon^\nu)=\langle \rho,\mu-\nu_{\dom} \rangle.
\]
\end{theorem}
\begin{proof} In \ref{subsec.af.stable} we saw that
\begin{equation}
\dim X_\mu(\epsilon^\nu)=\dim Y_{\mu,\nu},
\end{equation} where $Y_{\mu,\nu}$ denotes the intersection
$X_\mu(\epsilon^\nu) \cap X_0$.  From \eqref{eq.def.Ymunu} we have
\begin{equation} Y_{\mu,\nu}=f_\nu^{-1}\bigl(K\epsilon^\mu K\epsilon^{-\nu}
\cap U(L)\bigr)/U(\mathfrak o).
\end{equation} Applying Proposition \ref{prop.main} with $\nu'=0$ (so that
$f_{\nu,\nu'}=f_\nu$), we see that
\begin{equation}
\dim Y_{\mu,\nu} =\dim K\epsilon^\mu K\epsilon^{-\nu} \cap U(L) +
\sum_{\alpha >0} \min\{\langle \alpha,\nu
\rangle, 0\}.
\end{equation} Using the Mirkovi\'c-Vilonen dimension formula (Proposition
\ref{prop.mirk.vil}) and Lemma \ref{lem.rho.dom}, we obtain
\begin{equation}
\dim Y_{\mu,\nu} =\langle \rho,\mu-\nu \rangle + \langle \rho,\nu-\nu_{\dom}
\rangle=\langle \rho,\mu-\nu_{\dom} \rangle,
\end{equation} as desired.
\end{proof}

\subsection{Remark} In the last proof we only used Proposition
\ref{prop.main} in the special case
$\nu'=0$. The reason we took arbitrary $\nu'$ in Proposition \ref{prop.main}
was to allow a reduction to the case of dominant coweights in the  proof of
that proposition. When $\nu'=0$ and $\nu$  is dominant, this reduction step is
not needed, and since $X_\mu(\epsilon^\nu)$ depends only on the
$\sigma$-conjugacy class of $\epsilon^\nu$, hence only on the $W$-orbit of
$\nu$, it would have been enough to consider only dominant $\nu$.

Thus one might think it a waste of effort to have formulated and proved
Proposition
\ref{prop.main} in the generality that we did. However, later in the paper,
when we switch from the affine Grassmannian to the affine flag manifold, we
will need Proposition \ref{prop.main} for $\nu'=0$ and arbitrary $\nu$. For
non-dominant $\nu$ no simplification occurs when $\nu'=0$, and the most
natural formulation of Proposition \ref{prop.main} seems to be the one we have
given.

\subsection{Equidimensionality}

Next we will prove that $X_\mu(\epsilon^\nu)$ is equidimensional, in a sense
we will make precise below.  The equidimensionality follows easily from the
corresponding fact about the intersections $X^\mu \cap X_\nu$ due to
Mirkovi\'c-Vilonen (see Lemma \ref{MV_equidim}), and from a basic lemma
concerning fpqc (or fppf) quotients of algebraic groups by stabilizer subgroup
schemes.

We say an ind-admissible set $Y \subset U(L)$ is {\em irreducible} if for all
$m \geq 0$, the intersection $Y \cap U(-m)$ is the full inverse image of an
irreducible locally closed subset $\overline{Y} \subset U(-m)/U(n)$, for some
$n \geq 0$.  Note that if $Y \subset U(-m)$ is the full inverse image of a
locally closed subset $\overline{Y}_n \subset U(-m)/U(n)$, then
$\overline{Y}_n$ is irreducible if and only if its inverse image
$\overline{Y}_{n'} \subset U(-m)/U(n')$ is irreducible for some (equivalently,
all)
$n' \geq n$.  This holds because the canonical projection $U(-m)/U(n')
\rightarrow U(-m)/U(n)$ can be realized in suitable coordinates as a
projection $\mathbb A^a \times \mathbb A^b \rightarrow \mathbb A^a$ for some
non-negative integers $a$ and $b$.

We say an ind-admissible set $Y \subset U(L)$ is {\em equidimensional} if for
each $m \geq 0$, the intersection $Y \cap U(-m)$ is the full inverse image of
an equidimensional locally closed subvariety of $U(-m)/U(n)$ for some
(equivalently, all) sufficiently large $n \geq 0$.

Finally, an $A(F)$-stable locally closed subset $Y \subset X$ will be termed
equidimensional provided that the full inverse image in $U(L)$ of $Y \cap X_0
\subset U(L)/U(\mathfrak o)$ is an equidimensional ind-admissible set.

\begin{proposition} \label{equidim.prop} The affine Deligne-Lusztig variety
$X_\mu(\epsilon^\nu)$ is equidimensional.
\end{proposition}

To prove this proposition, we will need a few lemmas.

To state the first lemma, suppose that $H$ is an algebraic group over $k$,
that is, a reduced and irreducible $k$-group scheme of finite type.  Suppose
$H$ acts (on the left) on a variety, or more generally on a finite-type
$k$-scheme $\mathfrak X$.  For a closed point $x \in \mathfrak X(k)$, we let
$\mathcal O$ denote the $H$-orbit of $x$, a locally-closed subset of
$\mathfrak X$ which we give the reduced subscheme
structure.  Let $H_x \subset H$ denote the subgroup scheme which stabilizes
the point $x$, and let ${H/H_x}$ denote the fppf-sheaf associated to
the fppf-presheaf $R \mapsto H(R)/H_x(R)$. Finally,  let $\mathfrak p_x:H \to
\mathcal O$ be the morphism $h \mapsto hx$.

\begin{lemma} \label{orbit_lemma}
\begin{enumerate}
\item [(1)]  The morphism ${\mathfrak p}_x: H \rightarrow \mathcal O$ is fppf,
hence it is an epimorphism of fppf-sheaves and induces an isomorphism of
fppf-sheaves
\[
{H/H_x} \xrightarrow{\sim} \mathcal O.
\]
\item[(2)]  For every irreducible subset $Y \subset \mathcal O$, the inverse
image scheme ${\mathfrak p}_x^{-1}(Y)$ is a finite union of irreducible
subsets $Z$, each having dimension
$$
\dim Z = \dim Y + \dim H_x.
$$
\end{enumerate}
\end{lemma}

\begin{proof} Part (1) is contained in \cite{DG}, II, $\S3$, 5.2.  For part
(2), we have by loc.~cit.~III, $\S3$, 5.5 that $\dim \mathcal O = \dim H -
\dim H_x$.     Since $H \to \mathcal O$  is a dominant morphism of irreducible
varieties,
 there exists a dense open subset
$U \subset \mathcal O$ with the following property: for every irreducible set
$Y$ which meets $U$ and every irreducible component $W$ of
${\mathfrak p}_x^{-1}(Y)$ which meets ${\mathfrak p}_x^{-1}(U)$, we have
$$
\dim W - \dim Y = \dim H - \dim \mathcal O = \dim H_x,
$$ see e.g.~\cite{Mum}, I $\S8$, Thm.~3.  Since $H$ acts transitively on
${\mathcal O}$, the latter is covered by open subsets having the same
property, and this implies (2).
\end{proof}

\begin{lemma} \label{inverse_image_equidim} Let $\nu,\nu' \in X_*(A)$.  If
$V \subset U(L)$ is an equidimensional ind-admissible set, then so is
$f^{-1}_{\nu,\nu'}(V)$.
\end{lemma}

\begin{proof} We may assume $V$ is irreducible. We may also assume that $V$
is admissible, since the intersection
$f^{-1}_{\nu,\nu'}(V) \cap U(-m)$ does not change when $V$ is replaced by $V
\cap U(-m')$ for suitably large $m'$.

As in the proof of Proposition \ref{prop.main}, we may replace $(\nu,\nu')$
with $(\lambda + \nu, \lambda + \nu')$ for a suitably dominant $\lambda$ and
thereby reduce to the case where $\nu$ and $\nu'$ are dominant and $f =
f_{\nu,\nu}$ satisfies $f(U_0) \subset U_0$. Also, we may replace $V$ with
$\epsilon^\lambda V\epsilon^{-\lambda}$ and thus we may assume $V \subset
U_{n_0}$, where $n_0$ is chosen as in Proposition \ref{prop.main}.  We claim
that for any $m \geq 0$, the intersection $U(-m) \cap f^{-1}V$ is a union of
finitely many irreducible admissible sets all having the same dimension.  As
in Proposition \ref{prop.main}, it is enough to treat the case $m = 0$.

Recall that $H := U_0/U_n$ and $f_0$ induces the map $\bar{f}: H \rightarrow
H$, which is the right action of $H$ on $1 \in H$: $\bar{f}(h) = 1 * h$. Now
as before for any $n \geq n_0$ we have the commutative diagram
$$
\xymatrix{ U_0 \ar[r]^{f_0} \ar[d]^p & U_0 \ar[d]^p \\ U_0/U_n
\ar[r]^{\bar{f}} & U_0/U_n }
$$ and $\overline{V} \subset U_{n_0}/U_n \subset \bar{f}(H)$, where all
notation is as in Proposition \ref{prop.main}, except that $\overline{V}$ is
now an irreducible subset of $U_0/U_n$.

Since $f_0^{-1}(V) = p^{-1}\bar{f}^{-1}(\overline{V})$, we need only show that
$\bar{f}^{-1}(\overline{V})$ is the union of finitely many irreducible locally
closed subsets of $U_0/U_n$, all having the same dimension.  This follows from
Lemma \ref{orbit_lemma}.
\end{proof}

The next lemma was established by Mirkovi\'c and Vilonen as part of their
proof of the geometric Satake isomorphism.

\begin{lemma} [Mirkovi\'c-Vilonen] \label{MV_equidim} Let $\mu$ be
$G$-dominant and let $\nu \in X_*(A)$ be an element such that $\nu_{\rm dom}
\leq \mu$.  Then the variety $X^\mu \cap X_\nu$ is equidimensional.
\end{lemma}

\begin{proof} We shall deduce the lemma from the following result.  Let
${\rm IC}_\mu$ denote the intersection complex on the closure of $X^\mu$,
normalized in such a way that ${\rm IC}_\mu [\langle 2\rho, \mu
\rangle](\langle \rho, \mu \rangle)$ is a self-dual perverse sheaf of weight
$0$.  Then the complex ${\rm R}\Gamma_c(X_\nu,{\rm IC}_\mu)$ is concentrated
in cohomological degree $\langle 2\rho, \mu + \nu \rangle$.  This is a part of
the geometric Satake isomorphism \cite{MV1,MV2}, and it was also proved in a
direct fashion by Ng\^o-Polo in \cite{NP}.  Our goal is simply to take it as
an input and show how it implies the equidimensionality statement.

Suppose $X^\mu \cap X_\nu$ is not equidimensional, and choose an irreducible
component $C$ of dimension $d < \langle \rho, \mu + \nu \rangle$.  Let $U$
denote the complement in $X^\mu \cap X_\nu$ of all the other irreducible
components.  Thus $U$ is open and dense in $C$ and it is also open in
$\overline{X^\mu} \cap X_\nu$.  Let $Z$ denote the complement of $U$ in
$\overline{X^\mu} \cap X_\nu$.

Consider the exact sequence
$$ H^{2d-1}_c(Z,{\rm IC}_\mu) \rightarrow H^{2d}_c(U,{\rm IC}_\mu)
\rightarrow H^{2d}_c(\overline{X^\mu}\cap X_\nu, \, {\rm IC}_\mu).
$$ The first map is identically zero for weight reasons: the complex ${\rm
IC}_\mu$ has weights $\leq 0$, and so by \cite{D}, the left-most group has
weights $\leq 2d-1$.  On the other hand the middle group is pure of weight
$2d$, since ${\rm IC}_\mu|_U$ is the constant sheaf.  This shows that the
non-zero middle group injects into $H_c^{2d}(X_\nu,{\rm IC}_\mu)$, a
contradiction.
\end{proof}

Finally, note that Proposition \ref{equidim.prop} follows immediately from
Lemmas \ref{inverse_image_equidim} and \ref{MV_equidim} and the following
identities we already proved:
\begin{align*} X_\mu(\epsilon^\nu) \cap X_0 &= f_{\nu}^{-1}(K\epsilon^\mu
K\epsilon^{-\nu} \cap U(L))/U(\mathfrak o) \\  X^\mu \cap X_\nu &=
(K\epsilon^\mu K\epsilon^{-\nu} \cap U(L))/ \epsilon^{\nu}U(\mathfrak o)
\epsilon^{-\nu}.
\end{align*}

\section{Dimension theory for admissible subsets of $L^n$}  Let $V$ be a
finite dimensional vector space over $L$. We discuss admissibility and
dimensions for subsets of $V$. This will be needed in the next section.
\subsection{Lattices} Recall that a \emph{lattice} $\Lambda$ in $V$ is a
finitely generated
$\mathfrak o$-submodule of $V$ that generates $V$ as vector space. If
$V$ is $n$-dimensional, then every lattice in $V$ is free of rank $n$ as
$\mathfrak o$-module. Given two lattices $\Lambda_1$, $\Lambda_2$ in
$V$,  there exist $j,l \in \mathbb Z$ such that $\epsilon^j
\Lambda_2 \subset
\Lambda_1
\subset \epsilon^l\Lambda_2$, and moreover $\Lambda_1 \cap \Lambda_2$ and
$\Lambda_1 + \Lambda_2$ are again lattices.

\subsection{Admissible subsets of $V$} We say that a subset $Y \subset V$ is
\emph{admissible} if there exist lattices $\Lambda_2 \subset
\Lambda_1$ and a Zariski locally closed subset $\bar Y$ of the finite
dimensional $\bar k$-vector space $\Lambda_1/\Lambda_2$ such that $Y$ is the
inverse image under $\Lambda_1 \twoheadrightarrow
\Lambda_1/\Lambda_2$ of $\bar Y$; in this situation we say that $Y$ comes from
$\bar Y$. We say that a subset $Y \subset V$ is
\emph{ind-admissible} if for every lattice $\Lambda$ the intersection $Y
\cap \Lambda$ is admissible.  Note that a subset $Y$ is ind-admissible if $Y
\cap \Lambda_\alpha$ is  admissible for all ${\Lambda_\alpha}$ in a cofinal
family of lattices.

\subsection{Dimension theory for admissible subsets of $V$} Let
$\Lambda$ be a lattice and let $Y$ be an admissible subset of $V$. Choose
lattices $\Lambda_2 \subset \Lambda_1$ such that $\Lambda_2
\subset \Lambda$ and $Y $ comes from $\bar Y \subset
\Lambda_1/\Lambda_2$, and put
\[
\dim_\Lambda Y:=\dim \bar Y-\dim \Lambda/\Lambda_2,
\] a number independent of the choice of $\Lambda_1$, $\Lambda_2$. This notion
of dimension depends of course on $\Lambda$, as the notation indicates. We
have normalized things so that $\dim \Lambda=0$. However, for any two
admissible subsets $Y_1$, $Y_2$ the difference
\[
\dim_\Lambda Y_1-\dim_\Lambda Y_2
\] is independent of $\Lambda$, and so we will permit ourselves to write such
differences as simply
\[
\dim Y_1-\dim Y_2.
\]

For any ind-admissible subset $Y \subset V$ we define $\dim_\Lambda Y$ to be
the supremum of the numbers $\dim_\Lambda Y\cap \Lambda_1$ as
$\Lambda_1$ ranges through all lattices.

\begin{lemma} Let $g$ be an $L$-linear automorphism of $V$. Then for any
admissible subset $Y \subset V$ the inverse image $g^{-1}Y$ is admissible and
\[
\dim g^{-1}Y - \dim Y=\val \det g.
\]
\end{lemma}
\begin{proof} Easy.
\end{proof}

\section{Relative dimension of certain morphisms $f:V \to V$}  This section
contains a generalization of Lemma \ref{lem.warm} that will be needed in the
next section.
\subsection{Review of $F$-spaces} Recall that an \emph{$F$-space} is a pair
$(V,\Phi)$, with
$V$   a finite dimensional
$L$-vector space  and $\Phi$ a $\sigma$-linear bijection $V \to V$. The
category of $F$-spaces is an $F$-linear abelian category in which every object
is a finite direct sum of simple objects. The isomorphism classes of simple
objects are in natural bijection with
$\mathbb Q$, and a simple $F$-space corresponding to $\lambda \in
\mathbb Q$ is said to have \emph{slope} $\lambda$. For
$\lambda
\in
\mathbb Q$ we write
$V_\lambda$ for the sum of all the simple sub-$F$-spaces of $V$ having slope
$\lambda$. We then have the \emph{slope decomposition}
\[ V=\bigoplus_{\lambda \in \mathbb Q} V_\lambda.
\] Suppose that $V$ is $n$-dimensional. We get an unordered family of $n$
rational numbers by including each $\lambda \in \mathbb Q$ in the family with
multiplicity equal to $\dim V_\lambda$. These $n$ rational numbers are called
the \emph{slopes of $V$}.

Let $\lambda=s/r \in \mathbb Q$, with $r \ge 1$ and $(r,s)=1$. We get  a simple
$F$-space of slope $\lambda$ by taking $V=L^r$ and
$\Phi_\lambda (a_1,\dots,a_r)=(\epsilon^s\sigma a_r,\sigma a_1,\dots,
\sigma a_{r-1})$.  It follows that for any $F$-space $V$ the number
$\lambda \dim V_\lambda$ is an integer.

\subsection{Definition of $f:V\to V$ and $d(V,\Phi)$} For any $F$-space
$(V,\Phi)$ we define an $F$-linear map $f:V \to V$ by
\[ f(v):=\Phi(v)-v.
\] Moreover we define an integer $d(V,\Phi)$ by
\[ d(V,\Phi):=\sum_{\lambda < 0} \lambda \dim V_\lambda.
\]

\begin{lemma}\label{lem.preVkey} The following statements hold for $V$,
$\Phi$, $f$, $d(V,\Phi)$ as above.
\begin{enumerate}
\item The map $f$ is surjective.
\item For any lattices $\Lambda_2 \subset \Lambda_1$ the group
\[ (\Lambda_1 \cap \ker f)/(\Lambda_2 \cap \ker f)
\]  is finite.
\item There exists a lattice $\Lambda_0$ such that
\begin{enumerate}
\item $\Lambda_0':=f\Lambda_0$ is a lattice, and
\item $\dim \Lambda_0 -\dim \Lambda_0' = d(V,\Phi)$.
\end{enumerate}
\item For any lattice $\Lambda$ there exist lattices $\Lambda_1$,
$\Lambda_2$ such that $\Lambda_1 \subset f\Lambda \subset \Lambda_2$.
\item For any lattice $\Lambda'$ the inverse image $f^{-1}\Lambda'$ contains a
lattice and is ind-admissible.
\item If $\Lambda$, $\Lambda'$ are lattices such that $f\Lambda \supset
\Lambda'$, then
\[
\dim \Lambda \cap f^{-1}\Lambda' -\dim \Lambda' = d(V,\Phi).
\]
\item For any lattice $\Lambda'$ there is an equality
\[
\dim  f^{-1}\Lambda' -\dim \Lambda' = d(V,\Phi).
\]
\end{enumerate}
\end{lemma}

\begin{proof} We begin by proving the first three parts of the lemma for the
simple $F$-spaces $(L^r,\Phi_\lambda)$ described above. Recall that
$\lambda=s/r$ with $(r,s)=1$. We consider three cases.

First suppose that $s>0$. Then $f$ preserves the standard filtration
$\epsilon^n\mathfrak o^r$ on $L^r$, and induces the obviously bijective map
\[
\bar f(a_1,\dots,a_r)=(-a_1,\sigma a_1-a_2,\dots,\sigma a_{r-1}-a_r)
\] on each successive quotient $\bar k^r$. Therefore $f$ is bijective and
$f(\mathfrak o^r)=\mathfrak o^r$. In this case (1) and (2) are clear, and (3)
holds for $\Lambda_0=\mathfrak o^r$.

Next suppose that $s=0$. Then $r=1$ and we are once again looking at
$f:L \to L $ defined by $f(x)=\sigma(x)-x$. From previous work we know that
$f(\mathfrak o)=\mathfrak o$ and hence (since $f$ commutes with multiplication
by $\epsilon^a$) that
$f(\epsilon^a\mathfrak o)=\epsilon^a\mathfrak o$ for all $a \in \mathbb Z$.
Therefore
$f$ is surjective. Its kernel is obviously $F$.  So (1) and (2) are clear and
(3) holds for $\Lambda_0=\mathfrak o$.

Finally suppose that $s < 0$. Put $\Lambda_0=\mathfrak o^r$ and
$\Lambda_0'=\epsilon^s\mathfrak o \oplus \mathfrak o^{r-1}$. It is then clear
that
\[ f\Lambda_0 \subset \Lambda_0'.
\] Again using that $f$ commutes with multiplication by powers of
$\epsilon$, we see that
\[ f(\epsilon^a\Lambda_0) \subset \epsilon^a\Lambda_0'
\] for all $a \in \mathbb Z$. Thus $f$ induces a map from the associated
graded group for the filtration $\epsilon^a\Lambda_0$ to the one for the
filtration $\epsilon^a\Lambda_0'$. All successive quotients for both
filtrations can be identified with $\bar k^r$ in an obvious way, and when this
is done the map induced by $f$ on the successive quotients is always the
obviously bijective map
\[
\bar f(a_1,\dots,a_r)=(\sigma a_r,\sigma a_1-a_2,\dots,\sigma a_{r-1}-a_r).
\] We conclude that $f$ is bijective and that $f\Lambda_0=\Lambda_0'$. In this
case (1) and (2) are clear and (3) holds for the lattice
$\Lambda_0$ we have chosen, since $\dim \Lambda_0-\dim
\Lambda_0'=s=d(V,\Phi)$.

Next we prove the first three parts of the lemma in general. Write $V$ as a
direct sum of simple $F$-spaces. Then $f$ decomposes accordingly as a direct
sum. Since each summand of $f$ is surjective, so is $f$. As for (2) we are
free to enlarge $\Lambda_1$ and shrink $\Lambda_2$, so we may assume that both
$\Lambda_1$ and $\Lambda_2$ are direct sums of lattices in the simple summands
of $V$, which allows us to reduce to the simple case that has already been
treated. To prove (3) we note that we can find a suitable
$\Lambda_0$ by taking  a direct sum of suitable lattices in the simple
summands.

Now we are going to deduce the remaining parts of lemma from the three parts
we have already proved. We will no longer need to reduce to the simple case.
We begin by choosing a lattice $\Lambda_0$ as in (3). In particular
$\Lambda_0':=f\Lambda_0$ is a lattice. To see that (4) is true  we pick
$a,b$ such that
\[
\epsilon^a \Lambda_0 \subset \Lambda \subset \epsilon^b \Lambda_0.
\] Then
\[
\epsilon^a \Lambda_0' \subset f\Lambda \subset \epsilon^b \Lambda_0'.
\]

Next we prove (5). Pick $j$ such that $\epsilon^j\Lambda_0' \subset
\Lambda'$. Then $f^{-1}\Lambda'$ contains the lattice
$\epsilon^j\Lambda_0$. To show that $f^{-1}\Lambda'$ is ind-admissible it is
enough to show that for $l \ge 0$ the set
$Y_l:=\epsilon^{j-l}\Lambda_0 \cap f^{-1}\Lambda'$ is admissible. This is
clear since $Y_l$ comes from the closed subset $\bar Y_l$ of
$\epsilon^{j-l}\Lambda_0/\epsilon^j \Lambda_0$ obtained as the kernel of the
homomorphism
\[
\epsilon^{j-l}\Lambda_0/\epsilon^j \Lambda_0 \to
\epsilon^{-l}\Lambda'/\Lambda'
\] induced by $f$.

Now we prove the key statement (6). Let $\mathcal L$ be the set of pairs
$(\Lambda,\Lambda')$ of lattices such that $f\Lambda \supset
\Lambda'$. For any pair $(\Lambda,\Lambda')\in \mathcal L$  we put
\[ d(\Lambda,\Lambda'):=\dim \Lambda \cap f^{-1}\Lambda' - \dim \Lambda'.
\] We must show that $d(\Lambda,\Lambda')=d(V,\Phi)$ for all
$(\Lambda,\Lambda') \in \mathcal L$. Since we chose $\Lambda_0$ so as to
satisfy (3), we see that $(\Lambda_0,\Lambda_0') \in \mathcal L$ and
$d(\Lambda_0,\Lambda_0')=d(V,\Phi)$.  Thus it is enough to show that
\[ d(\Lambda,\Lambda')=d(\Lambda_0,\Lambda_0').
\] In fact we claim that for any $j \le 0$ such that $\Lambda \subset
\epsilon^j \Lambda_0$ there is a chain of equalities
\[ d(\Lambda,\Lambda')=d(\epsilon^j\Lambda_0,\Lambda')
=d(\epsilon^j\Lambda_0,\Lambda'\cap
\Lambda_0')=d(\epsilon^j\Lambda_0,\Lambda_0')=d(\Lambda_0,\Lambda_0').
\] Each of these equalities is a consequence of one of the following two
statements.

\bf{Statement 1.} \rm If $(\Lambda,\Lambda') \in \mathcal L$ and
$\Lambda'' \subset \Lambda'$, then $(\Lambda,\Lambda'') \in \mathcal L$ and
$d(\Lambda,\Lambda')=d(\Lambda,\Lambda'')$.

\bf {Statement 2.} \rm If $(\Lambda,\Lambda') \in \mathcal L$ and
$\Lambda \subset \Lambda_1$, then $(\Lambda_1,\Lambda') \in \mathcal L$ and
$d(\Lambda,\Lambda')=d(\Lambda_1,\Lambda')$.

Statement 1 follows from the fact that
$f$ induces an isomorphism
\[
\frac{\Lambda \cap f^{-1}\Lambda'}{\Lambda \cap f^{-1}\Lambda''} \to
\frac{\Lambda'}{\Lambda''}.
\] Statement 2 follows from the fact that
\[
\frac{\Lambda_1 \cap \ker f}{\Lambda \cap \ker f}
 \simeq \frac{\Lambda_1 \cap f^{-1}\Lambda'}{\Lambda
\cap f^{-1}\Lambda'},
\] which implies (by part (2) of this lemma) that $\Lambda_1 \cap
f^{-1}\Lambda'$ is a finite union of cosets of $\Lambda
\cap f^{-1}\Lambda'$.

Part (7) of the lemma follows immediately from part (6).
\end{proof}

Now we come to the key proposition. Let $V$, $V'$ be finite dimensional
$L$-vector spaces of the same dimension. Let $\phi:V \to V'$ be an
$L$-linear bijection and let $\psi:V \to V'$ be a $\sigma$-linear bijection.
Put $f:=\psi-\phi$, an $F$-linear map $V \to V'$. Since
$(V,\phi^{-1}\psi)$ is an $F$-space, it makes sense to put
$d:=d(V,\phi^{-1}\psi)$.

In the proposition below we will be considering the kernel $\ker \bar f$ of a
certain homomorphism $\bar f$ of algebraic groups over $\bar k$. We take this
kernel in the naive sense, that is to say, we take the reduced subscheme
structure on the closed set $\ker\bar f$. We then write
$(\ker\bar f)^0$ for  the identity component of $\ker \bar f$.

\begin{proposition}\label{prop.MNkey} For any lattice $\Lambda$ in $V$ there
exists a lattice $\Lambda'$ in
$V'$ and a non-negative integer $j$ such that
\[
\epsilon^j \Lambda' \subset f\Lambda \subset \Lambda'.
\]
 For any such $\Lambda'$, $j$, and for any $l \ge j$ we consider the
homomorphism
\[
\bar f: \Lambda/\epsilon^l \Lambda \to \Lambda'/\epsilon^l \Lambda'
\] induced by $f$. Then
\begin{enumerate}
\item $\im \bar f \supset \epsilon^j \Lambda' / \epsilon^l \Lambda'$.
\item $\dim \ker \bar f= d + \dim \phi^{-1} \Lambda' - \dim \Lambda$.
\item $(\ker \bar f )^0 \subset \epsilon^{l-j}\Lambda/\epsilon ^l\Lambda$.
\end{enumerate}
\end{proposition}
\begin{proof} We may use $\phi$ to identify $V'$ with $V$. Thus it is enough
to prove this proposition when $V'=V$ and $\phi=\id_V$. Thus $(V,\psi)$ is an
$F$-space and $f$ now has the familiar form $\psi-\id$.

The existence of $\Lambda'$, $j$ is a consequence of part (4) of Lemma
\ref{lem.preVkey}. Part (1) of the lemma follows from our assumption that
$\epsilon^j \Lambda' \subset f\Lambda$. To prove (2) we note that
\[
\ker \bar f=\frac{\Lambda \cap f^{-1}(\epsilon^l \Lambda')}{\epsilon^l
\Lambda},
\] a set having dimension
$d+\dim \Lambda'-\dim \Lambda$ by part (6) of Lemma \ref{lem.preVkey}, which
can be used since $f\Lambda \supset \epsilon^l\Lambda'$.

To prove (3) it is enough to show that
\[ (\ker \bar f) \cap (\epsilon^{l-j}\Lambda/\epsilon^l\Lambda)
\] has the same dimension as $\ker\bar f$. For this we note that
\[ (\ker \bar f) \cap (\epsilon^{l-j}\Lambda/\epsilon^l\Lambda)=
\frac{\epsilon^{l-j}\Lambda \cap f^{-1}(\epsilon^l \Lambda')}{\epsilon^l
\Lambda},
\] a set again having dimension
$d+\dim \Lambda'-\dim \Lambda$ by part (6) of Lemma \ref{lem.preVkey}, which
can be used since $f(\epsilon^{l-j}\Lambda) \supset
\epsilon^l\Lambda'$.
\end{proof}

\section{Reduction step for more general elements $b$}
\label{reduction.section} So far we have treated only elements $b \in A(L)$.
Now we will prove Rapoport's conjectural formula for
$\dim X^G_\mu(b)$ assuming its validity  for basic elements.

\subsection{Notation related to $P=MN$} We continue with $G \supset B=AU$ and
$W$ as in the introduction. In addition we consider a parabolic subgroup
$P=MN$ containing $B$. Here
$N$ is the unipotent radical of $P$, and $M$ is the unique Levi factor of
$P$ that contains $A$. Put $B_M:=B\cap M$ and $U_M:=U\cap M$, so that
$B_M=AU_M$ is a Borel subgroup of $M$ containing $A$.  We write $W_M$ for the
Weyl group of $A$ in $M$ and view $W_M$ as a subgroup of $W$.

We write $A_M$ for the identity component of the center of $M$. Then $A_M$ is
a subtorus of $A$, and $X_*(A_M)$ is a subgroup of
$X_*(A)$. Put $\mathfrak a_M:=X_*(A_M)_ \mathbb R$, a linear subspace of
$\mathfrak a$.

Let $R_N$ denote the set of roots of $A$ in $\Lie(N)$. We write
$\rho_N \in X^*(A)_\mathbb Q$ for the half-sum of the roots in $R_N$.  Thus
$\rho=\rho_M+\rho_N$, where $\rho$ (respectively, $\rho_M$) is the half-sum of
the roots of $A$ in $\Lie (U)$ (respectively, $\Lie(U_M)$).

For a coweight $\mu$ of $A$ we will need to distinguish between the two
relevant notions of dominance. We say that $\mu$ is \emph{$G$-dominant}
(respectively, \emph{$M$-dominant}) if $\langle \alpha,\mu \rangle \ge 0$ for
every root $\alpha$ of $A$ in $\Lie(U)$ (respectively,
$\Lie(U_M)$).

Similarly we need to distinguish between the two relevant partial orders on
$X_*(A)$.   When
$\mu-\nu$ is a non-negative integral linear combination of simple coroots for
$G$ (respectively, $M$), we   write $\nu \le \mu$ (respectively, $\nu
\underset{M}{\le} \mu$).

Recall from the introduction the canonical surjections $p_G:X_*(A)
\twoheadrightarrow \Lambda_G$ and $\eta_G:X^G=G(L)/K \to\Lambda_G$.
We apply these definitions to $M$ as well as
$G$, obtaining $p_M:X_*(A) \twoheadrightarrow \Lambda_M$ and
\[
\eta_M:X^M=M(L)/K_M \twoheadrightarrow \Lambda_M,
\]  where $K_M$ denotes $M(\mathfrak o)$. There is an obvious embedding $X^M
\hookrightarrow X^G$. In particular we view $x_0$ as the base-point for
$X^M$ as well as for $X^G$.

Note that the composition
$X_*(A_M)\hookrightarrow X_*(A) \twoheadrightarrow \Lambda_M$ is injective and
identifies $X_*(A_M)$ with a subgroup of finite index in
$\Lambda_M$, so that $\mathfrak a_M$ can also be identified with
$\Lambda_M \otimes_\mathbb Z \mathbb R$; in this way we obtain a canonical
homomorphism
\begin{equation}\label{eq.lmtoam}
\Lambda_M \to \mathfrak a_M.
\end{equation}

\subsection{Dimensions for admissible and  ind-admissible subsets of
$N(L)$}\label{sub.dimaisN} In \ref{sub.subU(m)} we defined subgroups
\[
\dots \supset U(-2) \supset U(-1) \supset U(0) \supset U(1)
\supset U(2) \supset
\dots
\] of $U(L)$, and in \ref{sub.dimaisU} we used these subgroups to define
(ind)-admissible subsets of $U(L)$, as well as dimensions of such sets.

For $m \in \mathbb Z$ we now put
\[ N(m):=N(L) \cap U(m)
\] and define (ind)-admissibility and dimensions for subsets of $N(L)$ just as
we did for $U(L)$, using the subgroups $N(m)$ in place of $U(m)$.

\subsection{Relative dimension of $f_m$} As before, for $\nu \in
\mathfrak a$ we write $\nu_{\dom}$ for the unique $G$-dominant element in the
$W$-orbit of $\nu$.
\begin{proposition}\label{prop.reldimfP}
  Let
$m \in M(L)$ and let
$\nu$ be the $M$-dominant element of
$\mathfrak a$ defined as  the Newton point of the $\sigma$-conjugacy class of
$m$ in $M(L)$.  Define a map
$f_m:N(L)
\to N(L)$ by $f_m(n)=n^{-1}m\sigma(n)m^{-1}$. Let $Y$ be an admissible subset
of $N(L)$. Then $f_m^{-1}Y$ is ind-admissible and
\begin{align}\label{eq.reldimN}
\dim f_m^{-1}Y-\dim Y&=\sum_{\alpha \in R_N}\min\{\langle \alpha,\nu
\rangle,0\} \\ &=\langle \rho,\nu-\nu_{\dom} \rangle.
\end{align}   Moreover $f_m$ is surjective.
\end{proposition}
\begin{proof} It would be awkward to prove this proposition as it stands,
since conjugation by
$m$ need not preserve $N(\mathfrak o)$. Therefore we are going to reformulate
the proposition in a way that will make it easier to prove.
\end{proof}

For $m_1,m_2 \in M(L)$ we now define a map
$ f_{m_1,m_2}:N(L) \to N(L)
$
 by
\[ f_{m_1,m_2}(n):=m_1n^{-1}m_1^{-1}\cdot m_2\sigma(n)m_2^{-1},
\] which gives back $f_m$ when $(m_1,m_2)=(1,m)$. Put $\mathfrak n:=\Lie N$.
For any
$m \in M(L)$ we denote by $\Ad_\mathfrak n(m)$ the adjoint action of $m$ on
$\mathfrak n(L)$. Now $\Ad_\mathfrak n(m_1)^{-1}\Ad_\mathfrak n(m_2)\sigma$ is
a $\sigma$-linear bijection from $\mathfrak n(L)$ to itself, so
$(\mathfrak n(L),\Ad_\mathfrak n(m_1)^{-1}\Ad_\mathfrak n(m_2)\sigma)$ is an
$F$-space, and we may put
\[ d(m_1,m_2):=d\bigl(\mathfrak n(L),\Ad_\mathfrak n(m_1)^{-1}\Ad_\mathfrak
n(m_2)\sigma\bigr) + \val \det \Ad_\mathfrak n(m_1).
\] Note that
\begin{align*} d(1,m)&=d(\mathfrak n(L),\Ad_\mathfrak n(m)\sigma)\\
&=\sum_{\alpha \in R_N} \min \{ \langle \alpha,\nu \rangle,0 \}
\end{align*} since by its very definition the Newton point $\nu$ of $m$ has
the property that the slopes of the $F$-space $(\mathfrak n(L),\Ad_\mathfrak
n(m)\sigma)$ are the numbers $\langle \alpha,\nu \rangle$ for $\alpha  \in
R_N$. Moreover it follows from Lemma \ref{lem.rho.dom} that
\[
\sum_{\alpha \in R_N} \min \{ \langle \alpha,\nu \rangle,0 \}=\langle
\rho,\nu-\nu_{\dom} \rangle,
\] since $\langle \alpha,\nu \rangle \ge 0$ for all positive roots $\alpha$ in
$M$.   Therefore the next proposition does indeed generalize the previous one,
and it will be enough to prove it.

\begin{proposition} The map $f_{m_1,m_2}$ is surjective. Moreover, for any
admissible subset $Y$ of $N(L)$ the inverse image $f^{-1}_{m_1,m_2} Y$ is
ind-admissible and
\[
\dim f^{-1}_{m_1,m_2} Y - \dim Y=d(m_1,m_2).
\]
\end{proposition}
\begin{proof} For any $m \in M(L)$ we have
\[ f_{mm_1,mm_2}=c_m \circ f_{m_1,m_2},
\] where $c_m:N(L) \to N(L)$ denotes the conjugation map $n \mapsto
mnm^{-1}$.  It is easy to see that for any admissible subset $Y$ of $N(L)$ the
inverse image
$c_m^{-1}Y$ is admissible, and that
\begin{equation}\label{eq.cmshift}
\dim c_m^{-1} Y -\dim Y=\val \det \Ad_\mathfrak n(m).
\end{equation}  Therefore the proposition is true for $m_1,m_2$ if and only
if it is true for $mm_1,mm_2$. So this proposition is not more general than
the previous one in any serious way, but it is more convenient to prove, since
we are free to improve the given pair $(m_1,m_2)$ as follows.

Let $M(L)_+$ denote the union of all $K_M$-double cosets
$K_M\epsilon^{\mu_M}K_M$ as $\mu_M$ ranges over all $M$-dominant coweights
such that $\langle \alpha,\mu_M \rangle \ge 0$ for all $\alpha \in R_N$.

What is the point of this definition? For any algebraic group $I$ over $k$ and
any $j \ge 0$ the kernel of the homomorphism $I(\mathfrak o)
\twoheadrightarrow I(\mathfrak o/\epsilon^j \mathfrak o)$ is a normal subgroup
of $I(\mathfrak o)$ that we will denote by $I_j$. Thus
\[ I(\mathfrak o)=I_0 \supset I_1 \supset I_2 \supset \dots
\] with $I_0/I_j=I(\mathfrak o/\epsilon^j \mathfrak o)$. Applying this
definition to $N$ (later we'll apply it to some other groups), we get normal
subgroups $N_j$ of $N(\mathfrak o)$, and the point is that for any $m \in
M(L)_+$ and $j \ge 0$
\[ mN_jm^{-1} \subset N_j.
\]

There exists $a \in A_M(L)$ such that $am_1,am_2 \in M(L)_+$. Replacing
$(m_1,m_2)$ by $(am_1,am_2)$, we see that it is enough to prove the
proposition under the additional assumption that $m_1,m_2 \in M(L)_+$, an
assumption we will make from now on.

Define homomorphisms $\phi,\psi:N(L) \to N(L)$ by $\phi(n):=m_1nm_1^{-1}$ and
$\psi(n):=m_2\sigma(n)m_2^{-1}$. Note that $\phi(N_j) \subset N_j$ and
$\psi(N_j) \subset N_j$ for all $j \ge 0$. Thus we may define a right action
of $N(\mathfrak o)$ on itself, called the $*$-action, with $n \in N(\mathfrak
o)$ acting on $n' \in N(\mathfrak o)$ by
\[ n'*n:=\phi(n)^{-1}n'\psi(n).
\]

The homomorphisms
\[
\phi_j,\psi_j:N_0/N_j \to N_0/N_j
\] induced by $\phi$, $\psi$ respectively allow us to define a right
$*$-action of $N_0/N_j$ on itself by
\[ n'*n:=\phi_j(n)^{-1}n'\psi_j(n).
\] Denoting by $n \mapsto \bar n$ the canonical surjection $N_0
\twoheadrightarrow N_0/N_j$, we have the obvious compatibility
\[
\overline{n'*n}=\bar n' * \bar n.
\]

We need to linearize our problem by means of a suitable filtration
\[ N=N[1] \supset N[2] \supset \dots
\] of $N$ by normal subgroups of $P$ such that each successive quotient
\[ N\langle i\rangle:=N[i]/N[i+1]
\] is abelian. To this end we introduce the following coweight $\delta_N$ of
the quotient $A/Z$ of $A$ by the center $Z$ of $G$. In $X_*(A/Z)$ we have the
$\mathbb Z$-basis of fundamental coweights $\varpi_\alpha$, one for each
simple root $\alpha$ of $A$. Define the coweight $\delta_N$ to be the sum of
 the fundamental coweights $\varpi_\alpha$ with $\alpha$ ranging through the
simple roots of $A$ that occur in $\mathfrak n=\Lie N$. Then for $i \ge 1$ we
let $N[i]$ be the product of all root subgroups $U_\alpha$ for which
$\langle \alpha, \delta_N \rangle \ge i$.

Put $\mathfrak n[i]:=\Lie N[i]$ and $\mathfrak n\langle i \rangle:=\Lie
N\langle i \rangle $. Then there is a $P$-equivariant isomorphism
\begin{equation}\label{eq.eq.iso} N\langle i \rangle \cong \mathfrak n\langle
i \rangle
\end{equation} of algebraic groups over $k$, and the action of $P$ on both
groups factors through $P \twoheadrightarrow M$. (To construct this
isomorphism first work with split groups over $\mathbb Z$.)

The homomorphisms $\phi$, $\psi$ preserve the subgroups $N[i](L)$ and
therefore induce homomorphisms
\[
\phi\langle i \rangle,\psi \langle i \rangle:N\langle i \rangle(L) \to
N\langle i \rangle(L)
\] as well, so we get $*$-actions of $N[i](\mathfrak o)$ and $N\langle i
\rangle (\mathfrak o)$ on themselves. Moreover, under the identification
\eqref{eq.eq.iso} of $N\langle i \rangle(L)$ with $\mathfrak n\langle i
\rangle (L)$, the homomorphism $\phi\langle i \rangle$ (respectively,
$\psi\langle i \rangle $) goes over to the $L$-linear bijection
$\Ad_{\mathfrak n \langle i \rangle}(m_1)$ (respectively, the $\sigma$-linear
bijection $\Ad_{\mathfrak n\langle i \rangle}(m_2)\sigma$), where
$\Ad_{\mathfrak n\langle i \rangle}(\cdot)$ denotes the adjoint action of
$M$ on $\mathfrak n\langle i \rangle=\mathfrak n[i]/\mathfrak n[i+1]$.

Thus for each $i$ we are in the linear situation considered in Proposition
\ref{prop.MNkey} (with $V=V'=\mathfrak n\langle i \rangle(L)$), and as in that
proposition we define an $F$-linear map
\[ f\langle i \rangle :\mathfrak n\langle i \rangle (L) \to \mathfrak n\langle
i \rangle (L)
\] by $f\langle i \rangle:= \psi\langle i \rangle-\phi\langle i \rangle$.

Since $m_1,m_2 \in M(L)_+$, both $\phi\langle i \rangle$ and $\psi\langle i
\rangle$ carry the lattice $\mathfrak n\langle i \rangle(\mathfrak o)$ into
itself;  therefore $f\langle i \rangle$ does too. By Proposition
\ref{prop.MNkey} there exists $j \ge 0$ such that
\begin{equation}\label{eq.jincl} f\langle i \rangle\bigl(\mathfrak n\langle i
\rangle (\mathfrak o)\bigr)
\supset
\epsilon^j\bigl(\mathfrak n\langle i \rangle (\mathfrak o)\bigr).
\end{equation} Since $\mathfrak n\langle i \rangle=0$ for large $i$, we can
even choose $j$ so that \eqref{eq.jincl} holds for all $i$.

We also need to look at this slightly differently. For $a \in A_M(F)$ define
\[ c_a:N(L) \to N(L)
\]  by $c_a(n):=ana^{-1}$. The homomorphisms $c_a$ preserve the subgroups
$N[i](L)$ and hence induce homomorphisms
\[ c_a:N\langle i \rangle(L) \to N\langle i \rangle(L).
\]

For $a \in A_M(F)_+:=A_M(F) \cap M(L)_+$ we have $c_a N_0 \subset N_0$,
$c_a N[i]_0 \subset N[i]_0$, $c_a N\langle i \rangle_0 \subset N \langle i
\rangle_0$, and   it is easy to see that
\[ c_a(n'*n)=c_a(n')*c_a(n)
\]  for all $n,n' \in N_0$.

\bf Claim 1. \rm There exists $a \in A_M(F)_+$ such that for any $i \ge 1$
\begin{equation}\label{eq.ca.inc.0}  c_a N[i]_0 \subset  1*N[i]_0.
\end{equation}

\noindent Now we prove Claim 1. We begin by proving the weaker statement that
for each
$i \ge 1$ there exists $a_i \in A_M(F)_+$ such that
\eqref{eq.ca.inc.0} holds for $a_i$ and $i$, and this we prove  by descending
induction on
$i$, starting with any $i$
 such that $N[i]$ is trivial, so that the statement we need to prove becomes
trivial. Assume the statement is true for
$i+1$ and prove it for $i$.  Thus we are assuming that there exists $a_{i+1}
\in A_M(F)_+$ such that
\eqref{eq.ca.inc.0} holds for $a_{i+1}$ and $i+1$.

Applying Proposition
\ref{prop.MNkey} just as we did to find $j$, we see that there exists
$a' \in A_M(F)_+$ such that
\begin{equation}\label{eq.caincl} c_{a'} N\langle i \rangle _0
\subset 1*N\langle i \rangle _0.
\end{equation} We will now check that \eqref{eq.ca.inc.0} holds for
$a_i:=a'a_{i+1}$ and
$i$.

Indeed, consider an element $x=c_{a_i}n$ with $n \in N[i]_0$. The image of
$c_{a'}n$ under $N[i]_0 \twoheadrightarrow N\langle i \rangle_0$ lies in
$c_{a'} N\langle i \rangle_0$, and therefore \eqref{eq.caincl} guarantees the
existence of $n' \in N[i]_0$ such that $(c_{a'} n)*n' \in N[i+1]_0$, whence
\[ c_{a_{i+1}}\bigl ( (c_{a'}n)*n'\bigr) \in c_{a_{i+1}} N[i+1]_0.
\] Since  \eqref{eq.ca.inc.0} holds for $a_{i+1}$ and
$i+1$, there exists $n'' \in N[i+1]_0$ such that
\begin{equation}\label{eq.ca.mess}
\bigl (c_{a_{i+1}} ((c_{a'}n)*n')\bigr )*n''=1.
\end{equation}  The left side of equation \eqref{eq.ca.mess} works out to
\[ c_{a_i}n*(c_{a_{i+1}}(n')n''),
\]  showing that $c_{a_i}n \in 1*N[i]_0$, as desired. This proves the weaker
statement.

To prove the claim itself we first note that for $a,a' \in A_M(F)_+$ it is
clear that $aa' \in A_M(F)_+$ and that
\[ c_{aa'}N[i]_0=c_a(c_{a'}N[i]_0) \subset c_a N[i]_0.
\] Therefore if \eqref{eq.ca.inc.0} holds for $a$ and $i$,  it also holds for
$aa'$ and $i$. We  already know that for each of the finitely many values of
$i$ for which $N[i]$ is non-trivial, we can find an element $a_i$ for which
\eqref{eq.ca.inc.0} holds for $a_i$ and $i$, and taking the product of all
these elements $a_i$, we get an element $a$ such that \eqref{eq.ca.inc.0}
holds for all $i$.
 This proves the claim. \medskip

Now we fix $a \in A_M(F)_+$ as in Claim 1. Choose $j \ge 0$  as before, so that
\eqref{eq.jincl} holds for all $i$. For reasons that will soon become
apparent, we now consider any $l \ge j$ large enough that for all $i \ge 1$
\begin{equation}\label{eq.lja.in} N[i]_{l-j} \subset c_a N[i]_0.
\end{equation}

 We denote by $H$, $H[i]$, $H\langle i \rangle$ the groups of $\mathfrak
o/\epsilon^l \mathfrak o$-points of
$N$, $N[i]$, $N\langle i \rangle$ respectively,  regarded as algebraic groups
over $\bar k$. Note that $H=H[1]$ and $H\langle i
\rangle=H[i]/H[i+1]$. On $H[i]$, $H\langle i \rangle$ we have descending
filtrations $H[i]_\beta$, $H\langle i \rangle_\beta$ coming from the powers
$\epsilon^\beta$ of $\epsilon$; more precisely, for $0 \le \beta \le l$ we put
\[ H[i]_\beta:=\ker \bigr[N[i](\mathfrak o/\epsilon^l \mathfrak o) \to
N[i](\mathfrak o/\epsilon^\beta \mathfrak o)\bigr]
\] and similarly for $H\langle i \rangle$. The $*$-actions on
$N[i](\mathfrak o)$ and $N\langle i \rangle (\mathfrak o)$ induce compatible
$*$-actions on $H[i]$ and $H\langle i \rangle$.

The homomorphism $c_a$ preserves both $N_0$ and $N_l$ and hence induces a
homomorphism $c_a:H \to H$. Since Claim 1 holds for $a$, we see that for any
$i
\ge 1$
\begin{equation}\label{eq.44} c_a H[i] \subset 1*H[i].
\end{equation} Moreover \eqref{eq.lja.in} implies immediately that
\begin{equation}\label{eq.22} H[i]_{l-j} \subset c_a H[i].
\end{equation} Also note that
\begin{equation}\label{eq.33} H[i] \cap c_aH=c_a H[i],
\end{equation} as one sees easily from the description of $H$ as a product of
copies of
$\mathfrak o/\epsilon^l\mathfrak o$, one for each root in $R_N$.

Write $S[i]$ (respectively, $S\langle i \rangle$) for the stabilizer, for the
$*$-action, of $1 \in H[i]$ in $H[i]$ (respectively, of $1 \in  H\langle i
\rangle$ in $H\langle i \rangle$). Note that $S:=S[1]$ is the stabilizer, for
the $*$-action, of $1 \in H$ in $H$. Here we are taking naive stabilizers: put
the reduced subscheme structure on the set-theoretic stabilizers. \medskip

\bf Claim 2. \rm $\dim S=d(m_1,m_2)$.

\noindent We begin by proving that
\begin{equation}\label{eq.si+1}
\dim S[i]=\dim S[i+1] + \dim S\langle i \rangle.
\end{equation}

For this it is enough to prove that
\[ S\langle i\rangle^0 \subset \im[S[i] \to H\langle i \rangle] \subset
S\langle i \rangle.
\] The second inclusion being clear, it is enough to prove the first one. So
let $s \in S\langle i \rangle^0$. By Proposition \ref{prop.MNkey} we have
\[ S\langle i \rangle^0 =\bigl (\ker\overline{ f\langle i \rangle}\bigr)^0
\subset H\langle i \rangle_{l-j}.
\] Pick $h \in H[i]_{l-j}$ such that $h \mapsto s$. By \eqref{eq.22} there
exists $h_0 \in H[i]$ such that $h=c_ah_0$. Now
$1*h \mapsto 1*s=1 \in H\langle i \rangle$, showing that $1*h \in H[i+1]$.
Moreover $1*h=c_a(1*h_0)$, showing that $1*h \in c_aH$. From
\eqref{eq.33} and \eqref{eq.44} we conclude that
\[ 1*h \in c_aH[i+1] \subset 1*H[i+1]
\]  which means that there exists $h_1 \in H[i+1]$ such that $1*hh_1=1$. Thus
$hh_1 \in S[i]$, and clearly $hh_1 \mapsto s$, proving that $s$ lies in the
image of $S[i]$, as desired. This proves the equality \eqref{eq.si+1}.

It follows from \eqref{eq.si+1} that
\[
\dim S =\sum_{i \ge 1} \dim S\langle i \rangle.
\] Moreover, by part (2) of Proposition \ref{prop.MNkey}
\[
\dim S\langle i \rangle=d(\mathfrak n\langle i \rangle(L),\phi\langle i
\rangle^{-1}\psi\langle i \rangle)+\dim \phi\langle i
\rangle^{-1}(\mathfrak n\langle i\rangle(\mathfrak o))-\dim \mathfrak n\langle
i \rangle(\mathfrak o).
\] Therefore
\begin{align*}
\dim S&=d\bigl(\mathfrak n(L),\Ad_\mathfrak n(m_1)^{-1}\Ad_{\mathfrak
n}(m_2)\sigma\bigr)+\dim
\Ad_\mathfrak n(m_1)^{-1}\mathfrak n(\mathfrak o)-\dim
\mathfrak n(\mathfrak o)\\ &=d(m_1,m_2),
\end{align*} completing the proof of Claim 2. \medskip

Now consider the map $f_0:N(\mathfrak o) \to N(\mathfrak o)$ defined by
$f_0(n):=1*n$. Equivalently, $f_0$ is the restriction of $f_{m_1,m_2}$ to
$N(\mathfrak o)$.

\bf Claim 3. \rm Suppose that $Y$ is an admissible subset of $c_a N_0 $ (with
$a$ again as in Claim 1). Then $f_0^{-1}Y$ is admissible and
\[
\dim f_0^{-1} Y -\dim Y=d(m_1,m_2).
\]

\medskip
\noindent Now we prove Claim 3. The set $Y$ comes from a locally closed subset
$\bar Y$ of $H=N(\mathfrak o/\epsilon^l\mathfrak o)$ for some suitably large
$l$.  By increasing $l$ we may assume that Claim 2 holds for $l$.
 We have a commutative diagram
\[
\begin{CD} N(\mathfrak o) @>f_0>> N(\mathfrak o) \\ @VpVV @VpVV \\ H @>\bar
f>> H
\end{CD}
\] where $p$ is the canonical surjection and $\bar f$ is induced by $f_0$. By
Claim 1 $\bar Y$ is contained in the image of $\bar f$. Every non-empty
(reduced) fiber of $\bar f$ is isomorphic to the stabilizer group
$S$ considered above. Therefore
\begin{align*}
\dim \bar f^{-1} \bar Y &=\dim \bar Y +\dim S \\ &=\dim \bar Y +d(m_1,m_2).
\end{align*} Since $f_0^{-1} Y=p^{-1}\bar f^{-1}\bar Y$, we see that
$f_0^{-1} Y$ is admissible and
\[
\dim f_0^{-1} Y -\dim Y=d(m_1,m_2),
\] proving Claim 3.

Now we finish the proof. We abbreviate $f_{m_1,m_2}$ to $f$ and $d(m_1,m_2)$
to $d$. Let $Y$ be any admissible subset of $N(L)$.

For any $a' \in A_M(F)$ it is clear that $f$ commutes with $c_{a'}$. Thus it
follows from Claim 3 together with \eqref{eq.cmshift}
 that $f^{-1}Y \cap a'^{-1}N(\mathfrak o)a'$ is admissible of dimension
$\dim Y + d$ for all $a' \in A_M(F)$ such that
$a'Ya'^{-1}
\subset aN_0 a^{-1}$ (with $a$ as before, so that Claims 1 and 3 hold for it).

Pick a coweight $\lambda_0 \in X_*(A_M)$ such that $\langle \alpha,
\lambda_0 \rangle > 0$ for all $\alpha \in R_N$, and for $t \in \mathbb Z$ put
$a_t:=\lambda_0(\epsilon^t)$. There exists $t_0 \ge 0$ such that
$a_tYa_t^{-1} \subset aN_0 a^{-1}$ for all $t \ge t_0$. Therefore $f^{-1}Y
\cap a_t^{-1}N(\mathfrak o)a_t$ is admissible of dimension $\dim Y + d$ for all
$t \ge t_0$, and since any of the subgroups $N(i)$ ($i \in \mathbb Z$) used in
\ref{sub.dimaisN} to define dimensions is contained in
$a_t^{-1}N(\mathfrak o)a_t$ for sufficiently large $t$, we see that
$f^{-1}Y$ is ind-admissible and that
\[
\dim f^{-1}Y -\dim Y =d,
\] as desired.

The last point is the surjectivity of $f$. We already know by Claim 1 that
 $f(N_0)$ contains $aN_0 a^{-1}$. Therefore
$f(a_t^{-1}N_0 a_t)$ contains $a_t^{-1} aN_0 a^{-1} a_t$ for all $t \in
\mathbb Z$. Since the union of the subgroups $a_t^{-1} aN_0 a^{-1} a_t$ ($t
\in
\mathbb Z$)  is all of
$N(L)$, we see that $f$ is indeed surjective.
\end{proof}

\subsection{Dimensions of intersections of $N(L)$- and $K$- orbits on
$X^G$} As before we put $x_\lambda=\epsilon^\lambda x_0$. Let $\mu$ be a
$G$-dominant coweight.

For $m \in M(L)$ we are interested in the intersection
\begin{equation}\label{eq.int.KN} N(L)mx_0
\cap Kx_\mu
\end{equation}  of the
$N(L)$-orbit of $mx_0$ and the $K$-orbit of $x_\mu$ in the affine Grassmannian
$X^G$. Let $\mu_M$ be the unique $M$-dominant coweight such that $m \in
K_M\epsilon^{\mu_M} K_M$. Then there exists $k_M \in K_M$ such that
$mx_0=k_Mx_{\mu_M}$, and left multiplication by $k_M$ defines an isomorphism
\begin{equation}\label{eq.iso.twoorb} N(L)x_{\mu_M} \cap Kx_\mu \simeq
N(L)mx_0 \cap Kx_\mu.
\end{equation}

 For a given $G$-dominant coweight $\mu$ we write $S_M(\mu)$ for the set of
$M$-dominant coweights $\mu_M$ for which the intersection
$N(L)x_{\mu_M}
\cap Kx_\mu$ is non-empty.  For $\nu \in \Lambda_M$ we then put
\[ S_M(\mu,\nu):=\{\mu_M \in S_M(\mu):p_M(\mu_M)=\nu \}.
\]

The next lemma will help us understand  the finite sets
$S_M(\mu)$ and $S_M(\mu,\nu)$.
 Define a finite subset
$\Sigma(\mu)$ of $X_*(A)$ as follows: $\mu' \in \Sigma(\mu)$ if and only if
$\mu'_{\dom} \le \mu$, where, as before, $\mu'_{\dom}$ denotes the unique
$G$-dominant element in the $W$-orbit of $\mu'$. Then denote by
$\Sdom$ the set of $M$-dominant elements in $\Sigma(\mu)$. Finally,  denote by
$\Smax$ the set of all elements  in $\Sdom$ that are maximal in $\Sdom$ with
respect to the partial order $\underset{M}{\le}$ on $X_*(A)$.

\begin{lemma}\label{lem.3subsets}  For any $G$-dominant coweight $\mu$ there
are inclusions
\[
\Smax \subset S_M(\mu) \subset \Sdom.
\]  Moreover, all three sets $\Smax$, $ S_M(\mu)$, $\Sdom$ have the same image
under the map $p_M:X_*(A) \twoheadrightarrow \Lambda_M$. In particular
$S_M(\mu,\nu)$ is non-empty if and only if $\nu$ lies in the image of
$\Sdom$ under $p_M$.
\end{lemma}
\begin{proof} First recall
\cite{Mat} (see also \cite{Rap00,H}) that a coweight $\mu'$ lies in
$\Sigma(\mu)$ if and only if $U(L)\epsilon^{\mu'}$ meets
$K\epsilon^\mu K$. By definition $\mu' \in S_M(\mu)$ if and only if
$\mu'$ is $M$-dominant and $N(L)\epsilon^{\mu'}$ meets $K\epsilon^\mu K$.
Since $N \subset U$, it is clear that $S_M(\mu) \subset
\Sdom$.

Now let $\mu' \in \Smax$. Thus there exists $u
\in U(L)$ such that $u\epsilon^{\mu'} \in K\epsilon^\mu K$. Write
$u=nu_M$ with  $n \in N(L)$, $u_M \in U_M(L)$. Then
$N(L)u_M\epsilon^{\mu'}$ meets
$K\epsilon^\mu K$, so there exists $\mu_M \in S_M(\mu)$ such that
$u_M\epsilon^{\mu'} \in K_M\epsilon^{\mu_M} K_M$. Now applying the fact
recalled at the beginning of the proof to $M$ rather than $G$, we see that
$\mu' \underset{M}{\le} \mu_M$, and of course $\mu_M$ lies in
$\Sdom$, since we have already proved that $S_M(\mu) \subset
\Sdom$.  By maximality of
$\mu'$ we then conclude that $\mu'=\mu_M$, proving that $\mu' \in S_M(\mu)$,
as desired.

For the second statement of the lemma it is enough to show that
$\Smax$ and $\Sdom$ have the same image under the map $p_M:X_*(A)
\twoheadrightarrow \Lambda_M$. So suppose that  $p_M(\mu)=\nu$ for some $\mu
\in \Sdom$. Clearly there exists $\mu' \in \Smax$ such that $\mu
\underset{M}{\le}\mu'$. Since $\mu \underset{M}{\le}\mu'$ implies
$p_M(\mu')=p_M(\mu)$, we are done.

The last statement of the lemma follows from the second statement.
\end{proof}

For $G$-dominant $\mu$ and $\mu_M \in S_M(\mu)$ we denote by
$d(\mu,\mu_M)$ the dimension of the intersection $N(L)x_{\mu_M}
\cap Kx_\mu$. From \eqref{eq.iso.twoorb} we see that
\begin{equation}
\dim N(L)mx_0 \cap Kx_\mu=d(\mu,\mu_M)
\end{equation} for all $m \in K_M\epsilon^{\mu_M} K_M$.

Our next task is to give an estimate for the numbers $d(\mu,\mu_M)$.  At this
point we need to introduce some new notation.  Recall that we have fixed a
Borel $B = AU$ in $G$, and this gives rise to a based root system.  Let
$G^\vee \supset B^\vee \supset A^\vee$ denote the complex reductive group
together with a Borel and a maximal torus, whose root system is dual to that
determined by $G \supset B \supset A$.  The Levi subgroup $M$ is determined by
a subset $\Delta_M$ of the simple
$B$-positive roots for $G$.  We let $M^\vee$ denote the Levi subgroup of
$G^\vee$ determined by the set of simple $B^\vee$-positive roots $\{
\alpha^\vee ~ | ~ \alpha \in \Delta_M \}$.  Then $M^\vee$ is a dual group for
$M$.  Given a $G$-dominant coweight $\mu \in X_*(A)$, we will simultaneously
think of it also as a weight $\mu \in X^*(A^\vee)$.  Let
$\chi^G_\mu$ denote the character of the unique irreducible
$G^\vee$-module with highest weight $\mu$.  For an $M$-dominant coweight
$\mu_M \in X_*(A)$, the symbol $\chi^M_{\mu_M}$ has the analogous meaning.

\begin{proposition} \label{est.dim.inters}  Let $\mu$ be a $G$-dominant
coweight.  Then for all $\mu_M \in S_M(\mu)$ there is an inequality
\begin{equation} d(\mu,\mu_M) \le \langle \rho,\mu+\mu_M \rangle -2\langle
\rho_M,\mu_M \rangle,
\end{equation} and equality holds if and only if $\chi^M_{\mu_M}$ occurs in
the restriction of $\chi^G_\mu$ to $M^\vee$.  Moreover, the number of
irreducible components of the intersection $N(L)x_{\mu_M} \cap Kx_\mu$ having
dimension
$\langle \rho, \mu + \mu_M \rangle - 2\langle \rho_M, \mu_M \rangle$ is the
multiplicity $a_{\mu_M \, \mu}$ with which $\chi^M_{\mu_M}$ occurs in
$\chi^G_\mu$.
\end{proposition}

\begin{proof} To simplify notation, during this proof we temporarily use the
symbol $\lambda$ instead of $\mu_M$ to denote an $M$-dominant coweight.

Recall that $F = k((\epsilon))$ where $k$ denotes the finite field with
$q$ elements.  We will use the symbols $G,M,B,A$ etc. to denote the groups of
$F$-points $G(F)$, $M(F)$, $B(F)$, $A(F)$ etc., and furthermore  here we write
$K$ (resp. $K_M$) for the compact open subgroup
$G(k[[\epsilon]])$ (resp. $M(k[[\epsilon]])$ of $G(F)$ (resp. $M(F)$).  The
intersection $N(F)x_\lambda \cap Kx_\mu$ is the set of $k$-points of the
variety $N(L)x_\lambda \cap G(\mathfrak o)x_\mu$, and an estimate for the
dimension of the latter will follow from a suitable estimate for the growth of
the number of points
$$ n(\lambda,\mu)(q) = \#(N(F)x_\lambda \cap Kx_\mu)
$$ as a function of $q$.  In fact we will show that $n(\lambda, \mu)(q)$ is a
polynomial in $q$ with degree bounded above by $\langle \rho, \mu +
\lambda \rangle - 2\langle \rho_M, \lambda \rangle$.  This is enough to prove
the upper bound on the dimension of $N(L)x_{\lambda} \cap G(\mathfrak o)x_\mu$.

To prove this we will calculate $n(\lambda,\mu)(q)$ in terms of values of
constant terms of standard spherical functions, and we will estimate those
values by manipulating the Kato-Lusztig formula \cite{lusztig83,kato}.  An
exposition   of this key ingredient, as well as the Satake isomorphism used
below, can be found in
\cite{HKP}.

Let $H_K(G) = C_c(K \backslash G/K)$ and $H_K(M) = C_c(K_M \backslash M/K_M)$
denote the spherical Hecke algebras of $G$ and $M$ respectively.  Convolution
is defined using the Haar measures giving $K$ respectively
$K_M$ volume 1.  For a parabolic subgroup $P = MN$ of $G$, the constant term
homomorphism $c^G_M: H_K(G) \rightarrow H_K(M)$ is defined by the formula
$$ c^G_M(f)(m) = \delta_P(m)^{-1/2} \int_N f(nm) dn.
$$ Here $\delta_P(m) := |{\rm det}({\rm Ad}(m); {\rm Lie}(N))|$, and the Haar
measure on $N$ is such that $N \cap K$ has volume 1.  We define in a similar
way $\delta_B$,
$\delta_{B_M}$, $c^G_A$, and $c^M_A$.  Recall that $U = U_M \, N$, and so
$$
\delta_B(t) = \delta_P(t) \delta_{B_M}(t)
$$ for $t \in A$, and
$$ c^G_A(f)(t) = (c^M_A \circ c^G_M)(f)(t).
$$

Also, $c^G_A$ (resp. $c^M_A$) is the Satake isomorphism $S^G$ for $G$ (resp.
$S^M$ for $M$).  Thus, the following diagram commutes:
$$
\xymatrix{ R(G^\vee) \ar[r]^{\cong} \ar[d]_{rest.} & \mathbb C[X_*(A)]^W
\ar[d]_{incl.} & H_K(G) \ar[l]_{\,\,\,\,\,\,\,\,\,\, S^G}
\ar[d]_{c^G_M} \\ R(M^\vee) \ar[r]^{\cong} & \mathbb C[X_*(A)]^{W_M} & H_K(M)
\ar[l]_{\,\,\,\,\,\,\,\,\,\, S^M}, }
$$ where $R(G^\vee)$ denotes the representation ring for $G^\vee$.
\medskip

For a $G$-dominant coweight $\mu$, let $f^G_\mu = {\rm char}(K \mu K)$, and
let $f^M_\lambda$ have the analogous meaning.  For any $M$-dominant coweight
$\lambda$ of $M$, we define numbers $a_{\lambda \mu}$ and
$b_{\lambda \mu}(q)$ by the equalities
\begin{align*} rest.(\chi^G_\mu) &= \sum_{\lambda} a_{\lambda \mu} ~
\chi^M_\lambda,\\ c^G_M(f^G_\mu) &= \sum_{\lambda} b_{\lambda \mu}(q)
~f^M_\lambda.
\end{align*}

\begin{lemma} \label{est.const.term}  Assume that $\lambda \in S_M(\mu)$.

\begin{enumerate}
\item [(1)] We have $b_{\lambda \mu}(q) \in \Z[q^{1/2},q^{-1/2}]$, and
$q^{-\langle \rho, \mu \rangle + \langle \rho_M, \lambda \rangle} b_{\lambda
\mu}(q) \in \Z[q^{-1}]$.
\item[(2)] There is an equality
\begin{equation*} b_{\lambda \mu}(q) = a_{\lambda \mu} \, q^{\langle \rho,
\mu \rangle - \langle \rho_M, \lambda \rangle} + \{ \mbox{terms with lower
$q^{1/2}$-degree} \}.
\end{equation*}
\end{enumerate}
\end{lemma}

\begin{proof} Recall the Kato-Lusztig formula
\begin{equation*} \chi^G_{\mu} = \sum_{\mu' \le \mu} q^{-\langle \rho, \mu
\rangle} P_{w_{\mu'}, w_\mu}(q) \,\, S^G(f^G_{\mu'}),
\end{equation*} where $\mu'$ ranges over $G$-dominant elements in
$\Sigma(\mu)$, $t_{\mu'} \in \widetilde{W}$ denotes the corresponding
translation element, and $w_{\mu'}$ denotes the longest element in the  double
coset $Wt_{\mu'}W$.  Of course $P_{w_{\mu'},w_\mu}(q)$ is the Kazhdan-Lusztig
polynomial, and hence its $q$-degree is strictly bounded above by $\langle
\rho, \mu - \mu' \rangle$ if $\mu \neq \mu'$.   Also,
$P_{w_{\mu},w_{\mu}}(q) = 1$.

The Kato-Lusztig formula for the group $M$ can be written
$$
\chi^M_{\lambda'} = \sum_{\lambda} p^M_{\lambda \lambda'}(q^{-1}) ~
\dfrac{S^M(f^M_{\lambda})}{q^{\langle \rho_M, \lambda \rangle }},
$$  for some polynomials $p^M_{\lambda, \lambda'} \in \Z[q^{-1}]$ which vanish
unless $\lambda \underset{M}{\leq} \lambda'$.  This expresses the relation between two bases
of $R(M^\vee)$ given by an ``upper triangular'' invertible matrix.  Inverting
the corresponding relation for the group $G$ gives the formula
$$
\dfrac{S^G(f^G_\mu)}{q^{\langle \rho, \mu \rangle}} = \sum_{\mu'} q^G_{\mu'
\mu}(q^{-1}) ~ \chi^G_{\mu'},
$$ for some $q^G \in \Z[q^{-1}]$.   We find that
\begin{align*}
\dfrac{S(f^G_\mu)}{q^{\langle \rho, \mu \rangle}} &= \sum_{\mu'} q^G_{\mu'
\mu} \, \chi^G_{\mu'} \\ &= \sum_{\mu'} \sum_{\lambda'} a_{\lambda' \mu'}
\, q^G_{\mu' \mu} \, \chi^M_{\lambda'} \\ &= \sum_{\mu', \lambda',
\lambda} a_{\lambda' \mu'} \, q^G_{\mu' \mu} \, p^M_{\lambda \lambda'} \,
\dfrac{S(f^M_\lambda)}{q^{\langle \rho_M, \lambda \rangle}}.
\end{align*}

We also have the relation
$$
\dfrac{S(f^G_\mu)}{q^{\langle \rho, \mu \rangle}} =  \sum_\lambda
\dfrac{b_{\lambda \mu}(q)}{q^{\langle \rho, \mu \rangle - \langle \rho_M,
\lambda \rangle}} ~ \dfrac{S(f^M_{\lambda})}{q^{\langle \rho_M, \lambda
\rangle }}.
$$ Thus we see that for each $\lambda \in S_M(\mu)$,
$$
\dfrac{b_{\lambda \mu}(q)}{q^{\langle \rho, \mu \rangle - \langle \rho_M,
\lambda \rangle}} = \sum_{\mu', \lambda'} a_{\lambda' \mu'} \, q^G_{\mu'
\mu} \, p^M_{\lambda \lambda'}.
$$ Part (1) follows from this equality, as does (2) once one notes that
$q^G_{\mu', \mu}$ and $p^M_{\lambda, \lambda'}$ as polynomials in $\mathbb
Z[q^{-1}]$ have {\em non-zero} constant terms if and only if $\mu' = \mu$ and
$\lambda = \lambda'$ respectively.
\end{proof}

Now we continue with the proof of the proposition.  Using the integration
formula from the Iwasawa decomposition $G= KNM$  ($y = knm$), we see that the
number of points
$$ n(\lambda,\mu)(q) = \#(N \epsilon^{\lambda} K/K \cap K \epsilon^{\mu} K/K)
$$ is given by

\begin{align*} n(\lambda,\mu)(q) &=
\int_G 1_{NK}(\epsilon^{-\lambda} y^{-1}) ~ 1_{K \epsilon^{-\mu} K}(y) ~ dy
\\
 &= \int_M \int_N \int_K 1_{NK}(\epsilon^{-\lambda} m^{-1}n^{-1}k^{-1}) \,
1_{K\epsilon^{-\mu} K}(knm) \,  dk \, dn \, dm \\
  &= \int_M \int_N 1_{NK}(m^{-1}) \,  1_{K \epsilon^{-\mu} K}(n m
\epsilon^{-\lambda}) \,  dn \, dm \\  &= \int_N 1_{K
\epsilon^{-\mu} K}(n \epsilon^{-\lambda} ) \, dn \\ &=
\delta_P(\epsilon^{-\lambda})^{1/2} \, c^G_M(1_{K \epsilon^{-w^G_0 \mu}
K})(\epsilon^{-\lambda}) \\ &= q^{\langle \rho_N, \lambda \rangle} \,
c^G_M(1_{K
\epsilon^{-w^G_0\mu}K})(\epsilon^{-w^M_0 \lambda}).
\end{align*}  Here $\rho_N$ is the half-sum of the positive roots occurring in
${\rm Lie}(N)$, and $w^G_0$ (resp. $w^M_0$) is the longest element in the Weyl
group $W$ (resp. $W_M$).

Using Lemma \ref{est.const.term}, we see that
\begin{align*} {\rm dim}(N(L)x_\lambda \cap G(\mathfrak o)x_\mu) &\leq
\langle \rho_N,\lambda \rangle + \langle \rho, -w^G_0 \mu
\rangle - \langle \rho_M, -w^M_0 \lambda \rangle \\ &= \langle \rho -
\rho_M, \lambda \rangle + \langle \rho, \mu \rangle - \langle \rho_M,
\lambda \rangle \\ &= \langle \rho, \lambda + \mu \rangle - 2\langle
\rho_M, \lambda \rangle,
\end{align*} and equality holds if and only if $\chi^M_\lambda$ is contained
in $\chi^G_\mu$.  Moreover, in that case, the number of irreducible components
of the intersection having top dimension is the multiplicity
$a_{\lambda \mu}$.
\end{proof}

\begin{corollary}\label{cor.dim.inters} Let $\mu$ be a $G$-dominant coweight.
 Then for all  $\mu_M \in S_M(\mu)$ there is an inequality
\begin{equation} d(\mu,\mu_M) \le \langle \rho,\mu+\mu_M \rangle -2\langle
\rho_M,\mu_M \rangle,
\end{equation} and when $\mu_M \in \Smax$, this inequality is an equality.
\end{corollary}
\begin{proof} If $\mu_M \in \Smax$, then it is easy to see that
$\chi^M_{\mu_M}$ appears in the restriction of $\chi^G_{\mu}$ to $M^\vee$.
\end{proof}

\subsection{Remark} \label{alt.mirk.vil} Taking $M = A$ in Proposition
\ref{est.dim.inters}, we recover the dimension formula of Mirkovi\'c-Vilonen,
Proposition
\ref{orig.mirk.vil}.

\subsection{Reduction step} Let $b \in M(L)$ and assume that $b$ is basic in
$M(L)$.  The map $\eta_M:M(L) \to \Lambda_M$  is constant on
$\sigma$-conjugacy classes and induces a bijection between the set of basic
$\sigma$-conjugacy classes in $M(L)$ and the group $\Lambda_M$. Put
$\nu:=\eta_M(b) \in \Lambda_M$. We denote by $\bar\nu \in \mathfrak a_M$ the
image of $\nu$ under the canonical homomorphism $\Lambda_M \to
\mathfrak a_M$ (see \eqref{eq.lmtoam}); since $b$ is basic, $\bar\nu$ is the
Newton point of
$b$.

For any
$G$-dominant coweight $\mu$ we have the affine Deligne-Lusztig variety
\[ X^G_\mu(b)=\{x \in G(L)/K:x^{-1}b\sigma(x) \in K\epsilon^\mu K\},
\]  which may or may not be non-empty.  Moreover, for any $M$-dominant
coweight $\mu_M$ we have the  affine Deligne-Lusztig variety
\[ X^M_{\mu_M}(b)=\{x \in M(L)/K_M:x^{-1}b\sigma(x) \in K_M\epsilon^{\mu_M}
K_M\}.
\]
 Since $b$ is basic in $M(L)$
 the affine Deligne-Lusztig variety
$X^M_{\mu_M}(b)$ is non-empty (see \cite{kottwitz-rapoport02}) if and only if
$p_M(\mu_M)=\nu$.

Our goal is to express $\dim X^G_\mu(b)$ in terms of the numbers
$\dim X^M_{\mu_M}(b)$ and $d(\mu,\mu_M)$ for the various coweights $\mu_M
\in S_M(\mu,\nu)$, and then to use this expression to reduce Rapoport's
dimension conjecture to the basic case.

We will use the following map $\alpha:X^G \to X^M$. By the Iwasawa
decomposition we have $G(L)=P(L)K$. Therefore there is an obvious bijection
\begin{equation*} P(L)/K_P \simeq G(L)/K,
\end{equation*}  where $K_P:=P(\mathfrak o)$.  The canonical surjective
homomorphism $P(L) \twoheadrightarrow M(L)$ induces a surjection
\begin{equation*} P(L)/K_P \twoheadrightarrow M(L)/K_M.
\end{equation*} The map $\alpha$ is defined as the composition
\[
\alpha:X^G= G(L)/K \simeq P(L)/K_P \twoheadrightarrow M(L)/K_M=X^M.
\]  We hasten to warn the reader that $\alpha$ is not a morphism of
ind-schemes. However, for any connected component $Y$ of $X^M$, the inverse
image $\alpha^{-1}Y$ is a locally closed subset of $X^G$, and the map
$\alpha^{-1}Y \to Y$ obtained by restriction from $\alpha$ is a morphism of
ind-schemes.
 For $m \in M(L)$ we can identify the fiber $\alpha^{-1}(mx_0)$ with
$N(L)/N(\mathfrak o)$ via
\[ N(L)/N(\mathfrak o) \ni n \mapsto mnx_0 \in \alpha^{-1}(mx_0).
\]

\begin{proposition}\label{prop.red.basic}  Assume that the
$\sigma$-conjugacy class of $b$ in
$M(L)$ is basic, and let $\nu$ be the element $\nu = \eta_M(b) \in
\Lambda_M$.

\begin{enumerate}
\item The space $X^G_\mu(b)$ is non-empty if and only if $\nu$ lies in the
image of the set  $\Sdom$.
\item The image of $X^G_\mu(b)$ under $\alpha$ is contained in
\[
\bigcup_{\mu_M \in S_M(\mu,\nu)} X^M_{\mu_M} (b).
\] Denote by $\beta$ the map
\[
\beta:X^G_\mu(b) \to \bigcup_{\mu_M \in S_M(\mu,\nu)} X^M_{\mu_M} (b)
\] obtained by restriction from $\alpha$.
\item For $\mu_M \in S_M(\mu,\nu)$ and $x \in X^M_{\mu_M}(b)$ the fiber
$\beta^{-1}(x)$ is non-empty and
 \[
\dim \beta^{-1}(x)=d(\mu,\mu_M)+\langle
\rho,\bar\nu -\bar\nu_{\dom}\rangle-\langle 2\rho_N,\bar\nu \rangle.
\]  In particular $\beta$ is surjective.
\item For all $\mu_M \in S_M(\mu,\nu)$ and all $\lambda \in \Lambda_M$ the set
$\beta^{-1}\bigl( X^M_{\mu_M}(b)\cap \eta_M^{-1}(\lambda)\bigr)$ is locally
closed in
$X^G_\mu(b)$, and
\[
\dim \beta^{-1}\bigl( X^M_{\mu_M}(b)\cap \eta_M^{-1}(\lambda)\bigr) =\dim
X^M_{\mu_M}(b) +d(\mu,\mu_M)+\langle
\rho,\bar\nu -\bar\nu_{\dom}\rangle-\langle 2\rho_N,\bar\nu \rangle.
\]
\item If $X^G_\mu(b)$ is non-empty, its dimension is given by
\[
\sup\{ \dim X^M_{\mu_M}(b) +d(\mu,\mu_M): \mu_M \in S_M(\mu,\nu)\}+\langle
\rho,\bar\nu -\bar\nu_{\dom}\rangle-\langle 2\rho_N,\bar\nu \rangle.
\]
\end{enumerate}
\end{proposition}

\begin{proof} We begin by proving (2). Let $x \in X^G_\mu(b)$. Write
$x=mnx_0$ with $m
\in M(L)$, $n \in N(L)$. Then
\[ n^{-1}m^{-1}b\sigma(m)\sigma(n) \in K\epsilon^\mu K,
\]  from which it follows that $m^{-1}b\sigma(m)$ lies in
$K_M\epsilon^{\mu_M}K_M$ for a unique $\mu_M \in S_M(\mu)$. Thus
$\alpha(x)=mx_0 \in X^M_{\mu_M}(b)$, showing that $X^M_{\mu_M}(b)$ is
non-empty and hence that $\mu_M \in S_M(\mu,\nu)$. This proves (2).

Next we prove (3). Let $\mu_M \in S_M(\mu,\nu)$ and $x \in X^M_{\mu_M}(b)$.
Choose $m \in M(L)$ such that $x=mx_0$. Then
$b':=m^{-1}b\sigma(m) \in K_M\epsilon^{\mu_M} K_M$ and
\[
\beta^{-1}(x)=f_{b'}^{-1}\bigl(K\epsilon^\mu Kb'^{-1} \cap N(L)
\bigr)/N(\mathfrak o),
\] where $f_{b'}$ is the  morphism $N(L) \to N(L)$ defined by
\[ f_{b'}(n)=n^{-1}b'\sigma(n)b'^{-1}.
\] Since $f_{b'}$ is surjective, the fiber $\beta^{-1}(x)$ is non-empty.

Since $b'$ is $\sigma$-conjugate to $b$ in $M(L)$, its Newton point is also
$\bar\nu$. Therefore Proposition \ref{prop.reldimfP} tells us that
\[
\dim \beta^{-1}(x)=\dim(K\epsilon^\mu Kb'^{-1} \cap N(L))+\langle
\rho,\bar\nu -\bar\nu_{\dom}\rangle.
\] Furthermore
\[ N(L)b'x_0 \cap Kx_\mu=\bigl( K\epsilon^\mu Kb'^{-1} \cap N(L) \bigr)
/b'N(\mathfrak o) b'^{-1},
\] from which it follows that
\[
\dim \bigl( K\epsilon^\mu Kb'^{-1} \cap N(L)
\bigr)=d(\mu,\mu_M)+\dim(b' N(\mathfrak o)b'^{-1}),
\]  because $b' \in K_M \epsilon^{\mu_M}K_M$. Since $\dim b'N(\mathfrak
o)b'^{-1}$ is obviously equal to $-\langle 2\rho_N,\bar\nu \rangle$, we
conclude that
\[
\dim \beta^{-1}(x)=d(\mu,\mu_M)+\langle
\rho,\bar\nu -\bar\nu_{\dom}\rangle-\langle 2\rho_N,\bar\nu \rangle.
\]  This proves (3), and (1) follows from (2), (3), Lemma \ref{lem.3subsets}
and the fact that
$X^M_{\mu_M}(b)$ is non-empty for all $\mu_M \in S_M(\mu,\nu)$.

Since $\dim \beta^{-1}(x)$ is independent of $x \in X^M_{\mu_M}(b)$, part (4)
follows from part (3)  and the fact  (noted already in the introduction) that
\[
\dim \bigl ( X^M_{\mu_M}(b)\cap
\eta_M^{-1}(\lambda) \bigr)=\dim  X^M_{\mu_M}(b)
\] for all $\lambda \in \Lambda_M$.

Finally we prove (5). Assume that $X^G_\mu(b)$ is non-empty, so that
$S_M(\mu,\nu)$ is non-empty as well. Then $X^G_\mu(b)$ is set-theoretically
the union of the locally closed pieces
\begin{equation}
\beta^{-1}\bigl ( X^M_{\mu_M}(b)\cap
\eta_M^{-1}(\lambda) \bigr),
\end{equation} where $\mu_M$ ranges over the finite set $S_M(\mu,\nu)$ and
$\lambda$ ranges through $\Lambda_M$.
 The space $X^G_\mu(b)$ is preserved by the action of $A_M(F)$, and the action
of
$a \in A_M(F)$ on $X^G_\mu(b)$ sends
$\beta^{-1}\bigl ( X^M_{\mu_M}(b)\cap
\eta_M^{-1}(\lambda) \bigr)$ to $\beta^{-1}\bigl ( X^M_{\mu_M}(b)\cap
\eta_M^{-1}(\lambda+\lambda_a) \bigr)$, where $\lambda_a:=\eta_M(a)$.  Since
$\eta_M(A_M(F))$ has finite index in $\Lambda_M$,  the action of $A_M(F)$ has
finitely many orbits on the set of pieces in our decomposition. Reasoning as
in \ref{subsec.af.stable}, we conclude that
$\dim X^G_\mu(b)$ is equal to  the supremum of the dimensions of these pieces,
and so we obtain (5) from (4).
\end{proof}

\subsection{Remark} Part (1) of the proposition we just proved is not new;  it
was proved in \cite{kottwitz-rapoport02}. Moreover the strategy used to prove
the proposition is simply a refinement of the method used to prove
non-emptiness in \cite{kottwitz-rapoport02}.

\subsection{Reduction of Rapoport's dimension conjecture to the basic case} We
continue with $G$-dominant $\mu$ and a basic element $b \in M(L)$, as above.
We again use $b$ to get $\nu \in \Lambda_M$ and the Newton point
$\bar\nu \in \mathfrak a_M$.

We write $\defect_M(b)$ for the \emph{defect} of $b$, a non-negative integer
no greater than the semisimple rank of $M$. By definition
$\defect_M(b)$ is simply the $F$-rank of $M$ minus the $F$-rank of the inner
form $J$ of $M$ obtained by twisting $M$ by $b$. (The group $J(F)$ is the
$\sigma$-centralizer of $b$ in $M(L)$.) Rapoport's conjecture (applied to
$M$ and the basic element $b$) states that for all $M$-dominant $\mu_M$ with
$p_M(\mu_M)=\nu$ the dimension of $X^M_{\mu_M}(b)$ is given by
\begin{equation}\label{eq.rap.dim.f}
\dim X^M_{\mu_M}(b)=\langle \rho_M,\mu_M\rangle-\frac{1}{2}\defect_M(b).
\end{equation}
Note that
$\langle \rho_M,\mu_M
\rangle$ and $\defect_M(b)/2$ lie in $\frac{1}{2}\mathbb Z$, but that their
difference is always an integer.

Let's assume this formula is true and combine it with Proposition
\ref{prop.red.basic} in order to calculate $\dim X^G_\mu(b)$. Recall that
$X^G_\mu(b)$ is non-empty if and only if $\nu$ lies in the image of
$\Sdom$.

\begin{theorem}\label{thm581}
Suppose that $X^G_\mu(b)$ is non-empty.  Assume that
\eqref{eq.rap.dim.f} is true for all $\mu_M \in S_M(\mu,\nu)$. Then the
dimension of $X^G_\mu(b)$ is given by
\[
\dim X^G_\mu(b)=\langle \rho,\mu-\bar\nu_{\dom} \rangle
-\frac{1}{2}\defect_G(b).
\]
\end{theorem}
\begin{proof} By Proposition \ref{prop.red.basic} $\dim X^G_\mu(b)$ is the sum
of
\[
\langle
\rho,\bar\nu -\bar\nu_{\dom}\rangle-\langle 2\rho_N,\bar\nu \rangle
\] and  the supremum of the numbers
\begin{equation}\label{eq.presup}
\dim  X^M_{\mu_M}(b) +d(\mu,\mu_M)
\end{equation} as $\mu_M$ ranges through $ S_M(\mu,\nu)$.

Combining Rapoport's dimension conjecture \eqref{eq.rap.dim.f} with Corollary
\ref{cor.dim.inters}, we see that the number \eqref{eq.presup} is bounded
above by
\[
\langle \rho_M,\mu_M\rangle-\defect_M(b)/2+\langle \rho,\mu+\mu_M
\rangle -2\langle
\rho_M,\mu_M \rangle
\] with equality whenever $\mu_M \in \Smax$. Using the obvious equality
$\langle \rho_N,\mu_M\rangle=\langle \rho_N,\bar\nu \rangle$, the displayed
expression above simplifies to
\[
\langle \rho,\mu\rangle +\langle
\rho_N,\bar\nu \rangle -\defect_M(b)/2,
\] a number which is independent of $\mu_M$. Since we have supposed that
$X^G_\mu(b)$ is non-empty, $\nu$ lies in the image of
$\Sdom$. Therefore by Lemma \ref{lem.3subsets} $\nu$ is in the image of
$\Smax$, which shows that the supremum of the numbers
\eqref{eq.presup} is equal to the last displayed expression. Adding
$\langle
\rho,\bar\nu -\bar\nu_{\dom}\rangle-\langle 2\rho_N,\bar\nu \rangle$
to this, and using that $\langle \rho_M,\bar\nu \rangle=0$ since $\bar\nu
\in \mathfrak a_M$, we find that
\[
\dim X^G_\mu(b)=\langle \rho,\mu-\bar\nu_{\dom}\rangle -\defect_M(b)/2.
\]
Now it remains only to observe that $\defect_M(b)=\defect_G(b)$. Indeed, the
$F$-rank of a reductive group is the same as that of any of its Levi
subgroups, so the $F$-ranks of $G$ and $M$ are the same. Moreover, the group
$J$ attached to $(M,b)$ is a Levi subgroup of the group $J_G$ attached to
$(G,b)$, so their $F$-ranks are the same.
\end{proof}

\subsection{Reduction to the superbasic case}\label{subsec.superbasic}
Let $b \in G(L)$. As in the introduction $b$ determines a homomorphism
$\nu_b:\mathbb D \to G$ over $L$. As in \cite[(4.4.3)]{kottwitz85} we then
have
\begin{equation}\label{eq.conjbnu}
\Int(b) \circ \sigma(\nu_b)=\nu_b.
\end{equation}
Replacing $b$ by a $\sigma$-conjugate, we again assume that $\nu_b$ factors
through $A$, which guarantees that $\sigma(\nu_b)=\nu_b$. Equation
\eqref{eq.conjbnu} then says that $b \in M(L)$, where $M$ is the Levi
subgroup of $G$ over $F$ obtained as the centralizer of $\nu_b$ in $G$.
Since $\nu_b$ is central in $M$, the element $b$ is basic in $M(L)$. We have
seen that if Rapoport's dimension conjecture is true for $(M,b)$, then it is
true for
$(G,b)$.

Thus it is enough to prove Rapoport's conjecture when $b$ is basic. Of
course it may still happen that a basic element $b$ is contained in a proper
Levi subgroup of $G$. For example this happens when $b=1$, as long as $G$ is
not a torus. We say that an element $b \in G(L)$ is \emph{superbasic} if no
$\sigma$-conjugate of $b$ lies in a proper Levi subgroup of $G$ over $F$.
Clearly superbasic elements are basic.

Suppose that $b \in G(L)$ is basic. Then the group $J$ discussed in the
introduction is the inner form of $G$ obtained by twisting the action of
$\sigma$ by $b$, so that
\[
J(F)=\{g \in G(L): g^{-1}b\sigma(g)=b\},
\]
or in other words $J(F)$ is the centralizer of $b\sigma$ in $G(L)$.

\begin{lemma}\label{lem.superbasic}
The basic element $b$ is superbasic if and only if $J$ is anisotropic modulo
the center $Z(G)$ of $G$.
\end{lemma}
\begin{proof}
The group $J/Z(G)$ is not anisotropic if and only if there exists a
non-central homomorphism $\mathbb G_m \to J$ which is defined over $F$. This happens if and only if
there exist non-central $\mu \in X_*(A)$ and an element $g \in G(L)$ such
that $b\sigma$ and $\Int(g) \circ \mu$ commute. Now $b\sigma$ and $\Int(g)
\circ \mu$ commute if and only if $g^{-1}(b\sigma)g$ and $\mu$ commute, and
this happens if and only if $g^{-1}b\sigma(g)$ and $\mu$ commute (since
$\sigma$ and $\mu$ commute), which is the same as saying that
$g^{-1}b\sigma(g)$ lies in the Levi subgroup of $G$ obtained as the
centralizer of $\mu$.
\end{proof}

It follows from Theorem \ref{thm581} that it is enough to prove Rapoport's
conjecture when $b$ is superbasic. In the introduction we observed that the
intersections of $X_\mu(b)$ with the various connected components of the
affine Grassmannian are all isomorphic to each other. Let us denote by
$X_\mu(b)_\lambda$ the intersection of $X_\mu(b)$ with the connected
component of the affine Grassmannian indexed by $\lambda \in \Lambda_G$. It
is easy to see that
\[
X_\mu(b)_\lambda=X_{\mu_{\ad}}(b_{\ad})_{\lambda_{\ad}},
\]
where the subscript ``$\ad$'' is being used to indicate the corresponding
objects for the adjoint group $G/Z(G)$ of $G$. Thus it is enough to prove
Rapoport's conjecture for superbasic elements in adjoint groups.

Any adjoint group is a direct product of simple groups. Therefore it is
enough to prove Rapoport's conjecture for superbasic elements in simple
(split) groups. From Lemma \ref{lem.superbasic} we know that a basic element
is superbasic if and only if the inner form $J$ of our simple group $G$ is
anisotropic. This can happen only when $G$ is of type $A_n$, and then it is
easy to see which elements are superbasic. Thus we obtain

\begin{proposition}
If Rapoport's conjecture is true for all basic elements $b$ in $GL_n(L)$ such
that the valuation of $\det(b)$ is relatively prime to $n$, then Rapoport's
conjecture is true in general.
\end{proposition}
\begin{proof}
Observe that $GL_n(L) \to PGL_n(L)$ is surjective, and that a
basic element $b \in GL_n(L)$ is superbasic if and only if $\val(\det b)$ is
relatively prime to $n$.
\end{proof}

\section{Affine Deligne-Lusztig varieties inside the affine flag manifold}
\subsection{Statement of the problem}\label{Sect61} We continue with $G$, $A$, $W$ as
before, but make a few changes in notation. We write
$\mathcal B(A)$ for the set of Borel subgroups $B=AU$ containing $A$. We fix
an alcove in the apartment $\mathfrak a:=X_*(A)_\mathbb R$ associated to $A$,
and we denote by $I$ the corresponding Iwahori subgroup of $G(L)$. Note that
$A(\mathfrak o) \subset I$.  We write
$\tilde W$ for the affine Weyl group $W \ltimes X_*(A)=N_{G(L)}(A)/A(\mathfrak
o)$.

We will now let
$X$ denote the affine flag manifold $G(L)/I$ and $\mathbf a_1$ its base point.
For
$x \in \tilde W$ we put $\mathbf a_x:=x\mathbf a_1 \in X$. Recall that
\begin{equation} X =\coprod_{x \in \tilde W} I\mathbf a_x
\end{equation} and that for any $B=AU \in \mathcal B(A)$
\begin{equation} X =\coprod_{x \in \tilde W} U(L)\mathbf a_x.
\end{equation}

Given $b \in G(L)$ and $x \in \tilde W$, we get the affine Deligne-Lusztig
variety
\begin{equation} X_x(b):=\{g \in G(L)/I : g^{-1}b \sigma(g) \in IxI \},
\end{equation}  a locally closed subset of $X$, possibly empty, whose
dimension we are going to compute for all
$b$ of the form $b=\epsilon^\nu$ for some $\nu \in X_*(A)$. However the
calculation will in fact only express $\dim X_x(\epsilon^\nu)$ in terms of
the (unknown) dimensions of intersections of $I$- and $U(L)$-orbits in $X$.
This is at least satisfactory for computer experiments, since the dimensions
of these intersections can be determined by an algorithm involving only
$\tilde W$ and the set of affine roots.

Here is the algorithm (see \cite{Da} and
\cite[Prop.~2.3.12]{bruhat-tits72}). Consider the $I$-orbit of
$\mathbf a_x$, which is an affine space $\mathbb A$ of dimension
$\ell(x)$. Choose a reduced expression for $x$. Equivalently, choose
a minimal gallery $\mathcal{G}$ joining $\mathbf a_1$ to $\mathbf
a_x$. (Here $\mathbf a_1$ is the alcove fixed by our Iwahori $I$, and the notions of length and 
reduced expression are defined using the simple reflections through the walls of $\mathbf a_1$.)  
The affine space $\mathbb A$ is the product of a succession
of affine lines, one for each wall between two alcoves of
$\mathcal G$. We will now decompose this affine space into
finitely many pieces, each of which is contained in a single
$U(L)$-orbit, by specifying which $U(L)$-orbit a point $y \in
\mathbb A$ is in.
To the Borel $B=AU$ corresponds a Weyl chamber $C_U$, which is the unique one
with the property that $A(\mathfrak o)U(L)$ is the union of the fixers in $G(L)$ of all
"quartiers" of the form $x + C_U$, for $x \in X_*(A)_{\mathbb R}$ (see
\cite{bruhat-tits72} 4.1.5, 4.4.4).
We can use any alcove in $C_U$ to  retract $y$ into the
standard apartment. The result depends on the choice of alcove,
but for alcoves deep inside $C_U$, all retractions of $y$ are the
same (cf.~loc.~cit.~2.9.1). 
The alcove obtained as this common retraction of $y$ lies in
the $U(L)$-orbit of $y$, telling us which $U(L)$-orbit $y$ belongs
to. We work our way up the minimal gallery $\mathcal{G}$,
simultaneously considering the retractions of all galleries $I
\mathcal{G}$ relative to an alcove sufficiently deep in $C_U$. At
the step corresponding to the edge between the $i^{\text{th}}$ and
$i+1^{\text{st}}$ alcoves of $\mathcal{G}$, there are two
possibilities: first, it may happen that the whole affine line at
this step folds away from the alcove we used to get the
retraction; second, it may happen that one point in the affine
line folds towards this retraction alcove, and the rest fold away
from it. Keeping track of all these folding possibilities gives a
decomposition of $\mathbb A$ into finitely many pieces, each of
which is contained in a single $U(L)$-orbit and is a product of
affine lines and affine lines minus a point. Taking the supremum
of the dimensions of the pieces lying in a given $U(L)$-orbit, we
obtain the dimension of the intersection of $\mathbb A$ with that
$U(L)$-orbit.

\subsection{Intersections of $I$- and $U(L)$-orbits} Let $B=AU \in \mathcal
B(A)$ and let $x,y \in \tilde W$.   Define a non-negative integer $d(x,y,B)$ by
\begin{equation} d(x,y,B):=\dim I\mathbf a_x \cap U(L)\mathbf a_y.
\end{equation}
(We set $d(x, y, B) = -\infty$ if this intersection is empty.)
It is clear that
\begin{equation} I\mathbf a_x \cap U(L)\mathbf a_y=\bigl( IxIy^{-1} \cap
U(L)\bigr)/(U(L) \cap yIy^{-1}),
\end{equation} and hence that
\begin{equation} d(x,y,B)=\dim IxIy^{-1} \cap U(L)-\dim U(L) \cap yIy^{-1}.
\end{equation} Here we are again using our dimension theory for $U(L)$, it
being obvious that both $ IxIy^{-1} \cap U(L)$ and $U(L) \cap yIy^{-1}$ are
admissible subsets of $U(L)$.

In the special case $y=\epsilon^\nu$ with $\nu \in X_*(A)$,  we have
\[ U(L)\cap yIy^{-1}=\epsilon^\nu(U(L) \cap I)\epsilon^{-\nu}
\] and hence
\[
\dim U(L)\cap yIy^{-1}=\dim  U(L) \cap I -2\langle \rho_B,\nu \rangle,
\] where $\rho_B$ is the half-sum of the roots of $A$ that are positive for
$B$.  Therefore we conclude that
\begin{equation}\label{eq.owiejfw}
\dim IxI\epsilon^{-\nu} \cap U(L)-\dim U(L) \cap I=d(x,\epsilon^\nu,B)-2\langle
\rho_B,\nu \rangle.
\end{equation}

\subsection{Computation of $\dim X_x(\epsilon^\nu)$}\label{Sect63} Pick $B=AU \in \mathcal
B(A)$.  As before we use the fact that $X_x(\epsilon^\nu)$ is $A(F)$-stable to
compute its dimension. The only difference is that $A(F)$ does not permute the
$U(L)$-orbits
$U(L)\mathbf a_y$ transitively. Nevertheless we have the only slightly more
complicated formula
\begin{equation}
\dim X_x(\epsilon^\nu)=\sup_{w \in W} \dim X_x(\epsilon^\nu) \cap U(L)\mathbf
a_w.
\end{equation} Clearly left multiplication by $w^{-1}$ gives an isomorphism
from
\[ X_x(\epsilon^\nu) \cap U(L)\mathbf a_w
\] to
\[ X_x(\epsilon^{w^{-1}\nu}) \cap w^{-1}U(L)w\mathbf a_1,
\] and hence
\begin{equation}
\dim X_x(\epsilon^\nu)=\sup_{w \in W} \dim X_x(\epsilon^{w\nu}) \cap
wUw^{-1}(L)\mathbf a_1.
\end{equation}

So we need to understand $ \dim X_x(\epsilon^{\nu}) \cap U(L)\mathbf a_1$ for
all $\nu \in X_*(A)$ and all $B=AU \in \mathcal B(A)$.  Again writing
$f_\nu$ for the map $u \mapsto u^{-1}\epsilon^\nu \sigma (u)\epsilon^{-\nu}$
from $U(L)$ to itself, we have
\[
 X_x(\epsilon^{\nu}) \cap U(L)\mathbf a_1=f_\nu^{-1}\bigl(IxI\epsilon^{-\nu}
\cap U(L)\bigr)/U(L)\cap I
\] and hence
\begin{equation}
\dim X_x(\epsilon^{\nu}) \cap U(L)\mathbf a_1=\dim
f_\nu^{-1}\bigl(IxI\epsilon^{-\nu} \cap U(L)\bigr)-\dim U(L)\cap I
\end{equation} Using Proposition \ref{prop.main}, Lemma \ref{lem.rho.dom} and
\eqref{eq.owiejfw}, we see that
\begin{equation}
\dim f_\nu^{-1}\bigl(IxI\epsilon^{-\nu} \cap U(L)\bigr)-\dim U(L)\cap
I=d(x,\epsilon^\nu,B)-\langle \rho_B,\nu+\nu_B \rangle,
\end{equation} where $\nu_B$ is the unique element in the $W$-orbit of $\nu$
that is dominant for
$B$. Thus we have proved
\begin{theorem}\label{thm.631}
For any $x \in \tilde W$, $\nu \in X_*(A)$ and $B \in
\mathcal B(A)$ there is an equality
\[
\dim X_x(\epsilon^\nu)=\sup_{w \in W} \{d(x,\epsilon^{w\nu},wBw^{-1})-\langle
\rho_{B},\nu+\nu_{B}\rangle\}.
\] In particular for $\nu=0$ we obtain
\[
\dim X_x(1)=\sup_{B \in \mathcal B(A)} d(x,1,B).
\]
\end{theorem}

This means in particular that $X_x(\epsilon^\nu) \ne \emptyset$ if and only if there
exists $w \in W$ such that $d(x,\epsilon^{w\nu},wBw^{-1}) \ge 0$.

\section{Computations}\label{sec.computations}
\subsection{The algorithm}

In this section we present computational results on the
non-emptiness and dimension of affine Deligne-Lusztig varieties in
the affine flag manifold for a few groups of low rank. These
results are computed according to Theorem~\ref{thm.631} of
Section~\ref{Sect63}, using the algorithm of Section~\ref{Sect61}
to compute the dimensions $d(x,\epsilon^{w\nu},wBw^{-1})$.

Computations such as the ones presented here can in principle be
done by hand, but their complexity leads us to use a computer
program. The program is capable of representing all elements of
the affine Weyl group of a root system up to a given length. It
can also multiply these elements. The program can represent affine
weights and coweights, as well as the action of the affine Weyl
group on these weights and coweights.

Even using a computer, the computations presented here took some
time. For instance, computing which alcoves correspond to
non-empty affine Deligne-Lusztig varieties together with their
dimensions in the $G_2$ case up to length $56$ took several days
on a medium-fast PC. To some extent, the program makes use of the
automatical structure of the affine Weyl group (see \cite{Ca}, for
instance). However, multiplication is carried out in the spirit of
the method explained on p.~96 in loc.~cit., not by using a
multiplier automaton (see section 3 in loc.~cit.). Rewriting this
part of the program would certainly lead to a significant
speed-up. Memory consumption was never a problem.

\subsection{Reuman's conjecture}\label{sub.comp.2}

We fix a Borel subgroup $B\subset G$ containing $A$, and let $\alpha_i$ denote the simple 
$B$-positive roots of $A$, i.e.~the simple roots in ${\rm Lie}(B)$.  
Let $C_0$ denote the dominant Weyl 
chamber, which by definition is the set 
$\{ x \in X_*(A)_{\mathbb R} ~ | ~ \langle \alpha_i, x \rangle > 0, \,\, \forall i \}$.  
We call the
unique alcove in the dominant Weyl chamber whose closure contains the
origin the base alcove. As Iwahori subgroup $I$ we choose the Iwahori
fixing this base alcove;  $I$ is the inverse image of the opposite Borel
group of $B$ under the projection $K \longrightarrow G(\bar{k})$.
We recall Reuman's conjecture from \cite{Reu04}, starting with a
definition.

\begin{defn}
For $w \in W$ let $w(C_0)$ denote the Weyl chamber corresponding to $w$.
For each simple positive root $\alpha_i$ of our root
system, let $\varpi^\vee_i$ denote the corresponding fundamental
coweight.  Then
we define the corresponding shrunken Weyl chamber $w(C_0)^{shr} \subset w(C_0)$ by
\[
w(C_0)^{shr} = w(C_0) + \sum_{\gfrac{i}{w\alpha_i > 0}} w\varpi^\vee_i.
\]
\end{defn}

\noindent The union of the shrunken Weyl chambers is the set of
all alcoves $\mathbf a$ in the standard apartment such that for
each positive finite root $\alpha$, either the base alcove and
$\mathbf a$ lie on opposite sides of the hyperplane
$\{\alpha=0\}$, or they lie on opposite sides of the hyperplane
$\{\alpha=1\}$.

We denote by $\eta_1 \colon W_a \longrightarrow W$ the usual
surjective homomorphism from the affine Weyl group to the finite
Weyl group, and by $\eta_2 \colon W_a \longrightarrow W$ the map
which associates to each alcove the finite Weyl chamber in which
it lies.

\begin{conjecture} \label{conj_b1}
If $x$ lies in the shrunken Weyl chambers, then $X_x(1) \ne \emptyset$ if and
only if
\begin{equation} \label{RC}
\eta_2(x)^{-1}\eta_1(x)\eta_2(x) \in W \setminus \bigcup_{T\subsetneq S} W_T,
\end{equation}
and in this case the dimension is given by
\[ \frac{1}{2}(\ell(x) + \ell(\eta_2(x)^{-1}\eta_1(x)\eta_2(x))).
\]
\end{conjecture}

\noindent Here $S$ denotes the set of simple reflections, and for
$T\subset S$, $W_T$ denotes the subgroup of $W$ generated by $T$.
The set $W \setminus \bigcup_{T\subsetneq S} W_T$ is all elements
of $W$ for which any reduced expression contains all simple
reflections.

In loc.~cit.~the conjecture was proved for the root systems $A_2$
and $C_2$. It is also true in all cases that were checked
computationally: $G_2$ up to length 56, $A_3$ up to length 38,
$C_3$ up to length 31, $A_4$ up to length 15, $C_4$ up to length
16. However, in the rank 4 case `most' of the alcoves of `small'
length lie outside the shrunken Weyl chambers: it would be
desirable to expand the range of our computations to encompass
greater lengths.

\begin{figure}
\caption{Non-empty affine Deligne-Lusztig varieties and their dimensions,
type $A_2$, $b=1$}
\includegraphics[width=13cm]{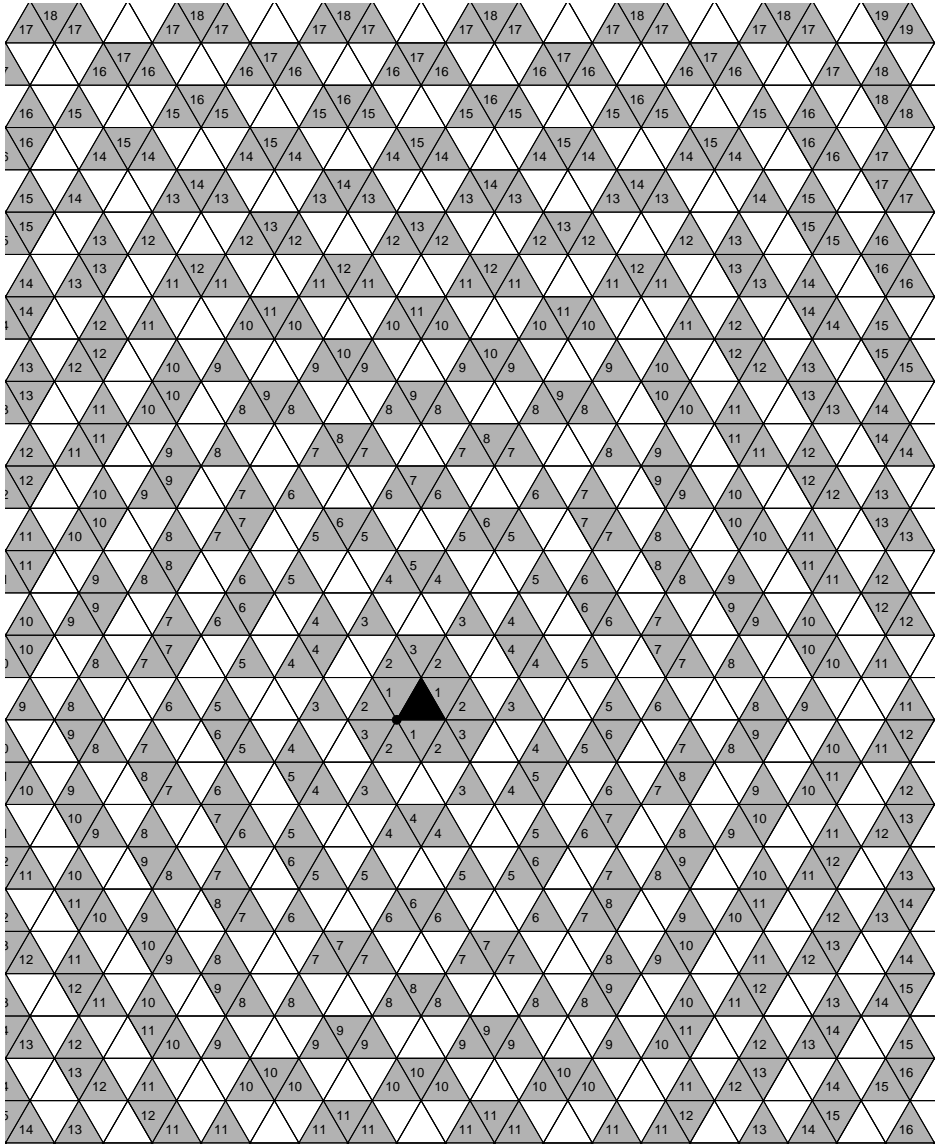} \label{1A}
\end{figure}

\begin{figure}
\caption{Non-empty affine Deligne-Lusztig varieties and their dimensions,
type $C_2$, $b=1$}
\includegraphics[width=13cm]{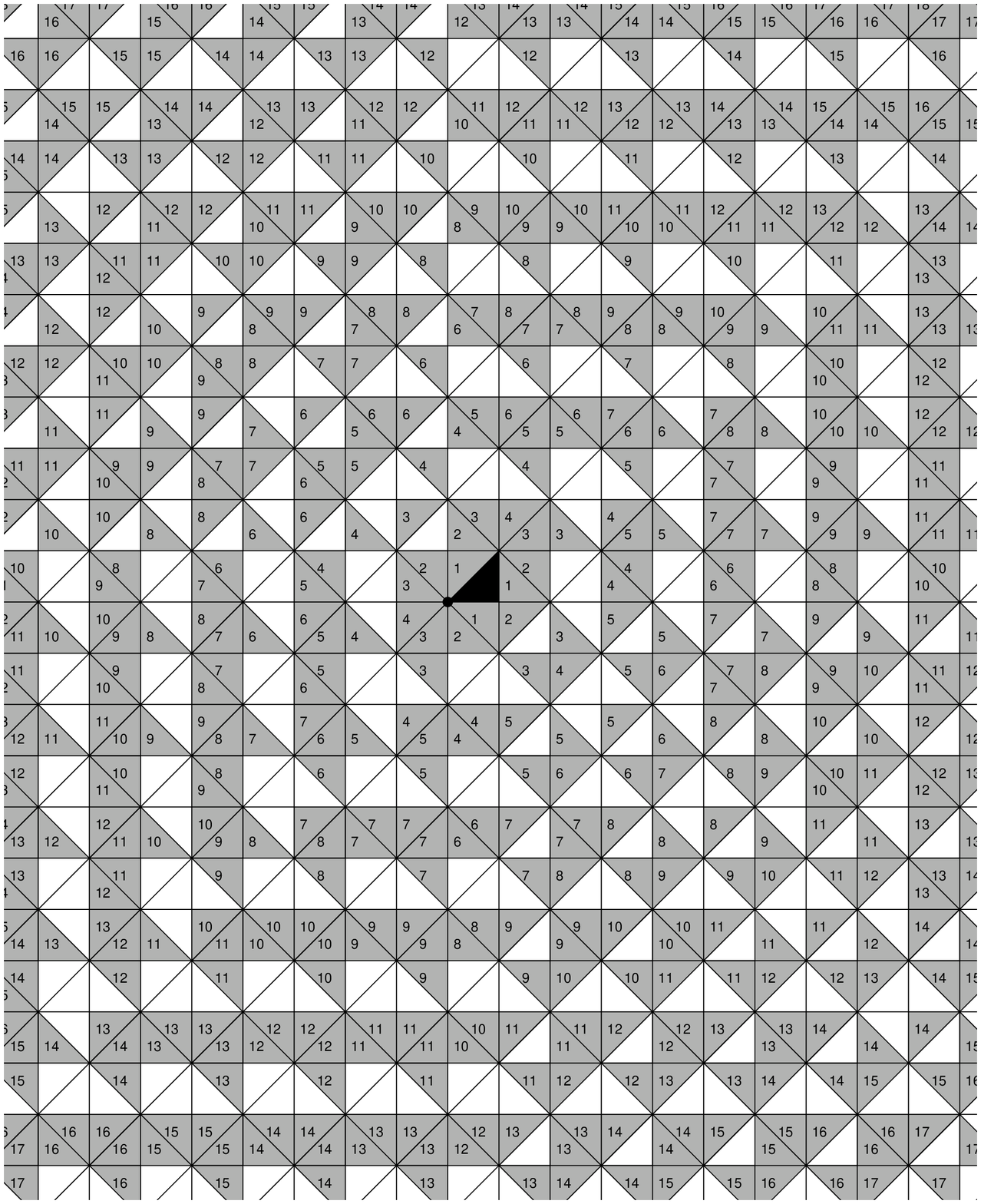} \label{1C}
\end{figure}

\begin{figure}
\caption{Non-empty affine Deligne-Lusztig varieties and their dimensions,
type $G_2$, $b=1$}
\includegraphics[width=13cm]{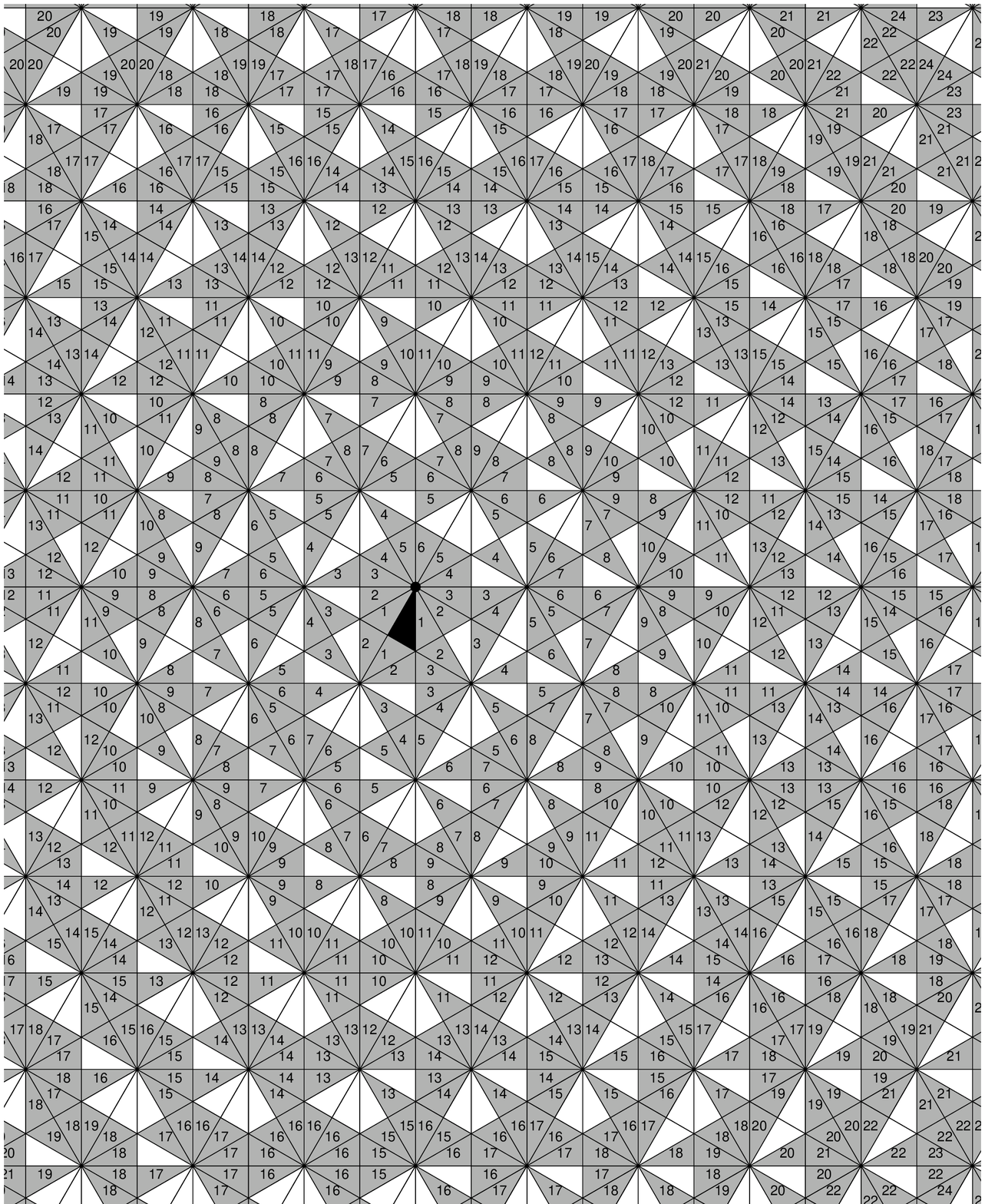}\label{1G}
\end{figure}

For root systems of rank 2, results can be represented by a
picture of the standard apartment (figure~\ref{1A} for $A_2$,
figure~\ref{1C} for $C_2$, and figure~\ref{1G} for $G_2$). Alcoves
corresponding to non-empty affine Deligne-Lusztig varieties are
gray; those corresponding to empty affine Deligne-Lusztig
varieties white; the base alcove is black. There is a dot at the
origin. Non-empty affine Deligne-Lusztig varieties have their
dimension on the corresponding alcove.

\begin{figure}
\caption{Comparison with Reuman's criterion, type $A_2$, $b=1$}
\includegraphics[width=13cm]{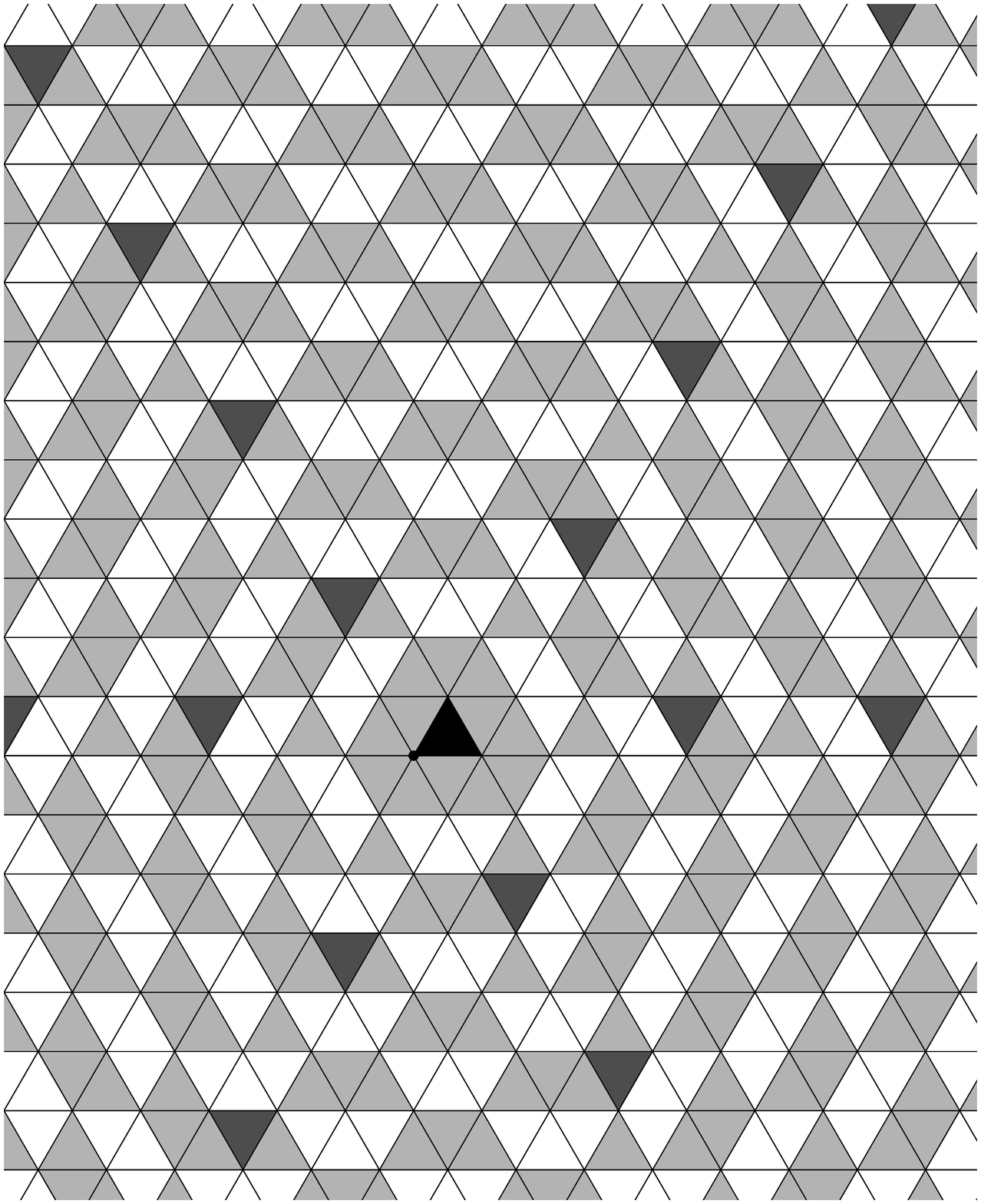} \label{2A}
\end{figure}

\begin{figure}
\caption{Comparison with Reuman's criterion, type $C_2$, $b=1$}
\includegraphics[width=13cm]{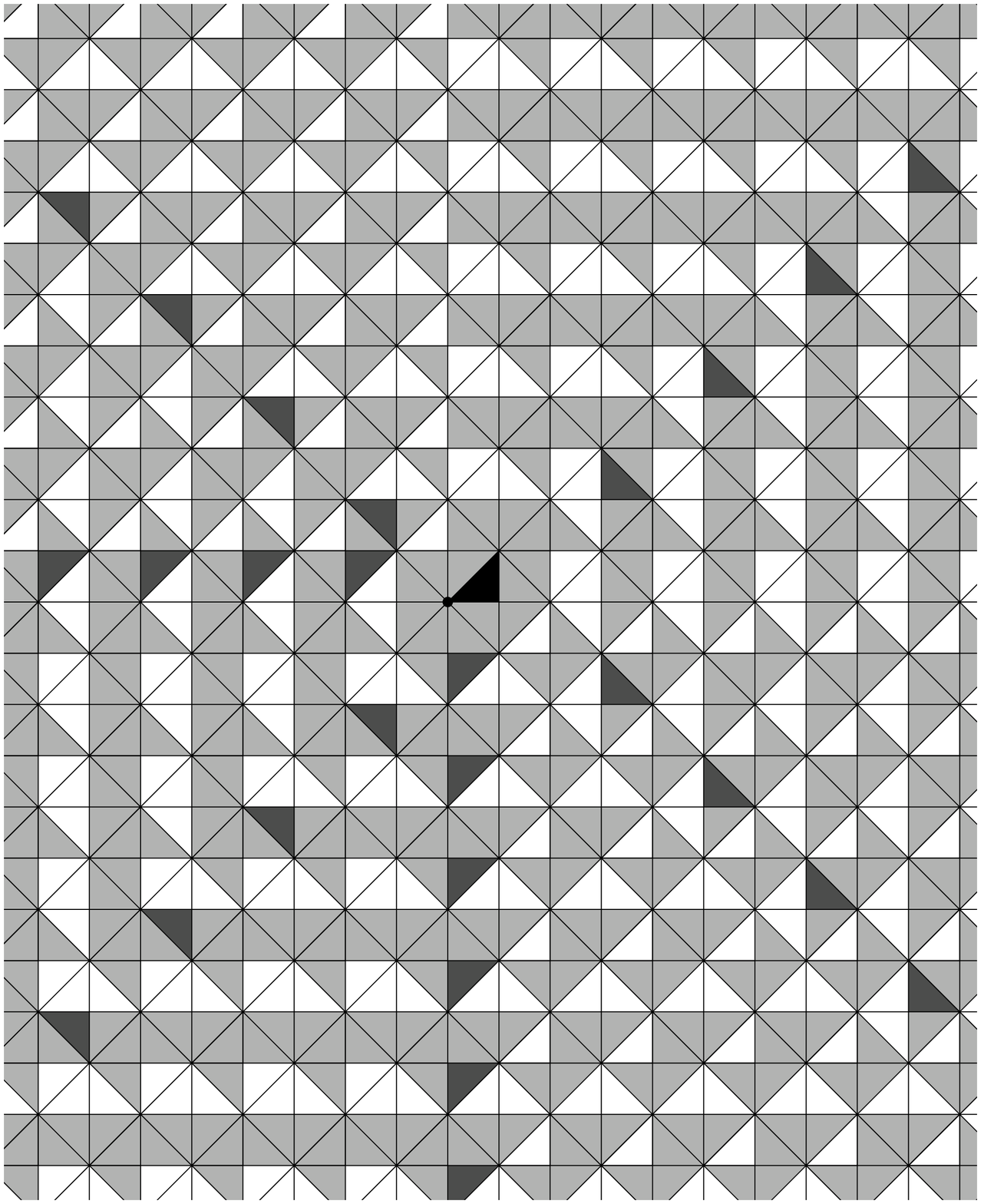} \label{2C}
\end{figure}

\begin{figure}
\caption{Comparison with Reuman's criterion, type $G_2$, $b=1$}
\includegraphics[width=13cm]{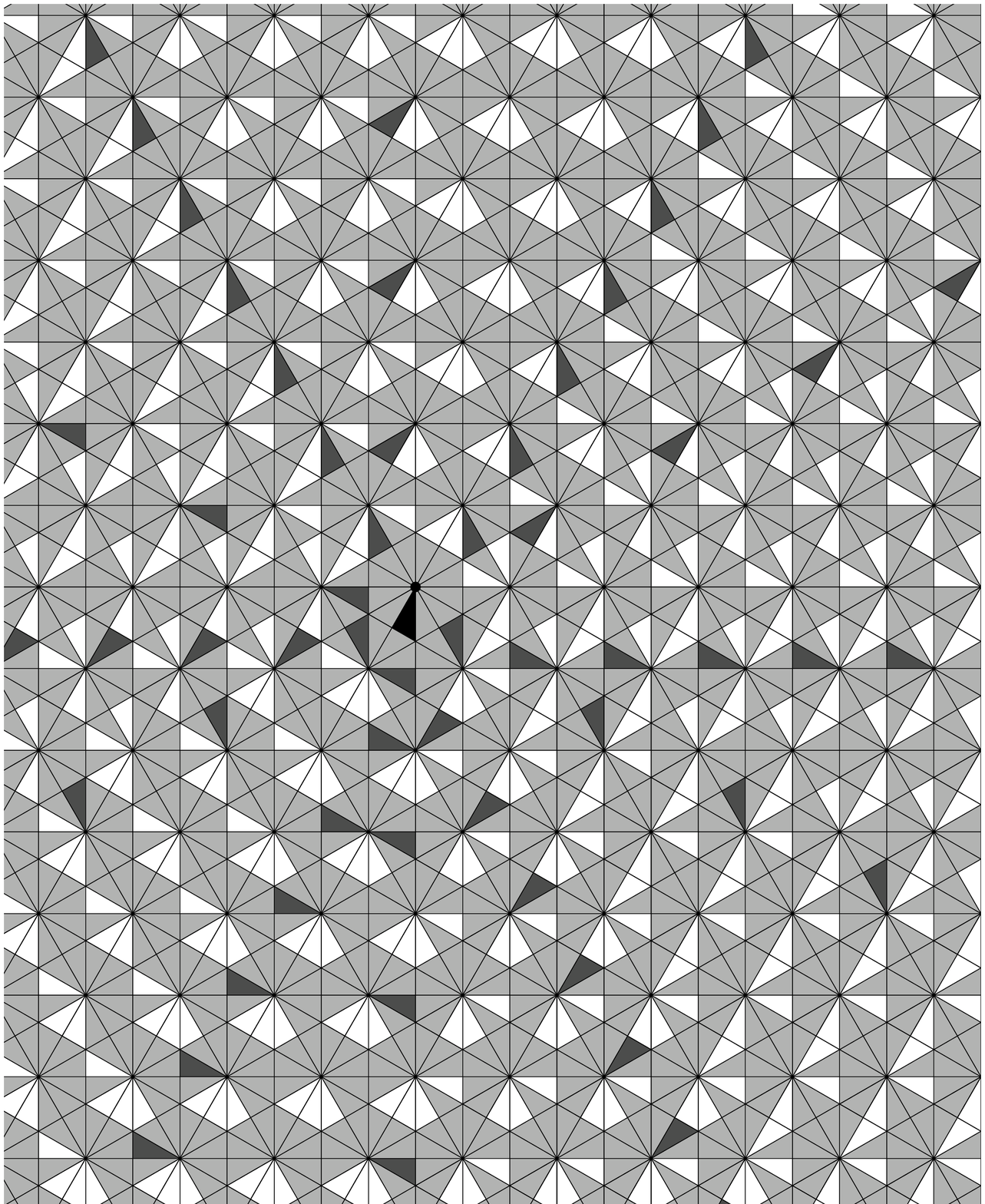} \label{2G}
\end{figure}

Outside the shrunken Weyl chambers, (\ref{RC}) is not equivalent
to $X_x(1) \ne \emptyset$. For instance, all $X_x(1)$ with $x \in
W$ are non-empty, even if $x \in \bigcup_{T\subsetneq S} W_T$.
However, for $x\not\in W$, (\ref{RC}) is a necessary criterion for
non-emptiness in all the cases we checked. In figures \ref{2A},
\ref{2C}, and \ref{2G} we show the standard apartments of the root
systems of type $A_2$, $C_2$ and $G_2$. The alcoves with non-empty
affine Deligne-Lusztig varieties are light gray, and those for
which (\ref{RC}) wrongly predicts non-emptiness are dark gray. All
dark gray alcoves lie outside the shrunken Weyl chambers. Except
in the $A_2$ case, even if the affine Deligne-Lusztig variety for
$x$ is non-empty, the predicted dimension may be wrong.

One can try to describe the set of alcoves for which the criterion
(\ref{RC}) is not equivalent to the non-emptiness of the
associated affine Deligne-Lusztig variety. For instance, in the
$A_2$ case, for $x\not\in W$, (\ref{RC}) wrongly predicts
non-emptiness of $X_x(1)$ if and only if the alcove corresponding
to $x$ lies outside the shrunken Weyl chambers and does not have
the same number of vertices as the base alcove lying on the wall
through the origin that it touches. We are not aware of a
description of this kind that holds for more general root systems.

The criterion (\ref{RC}) is not invariant under the symmetry of
root systems of type $A$ and $C$. It follows from \cite{Reu04} that 
for $A_2$, this realization
leads to a complete description of which affine Deligne-Lusztig
varieties are non-empty: $X_x(1) \ne\emptyset$ if and only if
(\ref{RC}) holds for $x$ and its images under rotation by 120 and
240 degrees about the center of the base alcove. This does not
work in the $C_2$ case, however.

\subsection{Lau's observation}

The following proposition is a direct consequence of \cite{kottwitz85}, Prop.
5.4.

\begin{proposition}
Let $G$ be semisimple and simply connected. Let $x\in W_a$, and
let $\underline{x} \in N_{G(L)}(A)$ be a representative of $x$
which lies in $G(F)$. Then $x$ has finite order if and only if
there exists an $y\in G(L)$ with $\underline{x}=y^{-1}\sigma(y)$.
In particular, in this case $X_x(1)$ is non-empty.
\end{proposition}

Starting from this proposition and analyzing Reuman's results in
the $A_2$ case, Eike Lau observed that one can get another
complete characterization of the non-empty affine Deligne-Lusztig
varieties for $A_2$:

\begin{proposition}
Let $G=SL_3$, and let $x \in W_a$, $\ell(x)>1$. Then
$X_x(1)\ne\emptyset$ if and only if there exists $n \ge 1$ with
$\ell(x^n) < \ell(x)-1$.
\end{proposition}

\noindent One can prove this proposition by systematically
analyzing Reuman's result for $SL_3$ (\cite{Reu04} or figure
\ref{1A}). It would be interesting to have a more conceptual
proof. For other root systems, a similar criterion does not seem
to hold (at least it is not sufficient to replace the offset $-1$
in the proposition by another integer). Nevertheless, this
criterion still seems to give a reasonable approximation to the
truth. It is possible that a minor modification of this criterion
might hold more generally.

\subsection{Partial folding sets}

\begin{figure}
\caption{Partial folding results, type $A_2$, $b=1$, $w={\rm id}$}
\includegraphics[width=13cm]{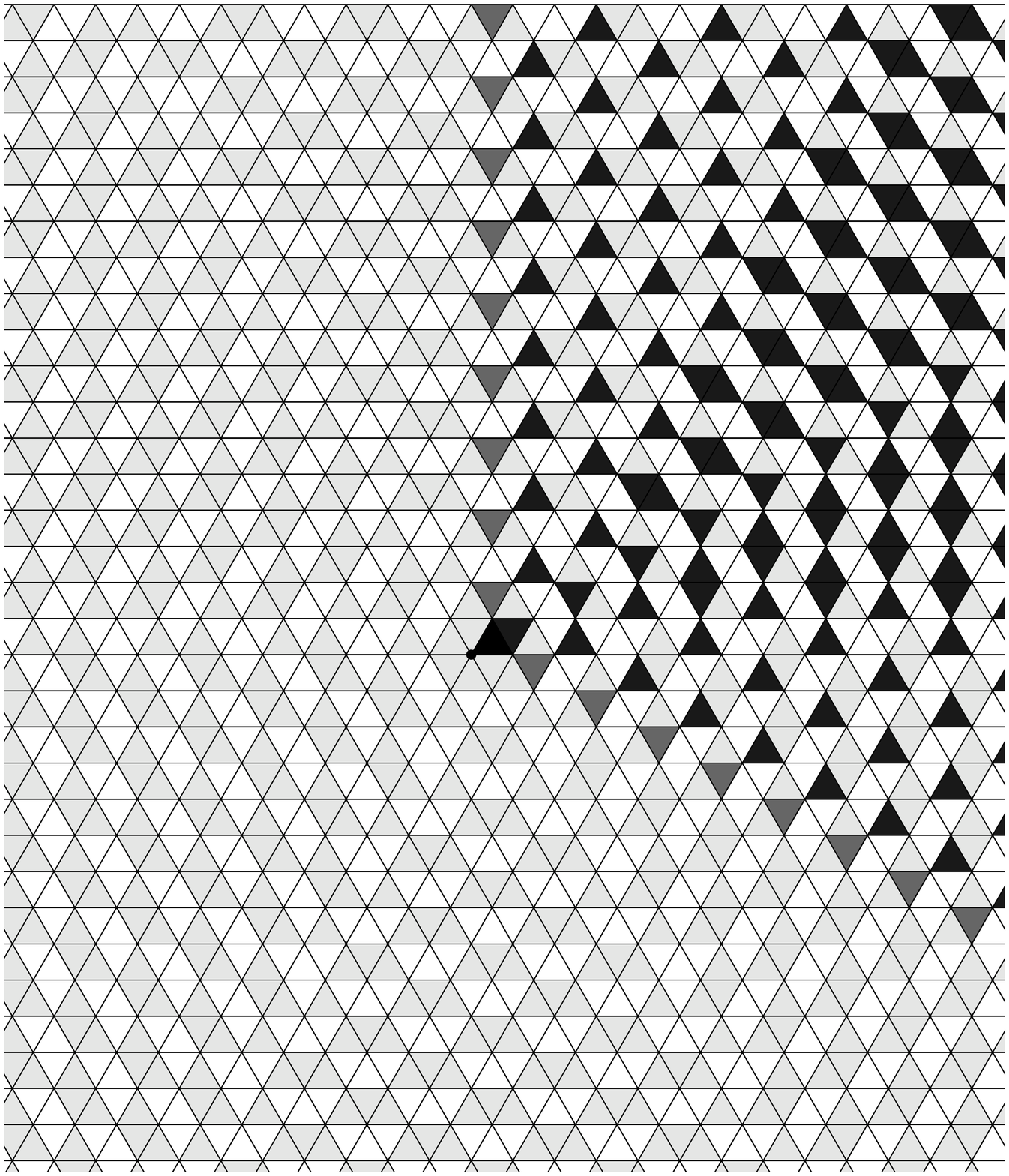} \label{3A.1}
\end{figure}

\begin{figure}
\caption{Partial folding results, type $A_2$, $b=1$, $w=s_1$}
\includegraphics[width=13cm]{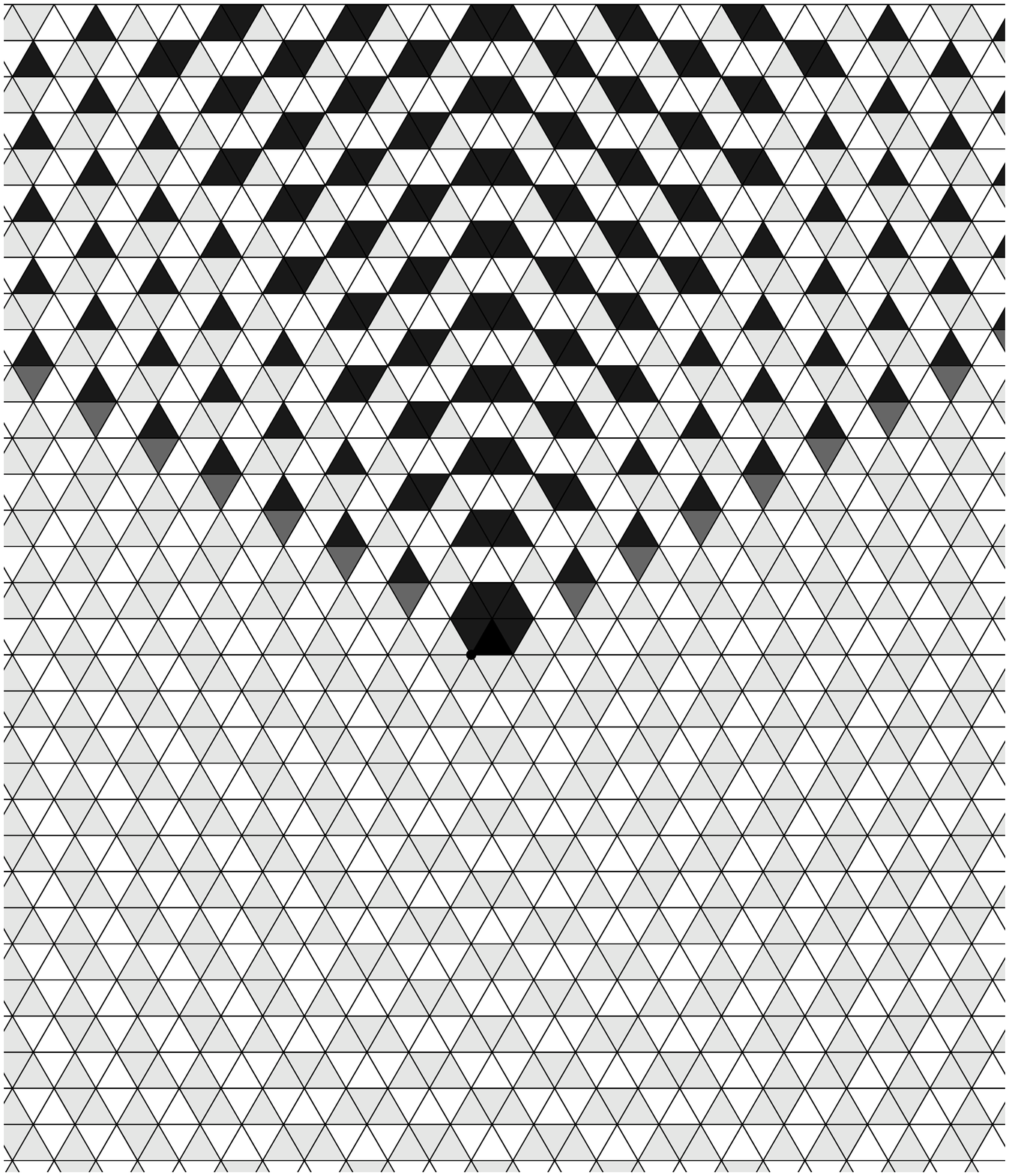} \label{3A.2}
\end{figure}

\begin{figure}
\caption{Partial folding results, type $A_2$, $b=1$, $w=s_1s_2$}
\includegraphics[width=13cm]{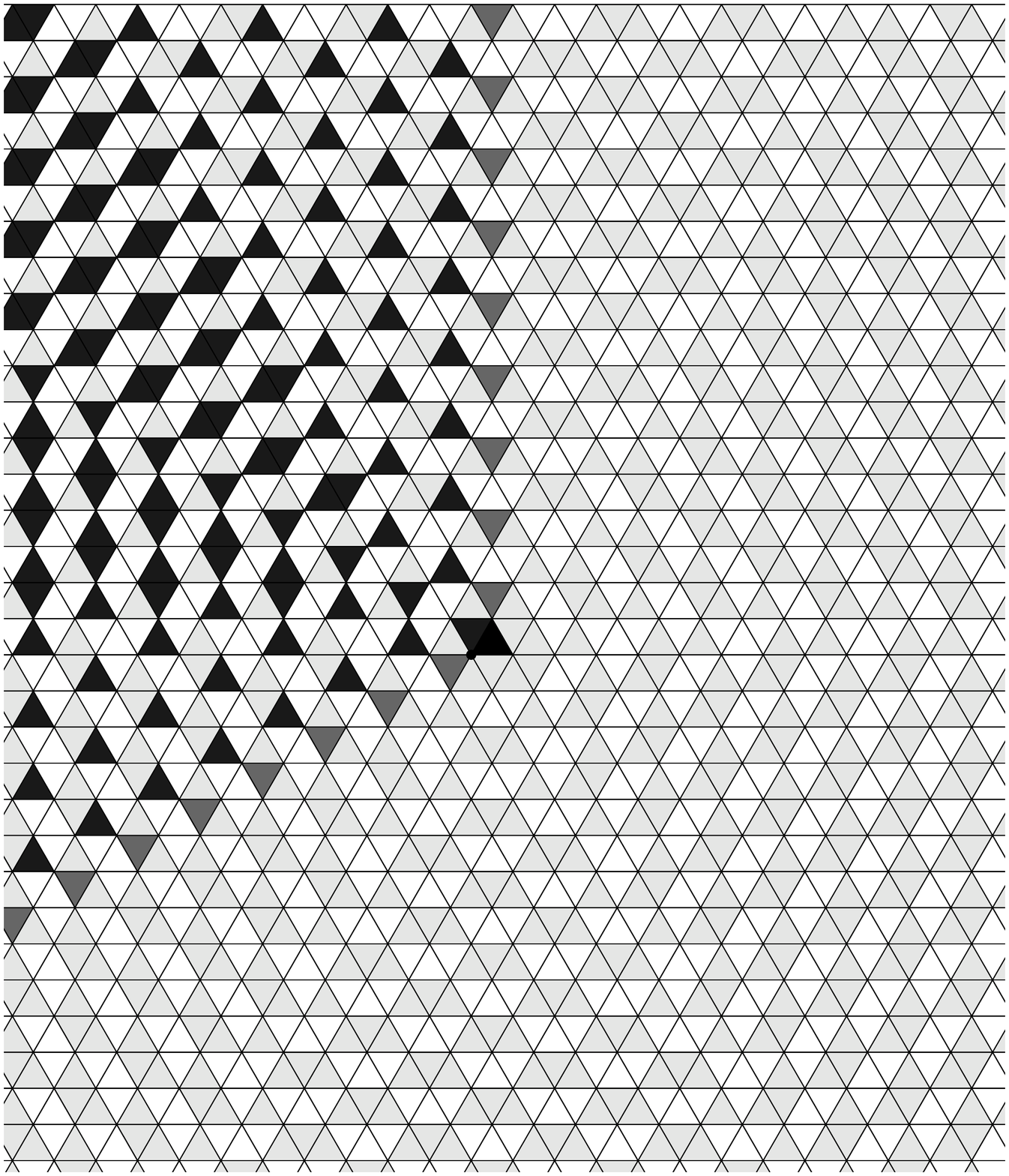} \label{3A.3}
\end{figure}

\begin{figure}
\caption{Partial folding results, type $A_2$, $b=1$, $w=s_1s_2s_1$}
\includegraphics[width=13cm]{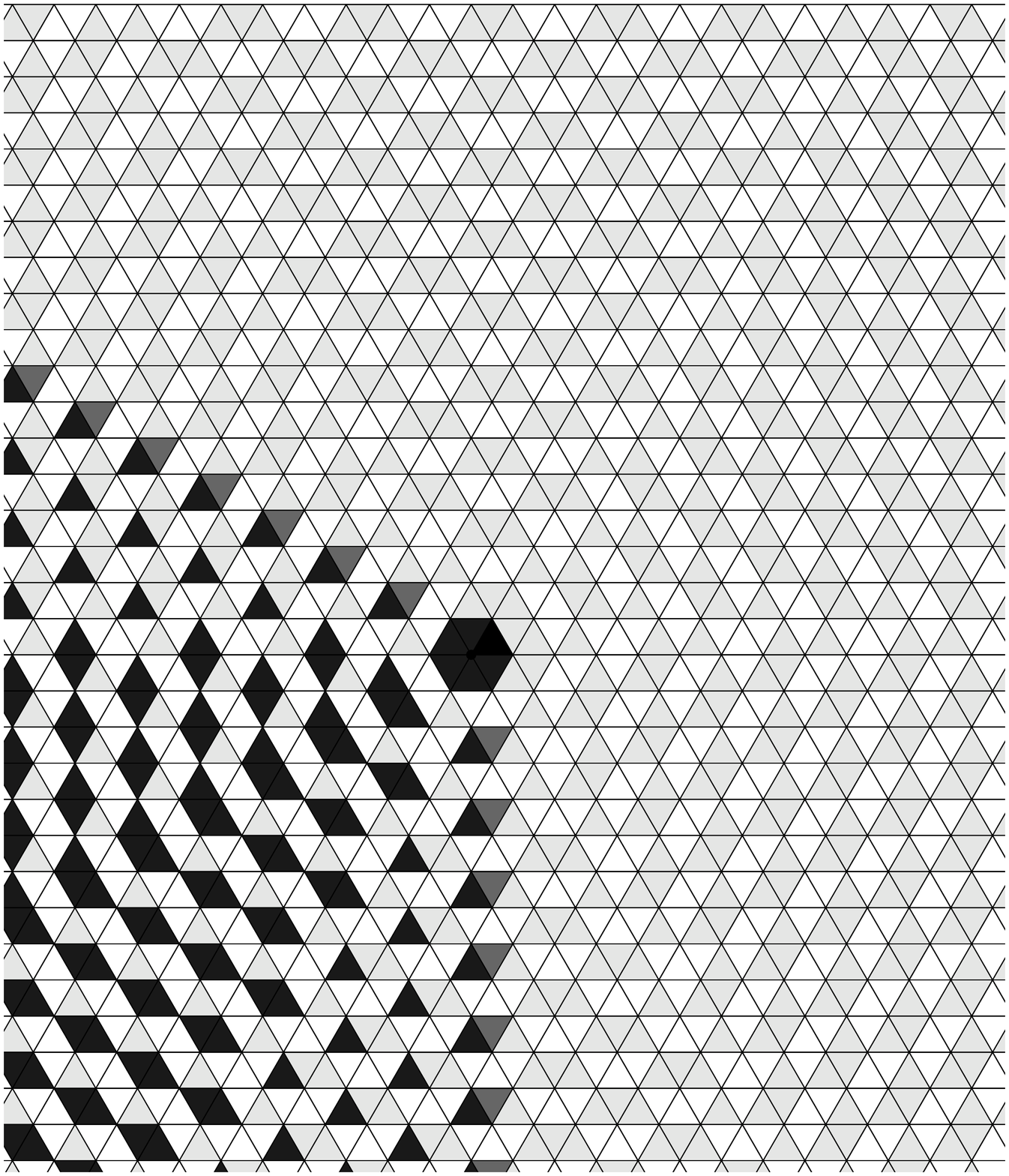} \label{3A.4}
\end{figure}

Another approach to understanding the set of $x$ with non-empty
affine Deligne-Lusztig variety is to consider the folding results
separately for each direction. More precisely, for fixed $w \in
W$, consider the set
\[ \{ x \in \tilde{W};\ d(x, 1, wB^-w^{-1}) \ge 0 \}, \]
where $B^-$ denotes the Borel group opposite to $B$.
We illustrate the partial results in the $A_2$ case in figures \ref{3A.1} ($w =
{\rm id}$), \ref{3A.2} ($w=s_1$), \ref{3A.3} ($w=s_1s_2$) and \ref{3A.4} ($w=s_1s_2s_1$).
In these figures, the alcove corresponding to $x \in \tilde{W}$ is colored

\begin{tabular}{rl}
white & \begin{tabular}{l}if $X_x(1) = \emptyset$ \end{tabular} \\ \hline
light gray & \begin{tabular}{l}if $X_x(1) \ne \emptyset$, but $d(x, 1, wB^-w^{-1}) = -\infty$
\end{tabular} \\ \hline
medium gray & \begin{tabular}{l} if $d(x, 1, wB^-w^{-1}) \ge 0$, but there
exists $w' \in W$ \\ such that $d(x, 1, w'B^-w'^{-1}) > d(x, 1, wB^-w^{-1})$
\end{tabular} \\ \hline
dark gray & \begin{tabular}{l} if $d(x, 1, wB^-w^{-1}) \ge 0$, and for all
$w'\in W$, \\ $d(x, 1, w'B^-w'^{-1}) \le d(x, 1, wB^-w^{-1})$ \end{tabular}
\end{tabular}

\noindent It is unclear whether these separate pieces are easier
to understand than the entire result.

\subsection{The case $b\ne 1$}

\begin{figure}
\caption{Non-empty affine Deligne-Lusztig varieties and their dimensions,
type $A_2$, $b=\epsilon^{(1,0,-1)}$}
\includegraphics[width=13cm]{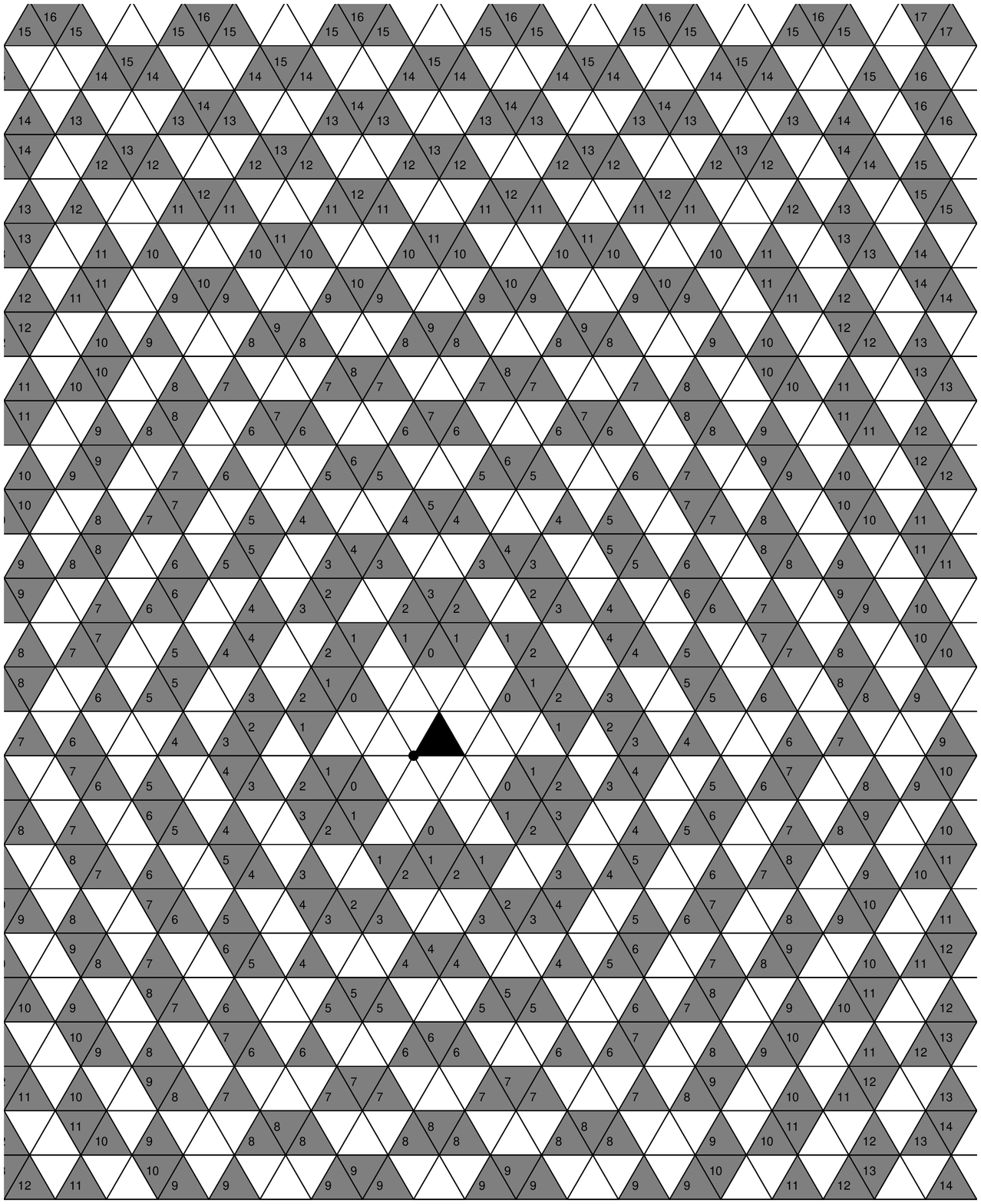} \label{4A}
\end{figure}

\begin{figure}
\caption{Non-empty affine Deligne-Lusztig varieties and their dimensions,
type $C_2$, $b=\epsilon^{(1,0)}$}
\includegraphics[width=13cm]{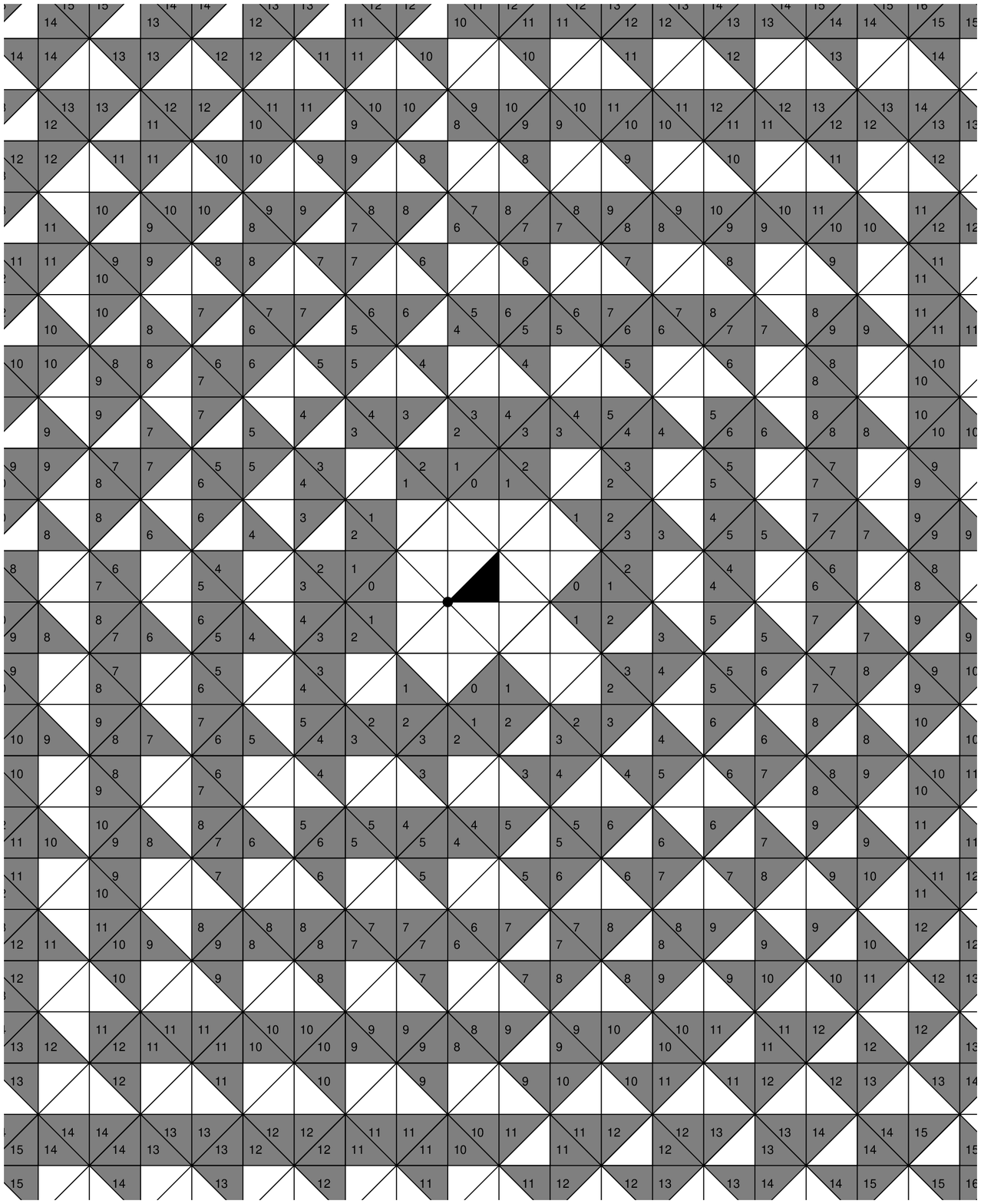} \label{4C}
\end{figure}

We give two examples with $b\ne 1$: figure \ref{4A} shows results
for $A_2$, $b=\epsilon^{(1,0,-1)}$; figure \ref{4C} shows results
for $C_2$, $b=\epsilon^{(1,0)}$.

The following conjecture describes the dimension of
$X_x(\epsilon^\nu)$ for $x$ of sufficient length.
\begin{conjecture}\label{conj.bne1}
Let $G$ be one of the groups $SL_2$, $SL_3$, $Sp_4$. Let $b =
\epsilon^\nu$, and write $\ell(b)$ for the length of the translation  
in $\widetilde{W}$ determined by $\nu$.  
Then there exists $n_0$ such that for all $x \in
\tilde{W}$ with $\ell(x) \ge n_0$, $X_x(b) \ne \emptyset$ if and
only if $X_x(1) \ne \emptyset$. In this case,
\[ \dim X_x(b) = \dim X_x(1)
- \frac{1}{2}\ell(b).
\]
\end{conjecture}
This conjecture describes $X_x(b)$ for $x$ in a larger region of
the standard apartment than Reuman's Conjecture 7.1 in
\cite{Reu04}, but it is less precise because $n_0$ is not
specified. The results in \cite{Reu02} support the new conjecture.
It is possible that the conjecture holds for other groups, too.
Various examples we have examined  seem to indicate that there is
a way to extend Reuman's conjecture to general $b$, but it is
still premature to formulate a precise statement.

\bibliographystyle{amsalpha}

\begin{thebibliography}{EGAIV}


\bibitem[BT72]{bruhat-tits72} F.~Bruhat and J.~Tits, \emph{Groupes
r{\'e}ductifs sur un corps local. {I}},
  Inst. Hautes {\'E}tudes Sci. Publ. Math. \textbf{41} (1972), 5--251.

\bibitem[Ca]{Ca} W. Casselman, \emph{Machine calculations in Weyl groups},
                 Invent. math. \textbf{116} (1994), 95--108.

\bibitem[Da]{Da} R.~Dabrowski, {\em Comparison of the Bruhat and
the Iwahori decompositions of a $\mathfrak p$-adic Chevalley group}, J.
Algebra {\bf 167} (1994), no.3, 704-723.

\bibitem[D]{D} P.~Deligne, {\em La conjecture de Weil II}, Inst.
Hautes \'Etudes Sci., Publ Math. {\bf 52} (1980).


\bibitem[DG]{DG} M.~Demazure and P.~Gabriel, {\em Groupes Alg\'ebriques: Tome
I.   G\'eom\'etrie Alg\'ebrique - G\'en\'eralit\'es -Groupes Commutatifs},
Masson and CIE, Paris (1970), 700 pp. + xxvi.



\bibitem[H]{H} T. J. Haines, \emph{On matrix coefficients of the Satake
isomorphism: complements to the paper of Rapoport}, manuscripta math.
\textbf{101} (2000), 167--174.


\bibitem[HKP]{HKP} T.~Haines, R.~Kottwitz, A.~Prasad, \emph{Iwahori-Hecke
algebras}, math.RT/0309168.


\bibitem[Ka82]{kato} S. Kato, {\em Spherical functions and a $q$-analogue of
Kostant's weight multiplicity formula}, Invent. Math. {\bf 66} (1982), 461-468.



\bibitem[Kot85]{kottwitz85} R.~Kottwitz, \emph{Isocrystals with additional
structure}, Compositio Math.
  \textbf{56} (1985), 201--220.

\bibitem[Kot97]{kottwitz97} R.~Kottwitz, \emph{Isocrystals with
additional structure. {II}}, Compositio
  Math. \textbf{109} (1997), 255--339.

\bibitem[Kot03]{kottwitz03} R.~Kottwitz, \emph{On the {H}odge-{N}ewton
decomposition for split groups}, Int. Math. Res. Not. \textbf{2003},
no. 26, 1433--1447.

\bibitem[Kot05]{kottwitz05} R.~Kottwitz, \emph{Dimensions of {N}ewton strata
in the adjoint quotient of reductive groups}, in preparation.

\bibitem[KR]{kottwitz-rapoport02} R.~Kottwitz and M.~Rapoport, \emph{On the
existence of {F}-crystals},  Comment. Math. Helv. \textbf{78} (2003),
153--184. arXiv:math.NT/0202229.

\bibitem[Lei02]{leigh02} C.~Leigh Lucarelli, \emph{A converse to {M}azur's
inequality for split classical groups},  J.  Inst.  Math. Jussieu,
\textbf{3} (2004), no. 2, 165--183. arXiv:math.NT/0211327.

\bibitem[Lu]{lusztig83} G. Lusztig,
\emph{Singularities, character formulas, and a $q$-analog of weight
multiplicities} In: \emph{Analyse et topologie sur les espaces singuliers,
I-II (Luminy, 1981)}, Soc. Math. France, Paris, 1983, pp. 208--229.

\bibitem[Mat]{Mat} H. Matsumoto, \emph{Analyse Harmonique dans les Syst\`emes
de Tits Bornologiques de Type Affine}, Springer Lecture Notes 590, Berlin,
1977.

\bibitem[Mie]{mier} E. Mierendorff, \emph{Moduli spaces of $p$-divisible
groups}, preprint, February, 2005. arXiv:math.AG/0502320.

\bibitem[MV1]{MV1} I. Mirkovic, K. Vilonen, {\em Perverse sheaves on affine
Grassmannians and Langlands duality}, Math. Res. Lett.  {\bf 7} (2000), no.1,
13-24.

\bibitem[MV2]{MV2} I. Mirkovic, K. Vilonen, {\em Geometric Langlands duality
and representations of algebraic groups over commutative rings}, preprint
2004, math.RT/0401222 v2.


\bibitem[Mum]{Mum} D.~Mumford, {\em The Red Book of Varieties and Schemes},
Lecture Notes in Math. {\bf 1358}, Springer-Verlag (1988).

\bibitem[NP]{NP} B.C. Ng\^{o} and P. Polo, {\em R\'{e}solutions de Demazure
affines et formule de Casselman-Shalika g\'{e}om\'{e}trique}, J. Algebraic
Geom. {\bf 10} (2001), no. 3., 515-547.

\bibitem[Rap00]{Rap00} M. Rapoport, \emph{A positivity property of the Satake
isomorphism}, manuscripta math. \textbf{101} (2000), 153--166.


\bibitem[Rap02]{rapoport02} M.~Rapoport, \emph{A guide to the reduction modulo
$p$ of {S}himura varieties},
  preprint, 2002, arXiv:math.AG/0205022.

\bibitem[RR96]{rapoport-richartz96} M.~Rapoport and M.~Richartz, \emph{On the
classification and specialization of
  {F}-isocrystals with additional structure}, Compositio Math. \textbf{103}
  (1996), 153--181.

\bibitem[RZ96]{rapoport-zink96} M.~Rapoport and T.~Zink, \emph{Period
Spaces for $p$-Divisible Groups}, Ann. of Math. Studies \textbf{141},
  Princeton University Press  (1996).


\bibitem[Reu02]{Reu02} D. Reuman, \emph{Determining whether certain affine
    Deligne-Lusztig sets are non-empty}, Thesis Chicago 2002, math.NT/0211434.

\bibitem[Reu04]{Reu04} D. Reuman, \emph{Formulas for the dimensions of some
    affine Deligne-Lusztig varieties}, Michigan Math. J. \textbf{52} (2004),
    no. 2, 435--451.
\end{thebibliography}
\providecommand{\bysame}{\leavevmode\hbox to3em{\hrulefill}\thinspace}

\end{document}